\documentclass[11pt]{article}

\usepackage{graphicx}
\usepackage{latexsym,amsmath,amsfonts,amscd, amsthm, dsfont}
\usepackage{mathrsfs}
\usepackage{bm,color}
\usepackage{epsfig,verbatim,epstopdf,graphics}
\usepackage{subfigure}
\usepackage{changebar}
\usepackage{multirow}
\usepackage{cases}
\usepackage{makecell}

\usepackage{algorithmic}

\usepackage{yhmath}
 \usepackage{booktabs} 
 \usepackage{tikz}
\usepackage{verbatim}
\usetikzlibrary{arrows,backgrounds,snakes,shapes}
 \numberwithin{equation}{section}

\graphicspath{{./}{./figure/}}
\allowdisplaybreaks

\newtheoremstyle{plainNoItalics}{}{}{\normalfont}{}{\bfseries}{.}{ }{}

\theoremstyle{plain}
\newtheorem{thm}{Theorem}[section]

\theoremstyle{plainNoItalics}

\newtheorem{rem}[thm]{Remark}

\newtheorem{exa}[thm]{Example}

\newcommand{\mdp}{$\mu$DP }

\newcommand{\be}{\begin{eqnarray}}
\newcommand{\ee}{\end{eqnarray}}
\newcommand{\beno}{\begin{eqnarray*}}
\newcommand{\eeno}{\end{eqnarray*}}

\makeatletter

\newcommand{\Rmnum}[1]{\expandafter\@slowromancap\romannumeral #1@}
\makeatother

\begin{document}

\baselineskip=1.8pc


\begin{center}
	{\Large  \bf High order finite difference WENO methods with unequal-sized sub-stencils  for the Degasperis-Procesi type equations}
\end{center}

\vspace{.2in}
\centerline
{
	Jianfang Lin\footnote{School of Mathematical Sciences, Zhejiang University, Hangzhou, Zhejiang 310058, P.R. China. {\tt jianfang.lin@zju.edu.cn}},
		Yan Xu\footnote{School of Mathematical Sciences, University of Science and Technology of China, Hefei, Anhui 230026, P.R. China. {\tt yxu@ustc.edu.cn}},
	Huiwen Xue\footnote{School of Mathematical Sciences, Zhejiang University, Hangzhou, Zhejiang 310058, P.R. China. {\tt xuehuiwen1105@163.com}},
	Xinghui Zhong\footnote{School of Mathematical Sciences, Zhejiang University, Hangzhou, Zhejiang 310058, P.R. China. {\tt zhongxh@zju.edu.cn}}}

\date{\today{}}

\bigskip
\noindent
{\bf Abstract.}
In this paper, we develop two finite difference weighted essentially non-oscillatory (WENO) schemes with unequal-sized sub-stencils for solving the Degasperis-Procesi (DP) and $\mu$-Degasperis-Procesi ($\mu$DP) equations, which contain nonlinear high order derivatives, and possibly peakon solutions or shock waves. By introducing auxiliary variable(s), we rewrite the DP equation as a hyperbolic-elliptic system, and the \mdp equation as a first order system. Then we choose a linear finite difference scheme with suitable order of accuracy for the auxiliary variable(s), and two finite difference WENO schemes with unequal-sized sub-stencils for the primal variable. One WENO scheme uses one large stencil and several smaller stencils, and the other WENO scheme is based on the multi-resolution framework which uses a series of unequal-sized hierarchical central stencils. Comparing with the classical WENO scheme which uses several small stencils of the same size to make up a big stencil, both WENO schemes with unequal-sized sub-stencils are simple in the choice of the stencil and enjoy the freedom of arbitrary positive linear weights. Another advantage is that the final reconstructed polynomial on the target cell  is a polynomial of the same degree as the polynomial over the big stencil, while the classical finite difference WENO reconstruction can only be obtained for specific points inside the target interval. Numerical tests are provided to demonstrate the high order accuracy and non-oscillatory properties of the proposed schemes. 

\vfill


\section{Introduction}
In this paper, we are interested in solving the Degasperis-Procesi (DP) equation
\begin{equation}\label{Eqn:DP}
u_{t}-u_{txx}+4f(u)_{x} = f(u)_{xxx},
\end{equation}
with $x\in \Omega\subset \mathbb{R}$ and
$f(u) = u^{2}/2$,
and the $\mu$-Degasperis-Procesi ($\mu$DP) equation 
\begin{equation}\label{Eqn:mu-DP_ori}
\mu(u)_{t}-u_{txx}+3\mu(u)u_{x}=3u_{x}u_{xx}+uu_{xxx},
\end{equation}
where 
$x\in\mathcal{S}^{1}=\mathbb{R}/\mathbb{Z}$ (the circle whose perimeter equals $1$), and $\mu(u)=\int_{\mathcal{S}^{1}}\! u\,\mathrm{d}x$ denotes the mean of $u$ on $\mathcal{S}^1$. We develop two finite difference weighted essentially non-oscillatory (WENO) schemes for solving \eqref{Eqn:DP} and \eqref{Eqn:mu-DP_ori} with unequal-sized sub-stencils, which provide a simpler way for the reconstruction procedure than the classical WENO schemes,  while still simultaneously maintaining  high order accuracy in smooth regions and controlling spurious numerical oscillations near discontinuities.

The DP equation was singled out first in \cite{degasperis_asymptotic_1999} by an asymptotic integrability test within a family of third order dispersive equations  in the  form of
\begin{equation}
\label{eq:thirdorder}
u_{t}+c_{0}u_{x}+\gamma u_{xxx}-\alpha^{2} u_{txx} = (c_{1}u^{2}+c_{2}u^{2}_{x}+c_{3}uu_{xx})_{x},
\end{equation}
with $\gamma$, $\alpha$, $c_{0}$, $c_{1}$ , $c_{2}$ and $c_{3}$ being real constants. The DP equation \eqref{Eqn:DP} can be transformed from \eqref{eq:thirdorder} with $c_{1} = -\frac{2c_{3}}{\alpha^{2}}$, $c_{2} = c_{3}$, see \cite{Degasperis.Holm.Hone_TMP2002} for more details. It is one of the only three equations that  satisfy the asymptotic integrability condition, besides  the Korteweg-De Vries (KdV) equation ($\alpha = c_{2} = c_{3} = 0$) and the Camassa-Holm (CH) equation ($c_{1} = -\frac{3c_{3}}{2\alpha^{2}}$, $c_{2} = \frac{c_{3}}{2}$).
The DP equation can be regarded as an approximate model of shallow water wave propagation in small amplitude and long wavelength regime \cite{Johnson_JNMP2003, Dullin.Gottwald.Holm_PD2004, Ivanov_PTRS2007, Constantin.Lannes_ARMA2009}, and its asymptotic accuracy is the same as  the CH equation (one order more accurate than the KdV equation). 
The well-posedness of the DP equation has been studied in 
\cite{yin_global_2003,yin_cauchy_2003,yin_global_2004-1,yin_global_2004,Coclite.Karlsen_JFA2006,Coclite.Karlsen_JDE2007,Coclite.Karlsen.Risebro_JNA2008} and the cited references therein.
The $\mu$DP equation is an extensive study of the DP equation. It can be regarded as  an evolution equation on the space of tensor densities over the Lie algebra of smooth vector fields on the circle.

One of the important features of the DP type equations is that they admit not only peakon solutions \cite{Degasperis.Holm.Hone_TMP2002}, 
but also shock waves \cite{Coclite.Karlsen_JFA2006, Lundmark_JNS2007}.  
Explicit expressions of multi-peakon and multi-shock solutions were provided in \cite{lundmark_multi-peakon_2003,lundmark_degasperis-procesi_2005,Lundmark_JNS2007} for the DP equation, and in \cite{Lenells.Misiolek_CMP2010} for the \mdp equation. 
Another feather of the DP type equations is that they satisfy those conservation laws which cannot guarantee the bound of the $H^{1}$-norm of the solution.
Due to these features, it is very difficult to design stable and high order accurate numerical methods for solving the DP and \mdp equations. 
For the DP equation, the existing numerical methods include
the particle method based on the multi-shock peakon solution \cite{Hoel_EJDE2007}, operator splitting finite difference methods \cite{Coclite.Karlsen.Risebro_JNA2008, Feng.Liu_JCP2009}, local discontinuous Galerkin (DG) methods \cite{xu_local_2011},
conservative finite difference methods \cite{Miyatake.Matsuo_JCAM2012},
compact finite difference methods with symplectic implicit Runge-Kutta (RK) time integration \cite{Yu.Sheu_JCP2013}, direct DG  methods \cite{Liu.Huang.Yi_MAA2014}, 
and Fourier spectral methods \cite{Xia_JSC2014, Cai.Sun.Wang_CCP2016}, etc.
Local DG method was developed for the \mdp equation in \cite{zhang_local_2019}.
Recently, classical WENO schemes were investigated for the DP equation in \cite{xia_weighted_2017} and for the \mdp equation in \cite{zhao_high_2020}.

WENO schemes, first designed in \cite{Liu.Osher.Chan_JCP1994} and improved and extended in \cite{jiang_efficient_1996}, were improved version of the essentially non-oscillatory (ENO) schemes \cite{Harten.Osher_JCP1987}. The key idea of ENO and WENO schemes  is actually an approximation procedure used to automatically choose the locally smoothest stencils, aimed at achieving arbitrarily high order accuracy in smooth regions and resolving shocks or other discontinuities sharply and in an essentially non-oscillatory fashion. They have been quite successful for solving hyperbolic and convection-diffusion equations with possibly discontinuous solutions or solutions with sharp gradient regions.  We refer readers to the lecture notes \cite{shu_essentially_1998} and  review papers \cite{Shu_SIAMRev2009, shu_essentially_2020} for the details and development of WENO schemes. 

In this paper, we are interested in solving the DP and $\mu$DP equations using finite difference WENO schemes with  unequal-sized sub-stencils. 
Different from the  classical WENO schemes in \cite{Liu.Osher.Chan_JCP1994, jiang_efficient_1996} which use several small stencils of the same size to make up the big stencil, WENO schemes with  unequal-sized sub-stencils use one big stencil and several smaller stencils, with linear weights chosen to be arbitrarily positive numbers. The idea of this type of WENO procedure first appeared in the context of central WENO schemes in \cite{levy_central_1999,levy_compact_2000,capdeville_central_2008}. Later, the so-called simple WENO  scheme based on this type of WENO reconstruction was constructed for hyperbolic conservation laws in \cite{zhu_new_2016,zhu_new_2017}. More recently, a class of multi-resolution WENO schemes based on this idea was developed for hyperbolic equations in \cite{zhu_new_2018,zhu_new_2019,zhu_new_2020},  in which a hierarchy of nested central stencils is used. The WENO schemes with unequal-sized sub-stencils are particularly attractive  because of their simplicity both in the choice of the stencil and in the freedom of arbitrary positive linear weights, especially for unstructured meshes. 
We refer reader to \cite{abedian_high-order_2013,jiang_high_2021} for WENO schemes with unequal-sized sub-stencils for solving degenerate parabolic equations which involves second order derivatives. 
In this paper, we generalize WENO schemes with unequal-sized sub-stencils to the DP type equations, which involve nonlinear high order derivatives ($>2$). To take care of these nonlinear high order derivatives, especially  dispersion terms $uu_{xxx}$ and $u_{x}u_{xx}$, we introduce auxiliary variable(s) and rewrite the original DP equation as a hyperbolic-elliptic system, and the original \mdp equation as a first order system. Then the primal variable is approximated by a finite difference scheme via the simple finite difference WENO procedure  \cite{zhu_new_2016, zhu_new_2018} or the multi-resolution WENO procedure \cite{zhu_new_2018,zhu_new_2019,zhu_new_2020}, while the auxiliary variable(s) is approximated by a linear finite difference scheme with suitable order of accuracy.  We test the accuracy of our proposed schemes with smooth solutions and non-oscillatory property with various peakon and shock solutions.


The remaining part of the paper is organized as follows: We lay out the details of two finite difference WENO schemes with the unequal-sized sub-stencils  for the DP equation in Section \ref{sec:DP} and for the \mdp equation in Section \ref{sec:mudp}. The numerical performance of the proposed schemes is shown in Section \ref{sec:numerical}, through extensive numerical tests for the DP  and  $\mu$DP equations. Concluding remarks are given in Section \ref{sec:conclusion}. 

\section{Finite difference WENO schemes  for the DP equation}
\label{sec:DP}
In this section, we present the details of the algorithm formulation of two finite difference WENO methods with unequal-sized sub-stencils for solving the DP equation \eqref{Eqn:DP} equipped with suitable initial and boundary conditions. 
We refer readers to \cite{xia_weighted_2017} for the discrete $L^{2}$ stability property of linear finite difference schemes  for the DP equation.

By introducing an auxiliary variable $q$, the DP equation \eqref{Eqn:DP} can be rewritten as a hyperbolic elliptic system
\begin{equation}\label{Eqns:DP}
\begin{cases}
u_{t}+f(u)_{x}+q=0, \\
q-q_{xx}=3f(u)_{x}.	
\end{cases}
\end{equation}

For simplicity, we consider a uniform grid  $\{x_{i}\}_{1, \cdots, N}$ with uniform mesh size $\Delta x = x_{i+1}-x_{i}$. Denote  $x_{i+\frac{1}{2}} =\frac{1}{2} (x_{i}+x_{i+1})$ as the half point. A semi-discrete  finite difference scheme for solving system \eqref{Eqns:DP} is given by
\begin{subequations}
\label{scheme:DP}
\begin{align}
\frac{\mathrm{d}u_{i}(t)}{\mathrm{d} t}+\frac{1}{\Delta x}\left(\hat{f}_{i+\frac{1}{2}}-\hat{f}_{i-\frac{1}{2}}\right)+q_{i}=0, \label{Schema:DP_a} \\
q_{i}-\frac{1}{\Delta x ^2}\left(\hat{q}_{i+\frac{1}{2}}-\hat{q}_{i-\frac{1}{2}}\right)=\frac{3}{\Delta x}\left(\hat{f}_{i+\frac{1}{2}}-\hat{f}_{i-\frac{1}{2}}\right), \label{Schema:DP_b}
\end{align}
\end{subequations}
where  $u_{i}(t)$ and $q_i$ are the numerical approximations to the point values  $u(x_{i}, t)$ and $q(x_i)$, respectively. Here $\hat{f}_{i+\frac{1}{2}}=\hat{f}(u_{i-r},\cdots,u_{i+s})$ and $\hat{q}_{i+\frac{1}{2}}=\hat{q}(u_{i-r},\cdots,u_{i+s})$ are the numerical fluxes which can be obtained by a reconstruction procedure. They are  chosen such that the flux difference approximates the derivatives with high order accuracy, i.e. 
\begin{subequations}
\label{eq:fluxappro}
\begin{align}
\frac{\hat{f}_{i+\frac{1}{2}}-\hat{f}_{i-\frac{1}{2}}}{\Delta x} &= f(u(x))_{x}|_{x_{i}}+\mathcal{O}(\Delta x^{\kappa}), 	\label{eq:fluxappro:f}\\
\frac{\hat{q}_{i+\frac{1}{2}}-\hat{q}_{i-\frac{1}{2}}}{\Delta x^2} &= q_{xx}|_{x_{i}}+\mathcal{O}(\Delta x^{\kappa+1}),	\label{eq:fluxappro:q}
\end{align}
\end{subequations} 
when the solution is smooth, and would generate non-oscillatory solutions when the solution contains possible discontinuities. The collection of grid points involved in the numerical flux, $S=\{x_{i-r},\cdots,x_{i+s}\}$, is called a stencil of the flux approximation.

For the numerical flux $\hat{q}_{i+\frac{1}{2}}$, we use a linear scheme where the flux is a linear combination of point values in the stencil
\begin{align}
\label{eq:linearscheme}
\hat{q}_{i+\frac{1}{2}}=\sum_{j=-r}^s a_j q_{i+j},
\end{align}
with the coefficients $a_j$ chosen to obtain suitable orders of accuracy. For example, a sixth order linear reconstruction
\begin{equation}
\label{6thlinear}
\hat{q}_{i+\frac{1}{2}}=\frac{1}{180}(-2q_{i-2}+25q_{i-1}-245q_{i}+245q_{i+1}-25q_{i+2}+2q_{i+3}),
\end{equation}
corresponds to the numerical flux with the stencil $S=\{x_{i-2},\cdots,x_{i+3}\}$ and gives a truncation error of  $\mathcal{O}(\Delta x^{6})$ in \eqref{eq:fluxappro:q}.
Similarly, a numerical flux with  eighth order  of accuracy in \eqref{eq:fluxappro:q} can be reconstructed as
\begin{equation}
\label{8thlinear}
\hat{q}_{i+\frac{1}{2}}=-\frac{1}{5040}(9q_{i-3}-119q_{i-2}+889q_{i-1}-7175q_{i}+7175q_{i+1}-889q_{i+2}+119q_{i+3}-9q_{i+4}).
\end{equation}

For the numerical flux $\hat{f}_{i+\frac{1}{2}}$, in order to ensure correct upwind biasing and stability, we first split the flux $f(u)$ as
\begin{align}
\label{flux:split}
f(u)=f^{+}(u)+f^{-}(u),
\end{align}
with
\begin{align}
\label{fluxsplitcondition}
\frac{\mathrm{d} f^{+}(u)}{\mathrm{d} u}\geq 0, \quad \frac{\mathrm{d} f^{-}(u)}{\mathrm{d} u}\leq 0.
\end{align}
The most commonly used flux splitting is the  Lax-Friedrichs splitting
\begin{equation}
\label{flux:LF}
f^{+}(u) = f(u)+\alpha u, ~~ f^{-}(u)=f(u)-\alpha u,
\end{equation}
with   $\alpha = \max\limits_{u}|f'(u)|$. More details can be found in the review papers \cite{Shu_SIAMRev2009, shu_essentially_2020}.
Then $f^+_{i+\frac{1}{2}}$ and $f^-_{i+\frac{1}{2}}$ can be reconstructed by a WENO procedure applied 
to $f^+(u)$ and $f^-(u)$ separately.
Finally the numerical flux $\hat{f}_{i+\frac{1}{2}}$ is formed as 
\begin{equation}
\label{hatf}
\hat{f}_{i+\frac{1}{2}}=\hat{f}^{+}_{i+\frac{1}{2}}+\hat{f}^{-}_{i+\frac{1}{2}}.
\end{equation}

Now we show the details of the reconstruction for $\hat{f}^\pm_{i+\frac{1}{2}}$ via the simple WENO procedure in Section \ref{sec:simplewenoDP} and via the multi-resolution WENO procedure in Section \ref{sec:mrwenoDP}. 
\subsection{Simple WENO scheme for $\hat{f}^{\pm}_{i+\frac{1}{2}}$}
\label{sec:simplewenoDP}
In this subsection, we lay out the procedure of the  simple WENO approximation to $\hat{f}^{\pm}_{i+\frac{1}{2}}$ with fifth order of accuracy, following the finite difference simple WENO method constructed in \cite{zhu_new_2016} for solving hyperbolic conservation laws. This class of WENO approximations uses one big stencil and several smaller stencils.   
\begin{figure}[!htbp]
\centering
\begin{tikzpicture}
\draw[line width=1pt] (1.5, 0)--(10.5, 0);

\draw[xshift = 3cm, line width=1pt] (0, 0)--(0, 0.2);
\node[below] at (3, -0.1) {$x_{i-2}$}; 

\draw[xshift = 4.5cm, line width=1pt] (0, 0)--(0, 0.2);
\node[below] at (4.5, -0.1) {$x_{i-1}$}; 

\draw[xshift = 5.25cm, blue, dash dot, line width=1pt] (0, 0)--(0, 0.5);
\node[above, dashed] at (5.25, 0.4) {$x_{i-\frac{1}{2}}$};

\draw[xshift = 6cm, line width=1pt] (0, 0)--(0, 0.2);
\node[below] at (6, -0.1) {$x_{i}$}; 

\draw[xshift = 6.75cm, blue, dash dot, line width=1pt] (0, 0)--(0, 0.5);
\node[fill, circle, scale = 0.3, red] at (6.75, 0) {};
\node[above, dashed] at (6.75, 0.4) {$x_{i+\frac{1}{2}}$}; 

\draw[xshift = 7.5cm, line width=1pt] (0, 0)--(0, 0.2);
\node[below] at (7.5, -0.1) {$x_{i+1}$}; 

\draw[xshift = 9cm, line width=1pt] (0, 0)--(0, 0.2);
\node[below] at (9, -0.1) {$x_{i+2}$}; 

\draw[decorate,decoration={brace, mirror, raise=8pt}, line width=1pt](4.5, -0.4)--(6, -0.4);
\node[below] at (5.25, -0.8) {$S_{2}$};

\draw[decorate,decoration={brace, mirror, raise=8pt}, line width=1pt](6, -1.2)--(7.5, -1.2);
\node[below] at (6.75, -1.6) {$S_{3}$};

\draw[decorate,decoration={brace, mirror, raise=8pt}, line width=1pt](3, -2)--(9, -2);
\node[below] at (6, -2.4) {$S_{1}$};
\end{tikzpicture}
\caption{Stencils used in the fifth order simple WENO scheme.}
\label{Fig:Stencil_WENO-ZQ}
\end{figure}
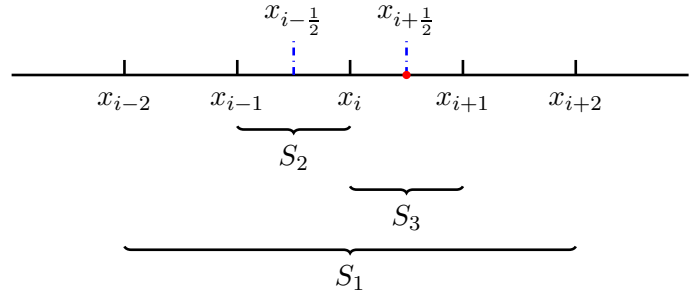

For the numerical flux $\hat{f}^{+}_{i+\frac{1}{2}}$, first choose a big  stencil $S_{1} = \{x_{i-2},x_{i-1},x_i,x_{i+1}, x_{i+2}\}$ and two smaller stencils $S_{2} = \{x_{i-1}, x_{i}\}$, $S_{3}=\{x_{i}, x_{i+1}\}$ as shown in  Figure \ref{Fig:Stencil_WENO-ZQ}, and reconstruct three polynomials: a fourth degree polynomial $p_1 (x)$ satisfying 
\begin{equation}
\label{eq:p1reconstruct}
\frac{1}{\Delta x}\int^{x_{j+\frac{1}{2}}}_{x_{j-\frac{1}{2}}}\! p_{1}(x)\, \mathrm{d} x=f^{+}({u_j}), \quad j=i-2,\ i-1,\ i,\ i+1,\ i+2,
\end{equation}
and two linear polynomials satisfying
\begin{equation}
\label{eq:p2reconstruct}
\frac{1}{\Delta x}\int^{x_{j+\frac{1}{2}}}_{x_{j-\frac{1}{2}}}{p}_{2}(x)\, \mathrm{d} x=f^{+}({u_j}), \quad j=i-1,\ i,
\end{equation}
and
\begin{equation}
\label{eq:p3reconstruct}
\frac{1}{\Delta x}\int^{x_{j+\frac{1}{2}}}_{x_{j-\frac{1}{2}}}{p}_{3}(x)\, \mathrm{d} x=f^{+}({u_j}), \quad j=i,\ i+1.
\end{equation}
The explicit expression of these reconstruction polynomials can be found in \cite{shu_essentially_1998} and thus is omitted here.
Then the WENO approximation is formed based on  the identity 
\begin{align}
\label{eq:p1identity}
p_1(x)=\gamma_1\left(\frac{1}{\gamma_1}p_1(x)
-\frac{\gamma_{2}}{\gamma_{1}}p_{2}(x)-\frac{\gamma_{3}}{\gamma_{1}}p_{3}(x)\right)+\gamma_2p_2(x)+\gamma_3p_3(x),
\end{align}
where   $\gamma_1$, $\gamma_2$ and $\gamma_3$ are three arbitrary positive linear weights. In fact, if we denote 
\begin{align}
\label{p1tilde}
\tilde{p}_{1}(x) &=\frac{1}{\gamma_1}p_1(x)
-\frac{\gamma_{2}}{\gamma_{1}}p_{2}(x)-\frac{\gamma_{3}}{\gamma_{1}}p_{3}(x),
\end{align}
which is also a polynomial of degree four, then the original high order reconstruction $p_1(x)$ on the big stencil $S_1$ can be rewritten as 
\begin{align}
p_1(x)=\gamma_1\tilde{p}_1(x)+\gamma_2p_2(x)+\gamma_3p_3(x),
\end{align}
which can be changed into a WENO reconstruction as
\begin{align}
\label{eq:px}
p(x)=\omega_1\tilde{p}_1(x)+\omega_2p_2(x)+\omega_3p_3(x).
\end{align}
Here $\omega_1,\omega_2$ and $\omega_3$ are nonlinear weights and are computed by the recipe
\cite{zhu_new_2016}:
\begin{equation}
\label{nonlinear}
\omega_{r}=\frac{\bar{\omega}_{r}}{\sum^{3}_{l=1}\bar{w}_{l}},\quad	\bar{\omega}_{r}=\gamma_{r}\left(1+\frac{\tau}{\beta_{r}+\varepsilon}\right),\quad \tau=\left(\frac{|\beta_{1}-\beta_{2}|+|\beta_{1}-\beta_{3}|}{2}\right)^{2},~~r=1, 2, 3, 
\end{equation}
where
$\beta_1,\beta_2$ and $\beta_3$ are the so-called smoothness indicators,  and $\varepsilon$ is a small positive number to avoid the denominator becoming zero. We take $\varepsilon=10^{-10}$ in our numerical experiments. 
For the smoothness indicators $\beta_{r}, r=1, 2, 3$, we use the similar recipe in \cite{jiang_efficient_1996,shu_essentially_1998}, given by
\begin{equation}
\label{eq:simple:beta1}
\beta_{1}=\sum\limits^{4}_{\ell=1}\Delta x^{2\ell-1}\int^{x_{i+\frac{1}{2}}}_{x_{i-\frac{1}{2}}}\!\left(\frac{d^{\ell}\tilde{p}_{1}(x)}{d x^{\ell}}\right)^{2}\mathrm{d}x, \quad	\beta_{r}=\Delta x\int^{x_{i+\frac{1}{2}}}_{x_{i-\frac{1}{2}}}\!\left(\frac{d{p}_{r}(x)}{d x}\right)^{2}\mathrm{d}x, \;\; r = 2, 3, \mathscr{l}
\end{equation}
which are the scaled square sum of the $L^2$-norms of derivatives of the reconstruction polynomials over the interval $I_i=[x_{i-\frac{1}{2}}, x_{i+\frac{1}{2} }]$, and measure how smooth the polynomials $\tilde{p}_1(x), \ p_2(x)$ and $\ p_3(x)$ are in the interval $I_i$, see e.g. \cite{shu_essentially_1998,shu_essentially_2020} for more details.
For the linear weights, any choice of positive numbers which satisfy $\gamma_1+\gamma_2+\gamma_3=1$ is adequate  for accuracy due to the identity \eqref{eq:p1identity}. Considering the balance between accuracy and ability to achieve essentially non-oscillatory shock transitions, we put a larger linear weight for $\gamma_1$ and smaller weights for $\gamma_2,\ \gamma_3$, following the practice in \cite{zhong_simple_2013,zhu_new_2016}. More discussion on this type of the linear weights can be found in \cite{zhong_simple_2013}.

Finally, $\hat{f}^{+}_{i+\frac{1}{2}}$ is given by  $\hat{f}^{+}_{i+\frac{1}{2}}={p}(x_{i+\frac{1}{2}}).$
The numerical flux $\hat{f}^{-}_{i-\frac{1}{2}}$ is obtained by using the above procedure on the same stencils with $f^-(u_{j})$ to obtain $p(x)$, and then by setting $\hat{f}^{-}_{i-\frac{1}{2}}=p(x_{i-\frac{1}{2}})$.

From the reconstruct procedure above, we can claim that beside the linear weights can be chosen as arbitrary positive numbers, another advantage of the simple WENO procedure is that the WENO reconstruction \eqref{eq:px} on the interval $I_i$ is a polynomial of the same degree as the polynomial $p_1(x)$ over the big stencil, while the classical finite difference WENO reconstruction \cite{jiang_efficient_1996} can only be obtained for specific points inside $I_i$.

\begin{rem}
The recipe \eqref{nonlinear} for computing the nonlinear weights $\omega_{1}, \omega_{2}, \omega_{3}$ through the linear weights $\gamma_{1}, \gamma_{2}, \gamma_{3}$ and the smooth indicators $\beta_{1}, \beta_{2}, \beta_{3}$ is different from the recipe used in the classical WENO procedure \cite{jiang_efficient_1996}. This is because $p_2(x)$ and $p_3(x)$ are two first order polynomials and  only second order approximations to $f^+(u)$. The requirement on the closeness of $w_r$ to the linear weights $\gamma_r$ ($\omega_r=\gamma_r+\mathcal{O}(\Delta x^4)$) in smooth regions is more stringent than that for the classical WENO scheme. Thus $\tau$ is introduced in \eqref{nonlinear} associated with the absolute difference of  $\beta_1$ from $\beta_2$ and $\beta_3$. More discussions about this recipe can be found in \cite{borges_improved_2008,castro_high_2011,don_accuracy_2013,zhu_new_2016}.	
\end{rem}
\begin{rem}
$\beta_1$ in \eqref{eq:simple:beta1} is computed using $\tilde{p}_1(x)$ defined in \eqref{p1tilde}, while $\beta_1$ in \cite{zhu_new_2016} is computed using $p_1(x)$ defined in \eqref{eq:p1reconstruct}.
Both choices obtain  high-order accuracy and equally good non-oscillatory results for all of our numerical simulations except for the \mdp equation with shock solutions, in which the numerical solution obtained with $\beta_1$ computed by $\tilde{p}_1(x)$  is essentially non-oscillatory while the numerical solution obtained with   $\beta_1$  computed by $p_1(x)$ has over- and under-shoots.
\end{rem}

\subsection{Multi-resolution WENO scheme for $\hat{f}^{\pm}_{i+\frac{1}{2}}$}
\label{sec:mrwenoDP}

In this subsection, we present the procedure of the  multi-resolution WENO approximation to  $\hat{f}^{\pm}_{i+\frac{1}{2}}$ with  $(2k+1)$-th order of accuracy, following the finite difference multi-resolution WENO method proposed for solving hyperbolic conservation laws in \cite{zhu_new_2018}. This type of WENO reconstructions uses a hierarchy of nested central stencils. 

The multi-resolution WENO scheme to obtain the approximation of $\hat{f}^{+}_{i+\frac{1}{2}}$ with ($2k+1$)-th order accuracy ($k\geq 1$) proceeds as follows:
\begin{figure}[!htbp]
\centering
\begin{tikzpicture}
\draw[line width=1pt] (0, 0)--(12, 0);

\draw[xshift = 1.5cm, line width=1pt] (0, 0)--(0, 0.2);
\node[below] at (1.5, -0.1) {$x_{i-3}$}; 

\draw[xshift = 3cm, line width=1pt] (0, 0)--(0, 0.2);
\node[below] at (3, -0.1) {$x_{i-2}$}; 

\draw[xshift = 4.5cm, line width=1pt] (0, 0)--(0, 0.2);
\node[below] at (4.5, -0.1) {$x_{i-1}$}; 

\draw[xshift = 5.25cm, blue, dash dot, line width=1pt] (0, 0)--(0, 0.5);
\node[above, dashed] at (5.25, 0.4) {$x_{i-\frac{1}{2}}$}; 

\draw[xshift = 6cm, line width=1pt] (0, 0)--(0, 0.2);
\node[below] at (6, -0.1) {$x_{i}$}; 

\draw[xshift = 6.75cm, blue, dash dot, line width=1pt] (0, 0)--(0, 0.5);
\node[fill, circle, scale = 0.3, red] at (6.75, 0) {};
\node[above, dashed] at (6.75, 0.4) {$x_{i+\frac{1}{2}}$}; 

\draw[xshift = 7.5cm, line width=1pt] (0, 0)--(0, 0.2);
\node[below] at (7.5, -0.1) {$x_{i+1}$}; 

\draw[xshift = 9cm, line width=1pt] (0, 0)--(0, 0.2);
\node[below] at (9, -0.1) {$x_{i+2}$}; 

\draw[xshift = 10.5cm, line width=1pt] (0, 0)--(0, 0.2);
\node[below] at (10.5, -0.1) {$x_{i+3}$}; 

\node[fill, circle, scale = 0.3, black, below] at (6, -0.7) {};
\node[below] at (6, -1) {$S_{1}$};

\draw[decorate,decoration={brace, mirror, raise=8pt}, line width=1pt](4.5, -1.4)--(7.5, -1.4);
\node[below] at (6, -1.8) {$S_{2}$};

\draw[decorate,decoration={brace, mirror, raise=8pt}, line width=1pt](3, -2.2)--(9, -2.2);
\node[below] at (6, -2.6) {$S_{3}$};

\draw[decorate,decoration={brace, mirror, raise=8pt}, line width=1pt](1.5, -3)--(10.5, -3);
\node[below] at (6, -3.4) {$S_{4}$};
\end{tikzpicture}
\caption{Stencils used in the multi-resolution WENO scheme.}
\label{Fig:Stencil_MRWENO}
\end{figure}
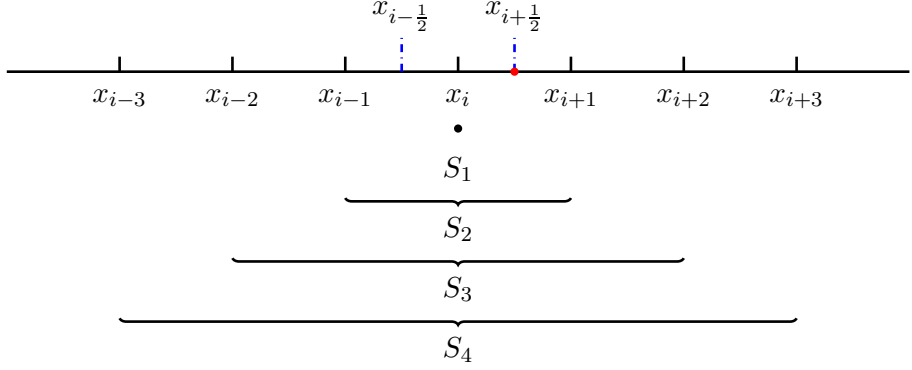
\begin{description}
\item[Step 1.] Choose a series of central stencils 
$		S_{r} = \{x_{i-r+1}, \cdots, x_{i+r-1}\}, r = 1, \cdots, k+1,$
as shown in Figure \ref{Fig:Stencil_MRWENO}. On each stencil $S_{r}$, reconstruct a $(2r-2)$-th degree polynomial $\tilde{p}_{r}(x)$ satisfying
\begin{equation}
\frac{1}{\Delta x}\int^{x_{j+\frac{1}{2}}}_{x_{j-\frac{1}{2}}}\! \tilde{p}_{r}(x)\, \mathrm{d} x=f^{+}({u_{j}}), \quad\quad j=i-r+1,\cdots,i+r-1.
\end{equation}
We explain them in details as follows:
\begin{enumerate}
\item[Step 1.1.] For a third-order approximation, choose two central stencils $S_1=\{x_i\}$ and $S_2=\{x_{i-1},x_i,x_{i+1}\}$, and reconstruct a polynomial $\tilde{p}_1(x)$ of degree zero and a polynomial $\tilde{p}_2(x)$ of degree two satisfying
\begin{align}
\frac{1}{\Delta x}\int^{x_{i+\frac{1}{2}}}_{x_{i-\frac{1}{2}}}\! \tilde{p}_{1}(x)\, \mathrm{d} x&=f^{+}({u_{i}}), \\
\frac{1}{\Delta x}\int^{x_{j+\frac{1}{2}}}_{x_{j-\frac{1}{2}}}\!\tilde{p}_{2}(x)\, \mathrm{d} x&=f^{+}({u_{j}}), \;\; j=i-1,\ i,\ i+1.
\end{align}
\item[Step 1.2.] For a fifth-order approximation, use the central stencil $S_3=\{x_{i-2},x_{i-1},x_i,x_{i+1},x_{i+2}\}$ and reconstruct a polynomial $\tilde{p}_3(x)$ of degree four satisfying
\begin{align}
\frac{1}{\Delta x}\int^{x_{j+\frac{1}{2}}}_{x_{j-\frac{1}{2}}}\! \tilde{p}_{3}(x)\, \mathrm{d} x=f^{+}({u_{j}}), \;\; j=i-2,\ i-1,\ i,\ i+1,\ i+2.
\end{align}
\item[Step 1.3.] For a seventh-order approximation, use the central stencil $S_4=\{x_{i-3}, \cdots,x_{i+3}\}$ and reconstruct a polynomial $\tilde{p}_4(x)$ of degree six satisfying
\begin{align}
\frac{1}{\Delta x}\int^{x_{j+\frac{1}{2}}}_{x_{j-\frac{1}{2}}}\! \tilde{p}_{4}(x)\, \mathrm{d} x=f^{+}({u_{j}}), \;\; j=i-3, i-2, \cdots,\ i+3.
\end{align}
\end{enumerate}
\item[Step 2.] Rewrite these reconstructed polynomials $\tilde{p}_r(x),r=1,\cdots,k+1$ as 
\begin{equation}
\label{eq:pk+1}
\tilde{p}_{r}(x) = \sum\limits^{r}_{l=1}\gamma_{r, l}p_{l}(x),
\end{equation}
similarly as \eqref{eq:px} in the simple WENO procedure discussed in Section \ref{sec:simplewenoDP}, where $\{\gamma_{r, l}\}$ are linear weights, satisfying $\gamma_{r, r}\neq 0,\  \gamma_{r, l}\geq 0$ and $\sum^{r}_{l=1}\gamma_{r, l}=1$.
Then 
\begin{equation}
\label{eq:prx}
p_{r}(x) = \frac{1}{\gamma_{r, r}}\tilde{p}_{r}(x)-\sum\limits^{r-1}_{l=1}\frac{\gamma_{r, l}}{\gamma_{r, r}}p_{l}(x), ~~r=1, \cdots, k+1.
\end{equation}
 We explain them in details as follows:
\begin{enumerate}
\item [Step 2.1.] For a third-order approximation, $p_2(x)$ is defined as
\begin{align}
p_1(x)=\tilde{p}_1(x), \quad	p_2(x)=\frac{1}{\gamma_{2,2}}\tilde{p}_2(x)-\frac{\gamma_{2, 1}}{\gamma_{2, 2}}p_{1}(x),
\end{align}
with $\gamma_{2, 1}+\gamma_{2,2}=1$ and $\gamma_{2,2}\neq 0$.
\item [Step 2.2.] For a fifth-order approximation, $p_3(x)$ is defined as
\begin{align}
p_3(x)=\frac{1}{\gamma_{3,3}}\tilde{p}_3(x)-\frac{\gamma_{3, 1}}{\gamma_{3, 3}}p_{1}(x)-\frac{\gamma_{3, 2}}{\gamma_{3, 3}}p_{2}(x),
\end{align}
with $\gamma_{3, 1}+\gamma_{3, 2}+\gamma_{3,3}=1$ and $\gamma_{3,3}\neq 0$.
\item [Step 2.3.] For a seventh-order approximation, $p_4(x)$ is defined as
\begin{align}
p_4(x)=\frac{1}{\gamma_{4,4}}\tilde{p}_4(x)-\frac{\gamma_{4, 1}}{\gamma_{4, 4}}p_{1}(x)-\frac{\gamma_{4, 2}}{\gamma_{4, 4}}p_{2}(x)-\frac{\gamma_{4, 3}}{\gamma_{4, 4}}p_{3}(x),
\end{align}
with $\gamma_{4, 1}+\gamma_{4, 2}+\gamma_{4, 3}+\gamma_{4,4}=1$ and $\gamma_{4,4}\neq 0$.
\end{enumerate}
We remark that the linear weights can be chosen as arbitrary positive numbers for the sake of accuracy in smooth region. To balance the need of essentially non-oscillatory shock transitions, the linear weights are usually taken in a way as in \cite{zhong_simple_2013,zhu_new_2018}. 
\item[Step 3.] Compute the smoothness indicators $\beta_{r}, r=1, \cdots, k+1$,
by the  recipe \cite{jiang_efficient_1996, shu_essentially_1998,shu_essentially_2020}: 
\begin{equation}\label{Def:smth_mrweno}
\beta_{r}=\sum\limits^{2r-2}_{\ell=1}\Delta x^{2\ell-1}\int^{x_{i+\frac{1}{2}}}
_{x_{i-\frac{1}{2}}}\!\left(\frac{d^{\ell}p_{r}(x)}{dx^{\ell}}\right)^{2}\, \mathrm{d} x,~~r=2, \cdots, k+1,\end{equation}
with slight modifications
made as follows:
\begin{itemize}
\item[Step 3.1.] For $\beta_{1}$, it would equal zero if we use 
the recipe \eqref{Def:smth_mrweno}, since $p_1(x)$ is a constant function. 
This leads to more smeared shock transitions for problems containing strong shocks or contact discontinuities, though it does not cause any problems in the accuracy test for problems with smooth solutions.   
Following \cite{zhu_new_2018}, we increase $\beta_1$ slightly, associated to the smoothness in $\{x_{i-1},x_i\}$ and $\{x_i,x_{i+1}\}$, measured by
$$\pi_{0} =(f^{+}(u_{i})-f^{+}(u_{i-1}))^{2},~\quad\pi_{1}=(f^{+}(u_{i+1})-f^{+}(u_{i}))^{2},$$
with more emphasis on the smaller one of these two measures, i.e. 
\begin{equation}\label{mrweno_beta1}
\begin{split}
\gamma_{1, 0} &=\begin{cases}
1/11, & \pi_{0}\geq \pi_{1}, \\
10/11, & \text{otherwise},
\end{cases} ~~ \gamma_{1, 1} = 1-\gamma_{1, 0}, \\
\theta_{0} &=\gamma_{1, 0}\left(1+\frac{|\pi_{0}-\pi_{1}|^{k}}{\pi_{0}+\varepsilon}\right),~~\theta_{1}=\gamma_{1, 1}\left(1+\frac{|\pi_{0}-\pi_{1}|^{k}}{\pi_{1}+\varepsilon}\right),~~\theta=\theta_{0}+\theta_{1}, \\
\beta_{1}&=\left(\frac{\theta_{0}}{\theta}\left(f^{+}(u_{i})-f^{+}(u_{i-1})\right)+\frac{\theta_{1}}{\theta}(f^{+}(u_{i+1})-f^{+}(u_{i}))\right)^{2}.
\end{split}
\end{equation}
Here $\varepsilon$ is a small positive number to avoid the denominator becoming zero and taken as $\varepsilon=10^{-10}$ in our numerical experiments. We remark here that the modified recipe of $\beta_1$ can be considered that $\beta_1$ in \eqref{mrweno_beta1} is computed by the classical recipe \eqref{Def:smth_mrweno} with a linear function instead of the constant function $p_1(x)$. An example of such a linear function is given by
\begin{align}
\label{eq:modifiedp1}
P_{1}(x) &= \left(\frac{\theta_{0}}{\theta}\left(f^{+}(u_{i})-f^{+}(u_{i-1})\right)+\frac{\theta_{1}}{\theta}(f^{+}(u_{i+1})-f^{+}(u_{i}))\right)\frac{x-x_i}{\Delta x}.
\end{align}

\item[Step 3.2.] For $\beta_{r}, r=2, \cdots, k+1$, we use the  recipe \eqref{Def:smth_mrweno} with $p_r(x)$  in \eqref{eq:prx} modified by replacing the constant function $p_1(x)$ with the linear function $P_1(x)$. That is, if we denote the modified version of $p_r(x)$ as $P_r(x)$, then $P_r(x)$ satisfies
\begin{equation}
\label{eq:Prx}
\begin{split}
P_{r}(x) &= \frac{1}{\gamma_{r, r}}\tilde{p}_{r}(x)-\sum\limits^{r-1}_{l=1}\frac{\gamma_{r, l}}{\gamma_{r, r}}P_{l}(x), ~~r=2, \cdots, k+1,
\end{split}
\end{equation}
with $\tilde{p}_{r}$ defined in \eqref{eq:pk+1} and $P_1(x)$ given in \eqref{eq:modifiedp1}.

\end{itemize}
\item [Step 4.] Compute the nonlinear weights based on the linear weights and smoothness indicators with the recipe \cite{borges_improved_2008,castro_high_2011,don_accuracy_2013,zhu_new_2016}, similar as \eqref{nonlinear}:
\begin{equation}
\omega_{r}=\frac{\bar{\omega}_{r}}{\sum^{k+1}_{l=1}\bar{\omega}_{l}},~~r=1, \cdots, k+1,
\end{equation}
with
\begin{equation}
\bar{\omega}_{r}=\gamma_{k+1, r}\left(1+\frac{\tau_{k+1}}{\beta_{r}+\varepsilon}\right),~~\text{and}~~\tau_{k+1}=\left(\frac{\sum^{k}_{l=1}|\beta_{k+1}-\beta_{l}|}{k}\right)^{k}.
\end{equation}
Here $\varepsilon$ is taken the same as in \eqref{mrweno_beta1}, i.e. $\varepsilon=10^{-10}$ in our numerical experiments.

\item[Step 5.] The final reconstructed polynomial $p(x)$ with $(2k+1)$-th order, is given by
\begin{equation*}
p(x) = \sum\limits^{k+1}_{r=1}\omega_{r}p_{r}(x).
\end{equation*}
\end{description}

The numerical flux $\hat{f}^{+}_{i+\frac{1}{2}}$ is obtained by setting
$\hat{f}^{+}_{i+\frac{1}{2}}={p}(x_{i+\frac{1}{2}}).$
The numerical flux $\hat{f}^{-}_{i-\frac{1}{2}}$ is obtained by using the above procedure on the same stencils with $f^-(u_{j})$ to  obtain $p(x)$, and then by setting $\hat{f}^{-}_{i-\frac{1}{2}}=p(x_{i-\frac{1}{2}})$.

In summary, to build a finite difference WENO scheme for solving \eqref{Eqns:DP} with unequal-sized sub-stencils, given the point values $\{u_i\}$, we proceed as follows:\\[2mm]
\underline{\bf Procedure I. Finite difference WENO scheme for the DP equation}
\begin{enumerate}
\item Find a smooth flux splitting \eqref{flux:split}, satisfying \eqref{fluxsplitcondition}.	
\item Compute $f^\pm(u_{j})$ and follow the simple WENO scheme in Section \ref{sec:simplewenoDP} or the multi-resolution WENO scheme in Section \ref{sec:mrwenoDP}, to obtain $\hat{f}^{\pm}_{i+\frac{1}{2}}$ for all $i$. Form $\hat{f}_{i+\frac{1}{2}}$  by \eqref{hatf} for all $i$.
\item Choose a linear scheme \eqref{eq:linearscheme} to compute the flux $\hat{q}_{i+\frac{1}{2}}$ for all $i$, satisfying \eqref{eq:fluxappro:q}. In our numerical tests, \eqref{6thlinear} is chosen when  $\hat{f}_{i+\frac{1}{2}}$ is obtained by the fifth-order simple WENO or the fifth-order multi-resolution scheme, and \eqref{8thlinear} is chosen when  $\hat{f}_{i+\frac{1}{2}}$  is obtained by the seventh-order multi-resolution scheme.
\item Form the scheme \eqref{scheme:DP}.   If we  denote  $ {\bf p}=(p_1,\cdots,p_N)^T$ and 
${\bf \tilde{f}}=(\tilde{f}_1,\cdots,\tilde{f}_N)^T$ with  $\tilde{f}_i=\frac{3}{\Delta x}(\hat{f}_{i+\frac{1}{2}}-\hat{f}_{i-\frac{1}{2}})$,  then \eqref{Schema:DP_b} can be written in the following matrix form
\begin{align}
{\bf A}{\bf p} = {\bf \tilde{f}}.
\end{align} 
Apply a linear solver with the matrix $\bf A$ and we get ${\bf p}={\bf A}^{-1}{\bf \tilde{f}}$, which can be used in \eqref{Schema:DP_a}, yielding the semi-discrete discretization
\begin{equation}
\label{semi:vector}
\frac{\mathrm d {u}_{i}}{\mathrm dt} ={L}( {u})_i.
\end{equation}
\item Apply any standard ODE solver for the time discretization of \eqref{semi:vector}, e.g. the third-order strong-stability-preserving (SSP) Runge-Kutta (RK) method \cite{shu_efficient_1988}:
\begin{align}
\label{tvdrk}
\begin{split}
&u^{(1)} = u^{n}+\Delta tL(u^{n}),\\ 
&u^{(2)} =\frac{3}{4}u^{n}+\frac{1}{4}\left(u^{(1)}+\Delta t L(u^{(1)})\right), \\
&u^{n+1} = \frac{1}{3}u^{n}+\frac{2}{3}\left(u^{(2)}+\Delta t L(u^{(2)})\right).
\end{split}
\end{align} 
\end{enumerate}

\section{Finite  difference WENO schemes for the \mdp equation}
\label{sec:mudp}
In this section, we present the algorithm formulation of two finite difference WENO methods with unequal-sized sub-stencils for solving the \mdp equation \eqref{Eqn:mu-DP_ori} equipped with suitable initial conditions. 

We consider an equivalent form of the $\mu$DP equation \eqref{Eqn:mu-DP_ori}, given by
\begin{equation}\label{Eqn:mu-DP}
u_{t}+uu_{x}+3\mu(u)(A^{-1}_{\mu}u)_{x} = 0,
\end{equation}
where $A_{\mu}(u)=\mu(u)-u_{xx}$ is an invertible operator. More details can be found in 
\cite{zhang_local_2019}.
By introducing auxiliary variables $q$ and $v$, Equation \eqref{Eqn:mu-DP} can be rewritten as a first-order system\begin{subequations}
\label{Eq:mu-DP}
\begin{align}
u_{t}+f(u)_{x}+3\mu(u)q &=0, \\
q-v_{x} &= 0, \\
\mu(v)-q_{x} &= u,
\end{align}	
\end{subequations}
with $f(u) = u^{2}/2$.   
For simplicity, we just take one period of the whole domain $\mathbb{R}$, i.e. [0,1], to represent 
the circle $\mathcal{S}^1$.
We once again consider a uniform grid: 
\begin{align}
\label{grid}
0=x_1<x_2<\cdots<x_{N+1}=1.
\end{align}
Denote $x_{i+\frac{1}{2}} =\frac{1}{2} (x_{i}+x_{i+1})$ as the half point, and $\Delta x = x_{i+1}-x_{i}$ as the mesh size. Then, a semi-discrete finite difference scheme for solving the system \eqref{Eq:mu-DP} is given by
\begin{subequations}
\label{scheme:mudp}
\begin{align}
\frac{\mathrm{d}u_{i}(t)}{\mathrm{d}t}+\frac{1}{\Delta x}\left(\hat{f}_{i+\frac{1}{2}}-\hat{f}_{i-\frac{1}{2}}\right)
+3\left(\Delta x\sum\limits^{N}_{j=1} u_{j}\right)q_{i} = 0,  \label{Schema:mu-DP_a}\\ 
q_{i}-\frac{1}{\Delta x}\left(\hat{v}_{i+\frac{1}{2}}-\hat{v}_{i-\frac{1}{2}}\right) = 0,  \label{Schema:mu-DP_b}\\
\left(\Delta x\sum\limits^{N}_{j=1}v_{j}\right)-\frac{1}{\Delta x}\left(\hat{q}_{i+\frac{1}{2}}-\hat{q}_{i-\frac{1}{2}}\right)=u_{i}, \label{Schema:mu-DP_c}
\end{align}
\end{subequations}
where $u_{i}(t), \ q_i$ and $\ v_i $ are the numerical approximations to the point values  $u(x_{i}, t), \ q(x_i)$ and $\ v(x_i)$, respectively. $\Delta x\sum\limits^{N}_{j=1} u_{j}$ and $\Delta x\sum\limits^{N}_{j=1}v_{j}$ are the discrete form of $\mu(u)=\int_{0}^1 u\mathrm{d}x$ and  $\mu(v)=\int_{0}^1 v\mathrm{d}x$ with the grid \eqref{grid}, respectively.
The numerical fluxes $\hat{f}_{i+\frac{1}{2}}$, $\hat{q}_{i+\frac{1}{2}}$ and $\hat{v}_{i+\frac{1}{2}}$ are chosen such that the flux difference approximates the derivatives with high order accuracy, i.e. 
\begin{subequations}
\label{eq:fluxappro:mu}
\begin{align}
&\frac{\hat{f}_{i+\frac{1}{2}}-\hat{f}_{i-\frac{1}{2}}}{\Delta x} = f(u(x))_{x}|_{x_{i}}+\mathcal{O}(\Delta x^{\kappa}), 	\label{eq:fluxappromu:f}\\
&\frac{\hat{v}_{i+\frac{1}{2}}-\hat{v}_{i-\frac{1}{2}}}{\Delta x} = v_{x}|_{x_{i}}+\mathcal{O}(\Delta x^{\kappa}),	\quad	\frac{\hat{q}_{i+\frac{1}{2}}-\hat{q}_{i-\frac{1}{2}}}{\Delta x} = q_{x}|_{x_{i}}+\mathcal{O}(\Delta x^{\kappa}),	\label{eq:fluxappromu:q}
\end{align}
\end{subequations} 
when the solution is smooth, and would generate non-oscillatory solutions when the solution contains possible discontinuities.

For the numerical fluxes $\hat{q}_{i+\frac{1}{2}}$ and $\hat{v}_{i+\frac{1}{2}}$, we take the simple choice given by
\begin{equation}
\label{alternating}
\hat{v}_{i+\frac{1}{2}}=v^{-}_{i+\frac{1}{2}}, ~~\quad \hat{q}_{i+\frac{1}{2}}=q^{+}_{i+\frac{1}{2}},
\end{equation}
where $v^{-}_{i+\frac{1}{2}}$ and $q^{+}_{i+\frac{1}{2}}$ are reconstructed by a linear scheme  \eqref{eq:linearscheme} with the coefficients chosen to obtain suitable order of accuracy.
For example, a fifth order linear reconstruction is given by
\begin{subequations}
\label{5thlinear}
\begin{align}
v^{-}_{i+\frac{1}{2}} &= \frac{2}{60}v_{i-2}-\frac{13}{60}v_{i-1}+\frac{47}{60}v_{i}+\frac{27}{60}v_{i+1}-\frac{3}{60}v_{i+2}, \\
q^{+}_{i+\frac{1}{2}} &=-\frac{3}{60}q_{i-1}+\frac{27}{60}q_{i}+\frac{47}{60}q_{i+1}-\frac{13}{60}q_{i+2}+\frac{2}{60}q_{i+3},
\end{align}
\end{subequations}
and a seventh order linear reconstruction is given by
\begin{subequations}
\label{7thlinear}
\begin{align}
v^{-}_{i+\frac{1}{2}} &= -\frac{3}{420}v_{i-3}+\frac{25}{420}v_{i-2}-\frac{101}{420}v_{i-1}+\frac{319}{420}v_{i}+\frac{214}{420}v_{i+1}-\frac{38}{420}v_{i+2}+\frac{4}{420}v_{i+3}, \\
q^{+}_{i+\frac{1}{2}} &= \frac{4}{420}q_{i-2}-\frac{38}{420}q_{i-1}+\frac{214}{420}q_{i}+\frac{319}{420}q_{i+1}-\frac{101}{420}q_{i+2}+\frac{25}{420}q_{i+3}-\frac{3}{420}q_{i+4}. 
\end{align}
\end{subequations}

For the numerical flux $\hat{f}_{i+\frac{1}{2}}$, the reconstructed procedure is the same as $\hat{f}_{i+\frac{1}{2}}$ in the  the DP equation discussed in Section \ref{sec:simplewenoDP} and \ref{sec:mrwenoDP}.

Now we summarize the procedure with a finite difference WENO scheme for solving \eqref{Eq:mu-DP} with unequal-sized sub-stencils, given the point values $\{u_i\}$, as follows:\\[2mm]
\underline{\bf Procedure II. Finite difference WENO scheme for the \mdp equation}
\begin{enumerate}
\item Apply steps 1-2 in Procedure I as we do for the DP equation to reconstruct$\hat{f}_{i+\frac{1}{2}}$ for all $i$.
\item Choose a linear scheme \eqref{eq:linearscheme} to compute the fluxes $\hat{v}_{i+\frac{1}{2}}$ and $\hat{q}_{i+\frac{1}{2}}$ defined in \eqref{alternating} for all $i$, satisfying \eqref{eq:fluxappromu:q}. In our numerical tests, \eqref{5thlinear} is chosen when $\hat{f}_{i+\frac{1}{2}}$  is obtained by the fifth-order simple WENO or the fifth-order multi-resolution scheme, and \eqref{7thlinear} is chosen when $\hat{f}_{i+\frac{1}{2}}$ is obtained by the seventh-order multi-resolution scheme.
\item Form the scheme \eqref{scheme:mudp}. If we denote ${\bf u}=(u_1,\cdots,u_N)^T,{\bf v}=(v_1,\cdots,v_N)^T$ and $ {\bf p}=(p_1,\cdots,p_N)^T$,  then \eqref{Schema:mu-DP_b} and \eqref{Schema:mu-DP_c} can be written in the following matrix form
\begin{align}
{\bf q}-{\bf A}{\bf v} &={\bf 0}, \label{Afc:a} \\
{\bf B}{\bf v}-{\bf C}{\bf q} &= {\bf u}.
\end{align}
Let ${\bf D} = {\bf B}-{\bf C}{\bf A}$. We can get ${\bf v}={\bf D}^{-1} {\bf u}$  by applying a linear solver with the matrix $\bf D$.
Then \eqref{Schema:mu-DP_a} with $\bf q$ now being expressed by $\bf u$ via \eqref{Afc:a}  yields 
the following semi-discrete discretization
\begin{equation}
\label{semi:vector2}
{\bf u}_{t}={\bf L}({\bf u}).
\end{equation}
\item Apply e.g. \eqref{tvdrk} for the time discretization of \eqref{semi:vector2}.
\end{enumerate}

%
%

\section{Numerical results}
\label{sec:numerical}
In this section, we present numerical tests to demonstrate the performance of the fifth order finite difference simple WENO scheme, the fifth order and seventh order finite difference multi-resolution WENO schemes, denoted as WENO5, MR-WENO5 and MR-WENO7, respectively for simplicity,  for the DP  and  $\mu$DP equations. Temporal discretization is carried out by the SSP RK method \eqref{tvdrk}, unless otherwise specified. The time step is set as
$\Delta t = \text{CFL}\cdot\Delta x$, where CFL is taken as  0.3.
 For the accuracy tests (Examples \ref{Ex:twave_sm} and \ref{Ex:mu-DP_twave_sm}), we adjust the time step $\Delta t$ as  $\Delta t = \text{CFL}\cdot\Delta x^{5/3}$ for the fifth-order WENO schemes, and $\Delta t = \text{CFL}\cdot\Delta x^{7/3}$ for the seventh order WENO scheme.   For examples with peakon or shock solutions (Example \ref{Ex:peakon_one}-\ref{Ex:twopeakon}, Example \ref{Ex:peakon+anti-peakon}-\ref{Ex: wavebreak}, Example \ref{Ex:mudp_peakon}), all three methods obtain equally good non-oscillatory solutions and thus we only show the numerical results obtained by  one method (randomly chosen) to save space.  We set linear weights as
\begin{equation*}
	\gamma_{1} = 0.98, ~~\gamma_{2} = 0.01, ~~ \gamma_{3} = 0.01,
\end{equation*}
for WENO5, and
\begin{equation*}
	\begin{array}{llll}
		\gamma_{2, 1} = 1/11, & \gamma_{2, 2} = 10/11, & & \\
		\gamma_{3, 1} = 1/111, & \gamma_{3, 2} = 10/111, & \gamma_{3, 3} = 100/111, &  \\
		\gamma_{4, 1} = 1/1111, & \gamma_{4, 2} = 10/1111, &\gamma_{4, 3} = 100/1111, &\gamma_{4, 4} = 1000/1111,
	\end{array}
\end{equation*}
for MR-WENO5 and MR-WENO7, unless otherwise specified. 


\subsection{Numerical results of the DP equation}
In this section, we present the numerical results to demonstrate the performance of WENO5, MR-WENO5 and MR-WENO7 for the DP equation with different initial conditions. The computation domain is chosen large enough such that the solution is small enough at the boundary of the domain for periodic boundary conditions to hold approximately at the level of truncation errors. 

\begin{exa}\label{Ex:twave_sm}
{\bf Accuracy test for the single smooth soliton solution}\\
In this example,	we consider the DP equation with the traveling wave solution $u(x , t) = U(x-ct)$, where $c$ is the wave speed. Denote $\xi = x-ct$, and assume  $\lim\limits_{\xi\rightarrow \infty} U(\xi)= A$. The smooth soliton solution of the DP equation can be given in an explicit formula \cite{Zhang.Qiao_MPAG2007} as 
\begin{equation*}
U(\xi) = A((4-\sqrt{5})-\frac{2\sqrt{5}}{X(\xi)^{2}-1}),
\end{equation*}
where $X(\xi)$ is defined by
\begin{align*}
X(\xi) &=\left(-\frac{7+3\sqrt{5}}{3}b+\frac{38+17\sqrt{5}}{27}b^{3}+\sqrt{\frac{2+\sqrt{5}}{27}+\frac{517+231\sqrt{5}}{54}b^{2}-\frac{521+233\sqrt{5}}{54}b^{4}}\right)^{\frac{1}{3}} \\
&+\left(-\frac{7+3\sqrt{5}}{3}b+\frac{38+17\sqrt{5}}{27}b^{3}-\sqrt{\frac{2+\sqrt{5}}{27}+\frac{517+231\sqrt{5}}{54}b^{2}-\frac{521+233\sqrt{5}}{54}b^{4}}\right)^{\frac{1}{3}} \\
&+\frac{2+\sqrt{5}}{3}b,
\end{align*}
with $b = \frac{1+e^{-|\xi|}}{1-e^{-|\xi|}}$.
We set $A = 1$ and $c = 5$, and take the computational domain as $[-50, 50]$. 
We list the $L^{1}$ and $L^{\infty}$ errors and  orders of accuracy with WENO5, MR-WENO5 and MR-WENO7 at $T=1$  in Table \ref{Tab:twave_sm}. We can see that all three methods achieve the desired order of accuracy, i.e. fifth-order accuracy for WENO5 and MR-WENO5 and seventh-order accuracy for MR-WENO7.
	
\begin{table}[!htbp]
\centering
\caption{The DP equation with the single smooth soliton solution in Example \ref{Ex:twave_sm} at $T=1$.}
\label{Tab:twave_sm}
\smallskip
\begin{tabular}{ccccccccc}
\toprule
& \multicolumn{4}{c}{WENO5} & \multicolumn{4}{c}{MR-WENO5}  \\ \cmidrule(lr){2-5}\cmidrule(lr){6-9}
$N$ & $L^{1}$ \mbox{error} & \mbox{Order} & $L^{\infty}$ \mbox{error} & \mbox{Order} & $L^{1}$ \mbox{error} & \mbox{Order} & $L^{\infty}$ \mbox{error} & \mbox{Order}\\ \hline
80 &      1.65E-02 &          &     9.27E-02 &          &     1.58E-02 &          &     8.92E-02 &         \\
160 &     7.30E-04 &     4.45 &     5.33E-03 &     4.08 &     7.20E-04 &     4.41 &     5.13E-03 &     4.08\\
320 &     2.23E-05 &     5.01 &     1.64E-04 &     5.00 &     2.37E-05 &     4.89 &     1.65E-04 &     4.94\\
640 &     6.24E-07 &     5.15 &     5.01E-06 &     5.02 &     6.15E-07 &     5.27 &     4.18E-06 &     5.29\\
1280 &    1.96E-08 &     4.99 &     1.26E-07 &     5.30 &     1.96E-08 &     4.97 &     1.26E-07 &     5.05\\
\midrule
& \multicolumn{8}{c}{MR-WENO7}  \\ \cmidrule(lr){2-9}
$N$ & &$L^{1}$ \mbox{error} && \mbox{Order} & &$L^{\infty}$ \mbox{error} && \mbox{Order} \\ \hline
80 &      &2.67E-02 &     &     &     &1.26E-01 &    &     \\
160 &     & 2.32E-04&     &6.79 &     &1.46E-03 &    &6.37\\
320 &     &1.78E-06 &     &6.99 &     &1.15E-05 &    &6.96\\
640 &     &1.43E-08 &     &6.94 &     &8.91E-08 &    &7.00\\
\bottomrule    
\end{tabular}
\end{table}

\end{exa}

\begin{exa}
\label{Ex:peakon_one}
{\bf Single peakon and anti-peakon solutions}\\
In this example, we consider  the wave propagation of the peakon solution \cite{Degasperis.Holm.Hone_TMP2002} given by
$$u(x, t) = ce^{-|x-ct|},$$ 
and the anti-peakon solution given by $$u(x, t) = -ce^{-|x+ct|}.$$ 
\begin{figure}[!htbp]
	\centering%
	\includegraphics[width=0.45\textwidth]{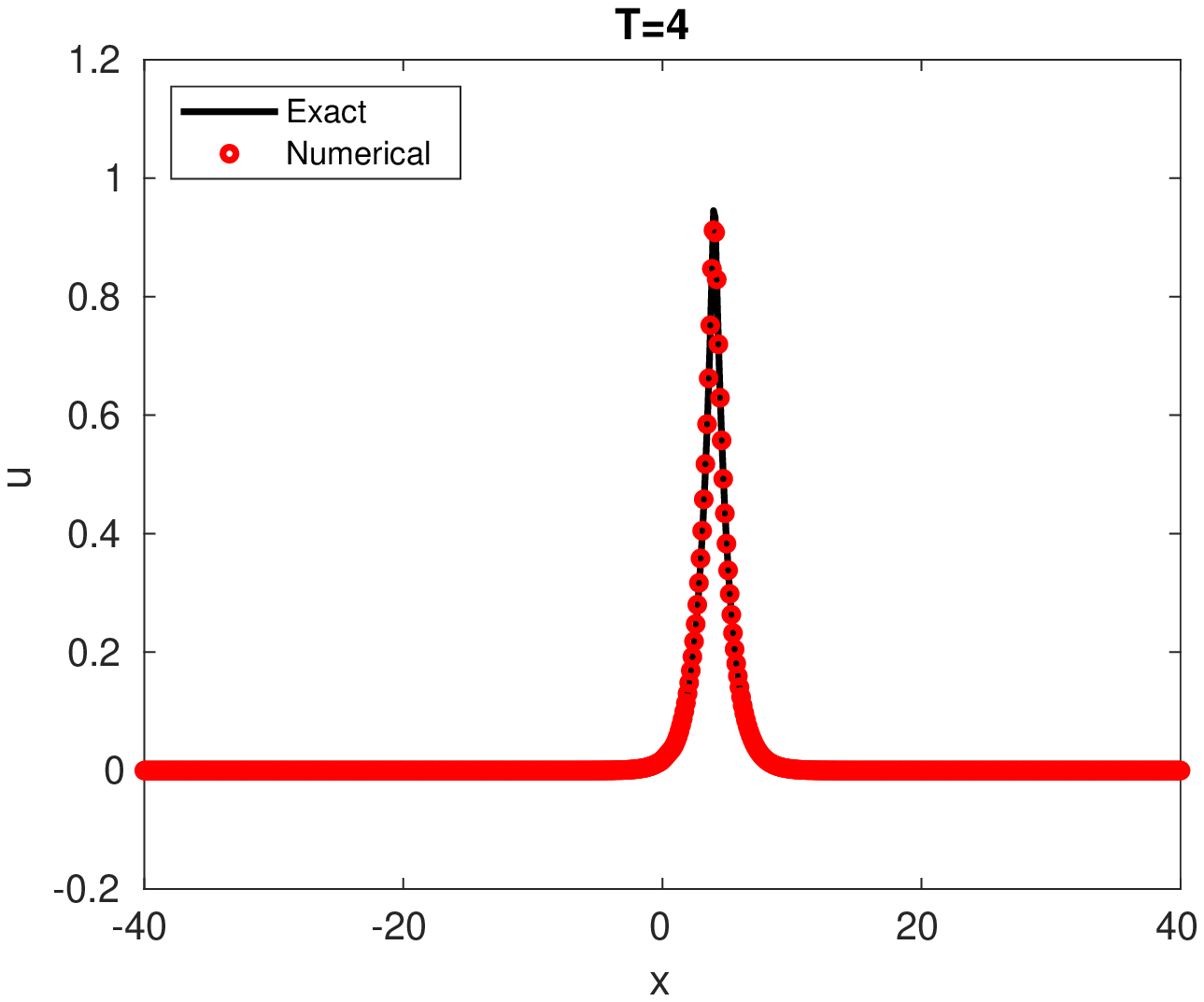}	
	\includegraphics[width=0.45\textwidth]{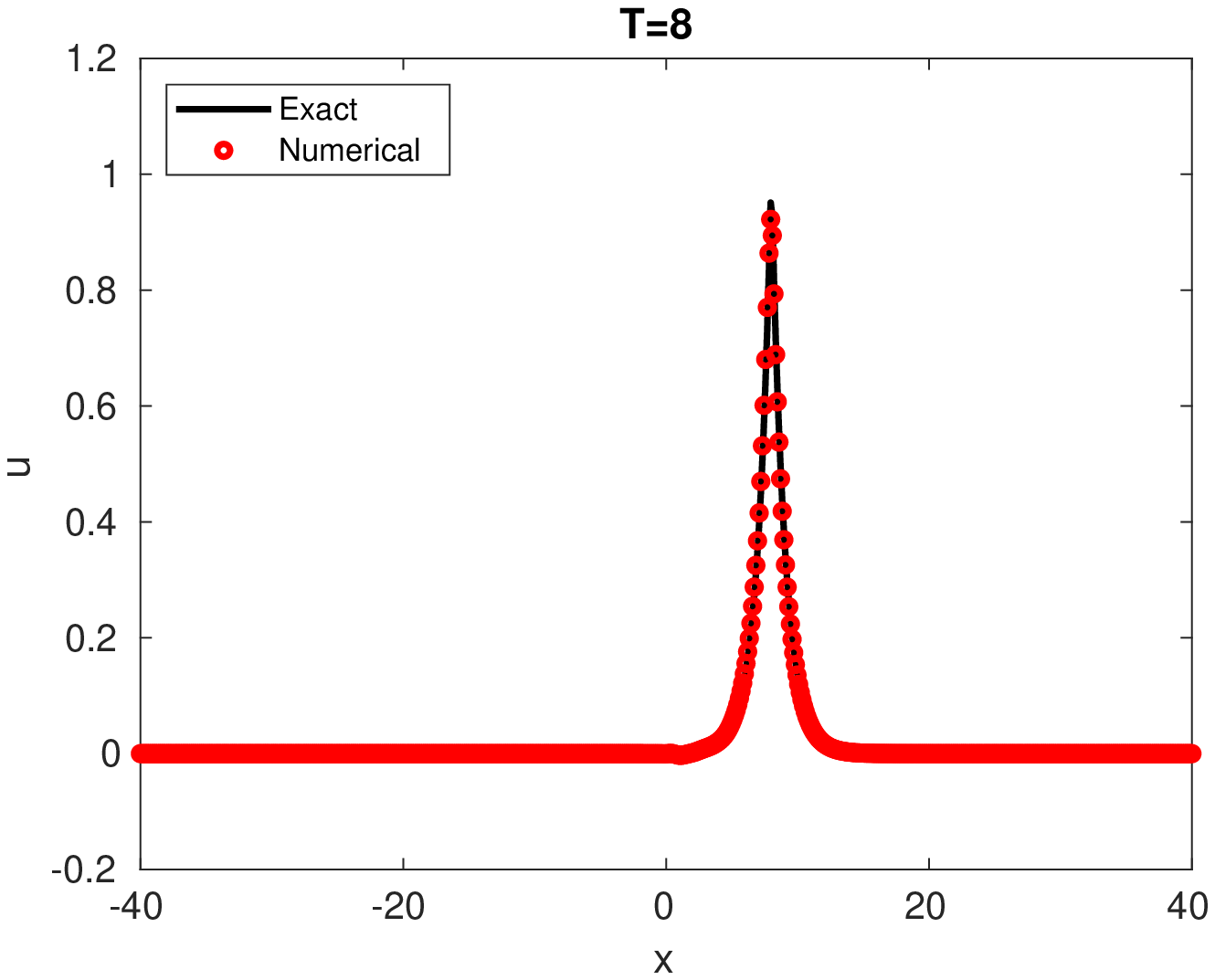}
	\includegraphics[width=0.45\textwidth]{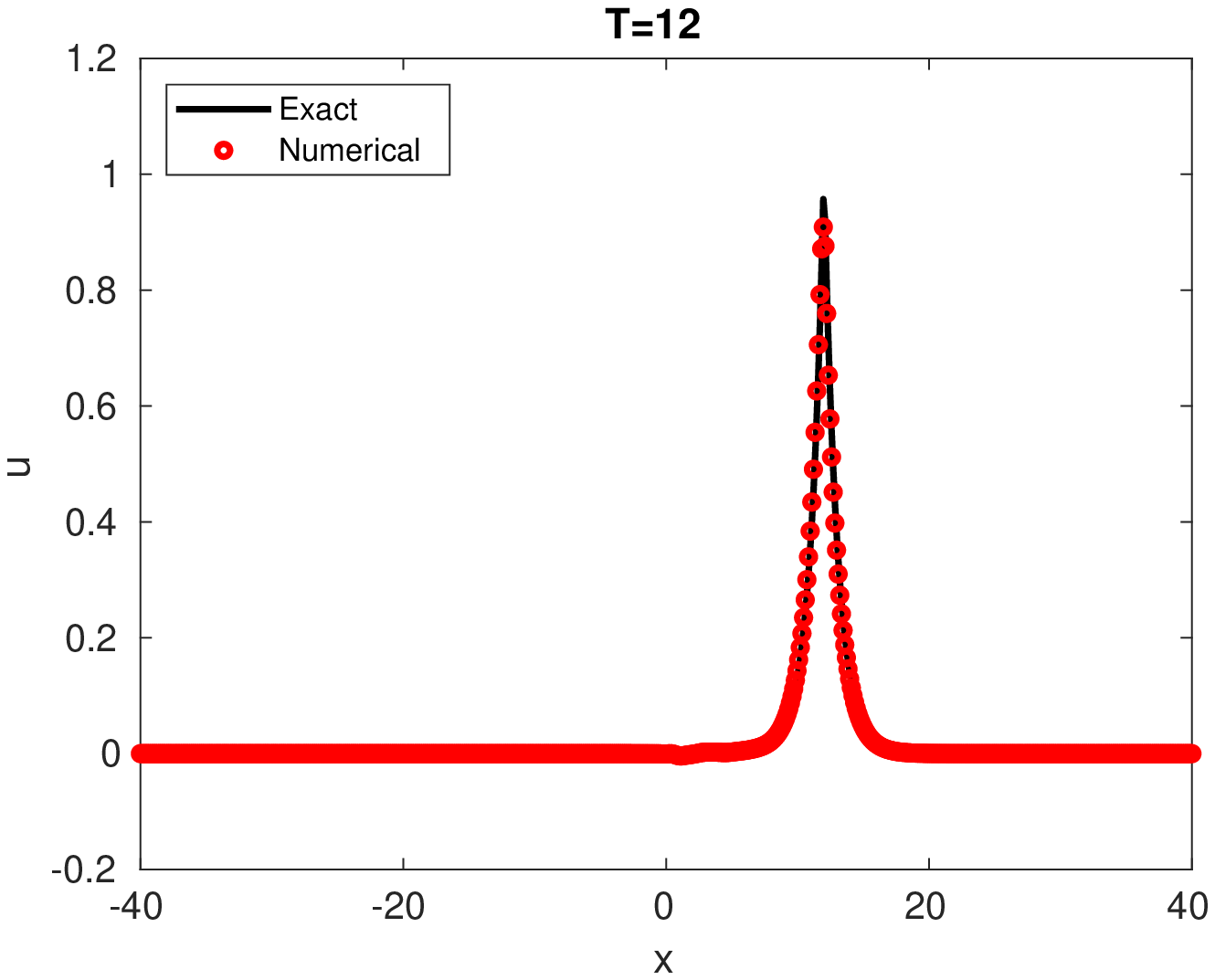} 
	\includegraphics[width=0.45\textwidth]{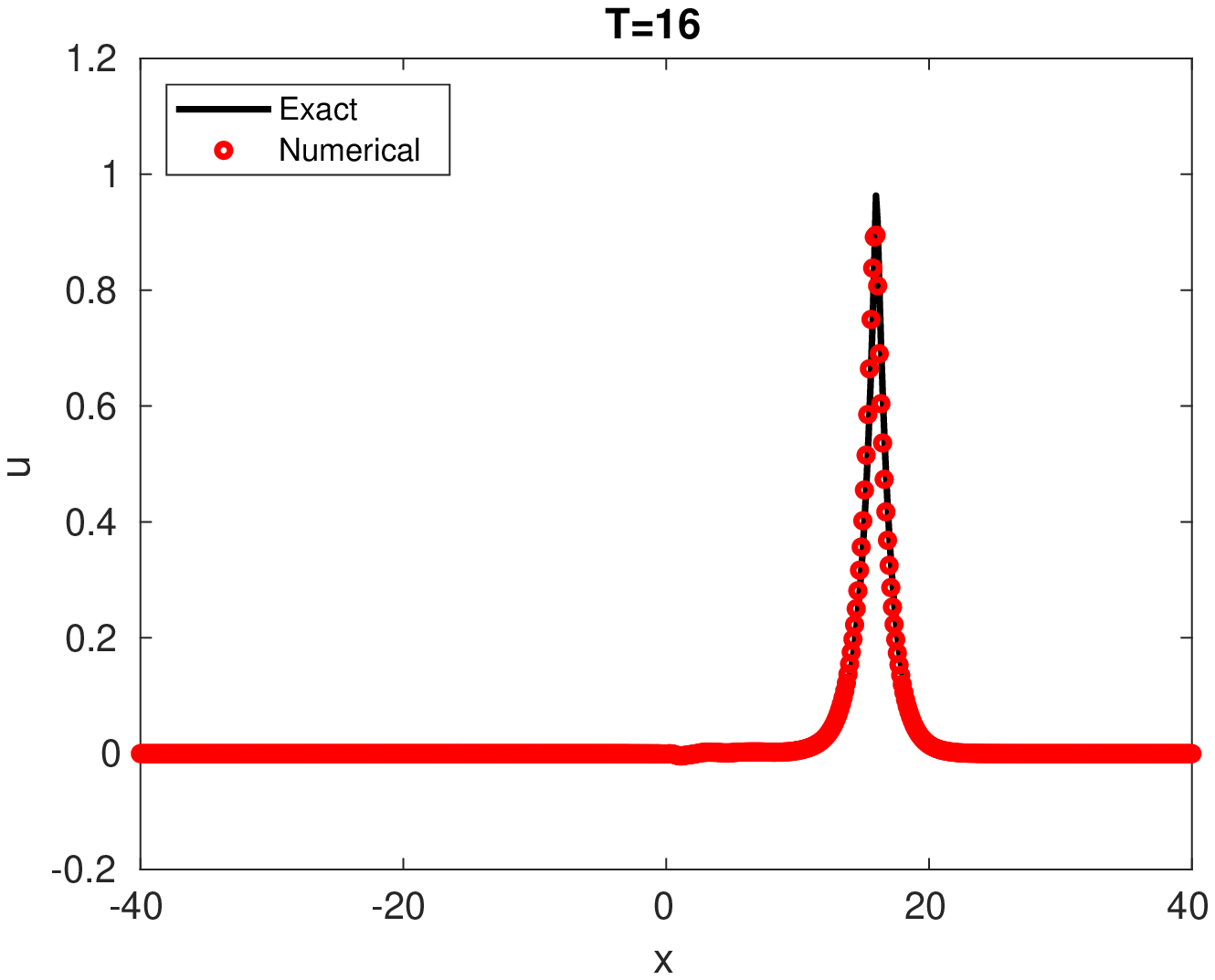}
	\caption{The single peakon solution of the DP equation in Example \ref{Ex:peakon_one}. $N = 640$. WENO5.}
	\label{Fig:peakon_NWENO5}
\end{figure} 	
We set the traveling speed $c = 1$ and the computational domain $[-40, 40]$.  In Figure \ref{Fig:peakon_NWENO5} 
and \ref{Fig:anti-peakon_NWENO5},
we show the peakon and anti-peakon profile at  $T = 4, 8, 12$ and $16$ with $N = 640$. We can see clearly that the moving peakon and anti-peakon profiles are well  resolved. There is no numerical oscillation near the wave crest.  {\color{red}We also observe
	 slight under-shoots for the peakon solution and slight over-shoots for the anti-peakon solution around $x=0$, which phenomenon is very common among many numerical methods
	 for the DP equation, such as those proposed in  \cite{Feng.Liu_JCP2009,  xu_local_2011, xia_weighted_2017}, etc.} 

%
%
%
%
\begin{figure}[!htbp]
\centering%
\includegraphics[width=0.45\textwidth]{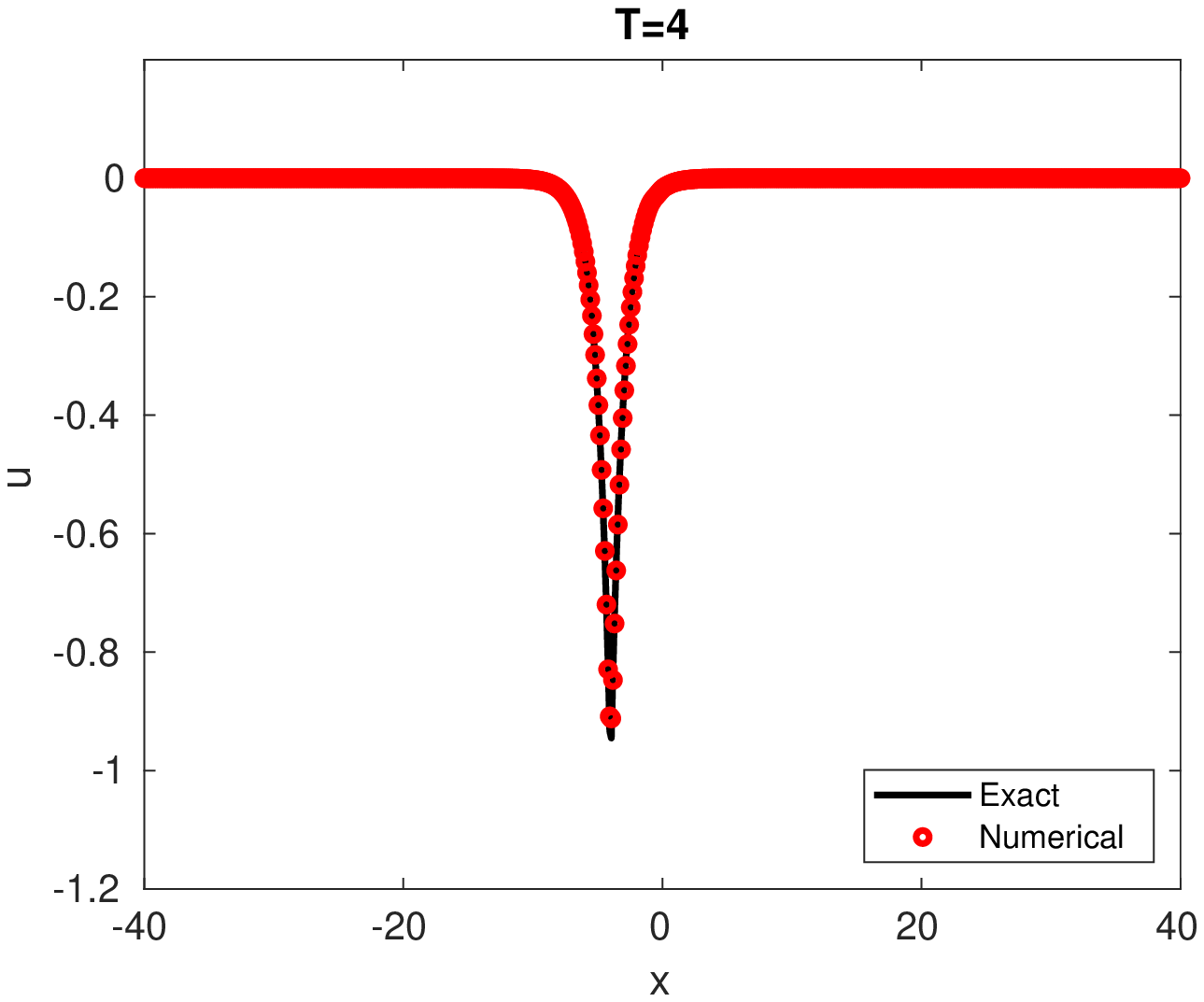} 
\includegraphics[width=0.45\textwidth]{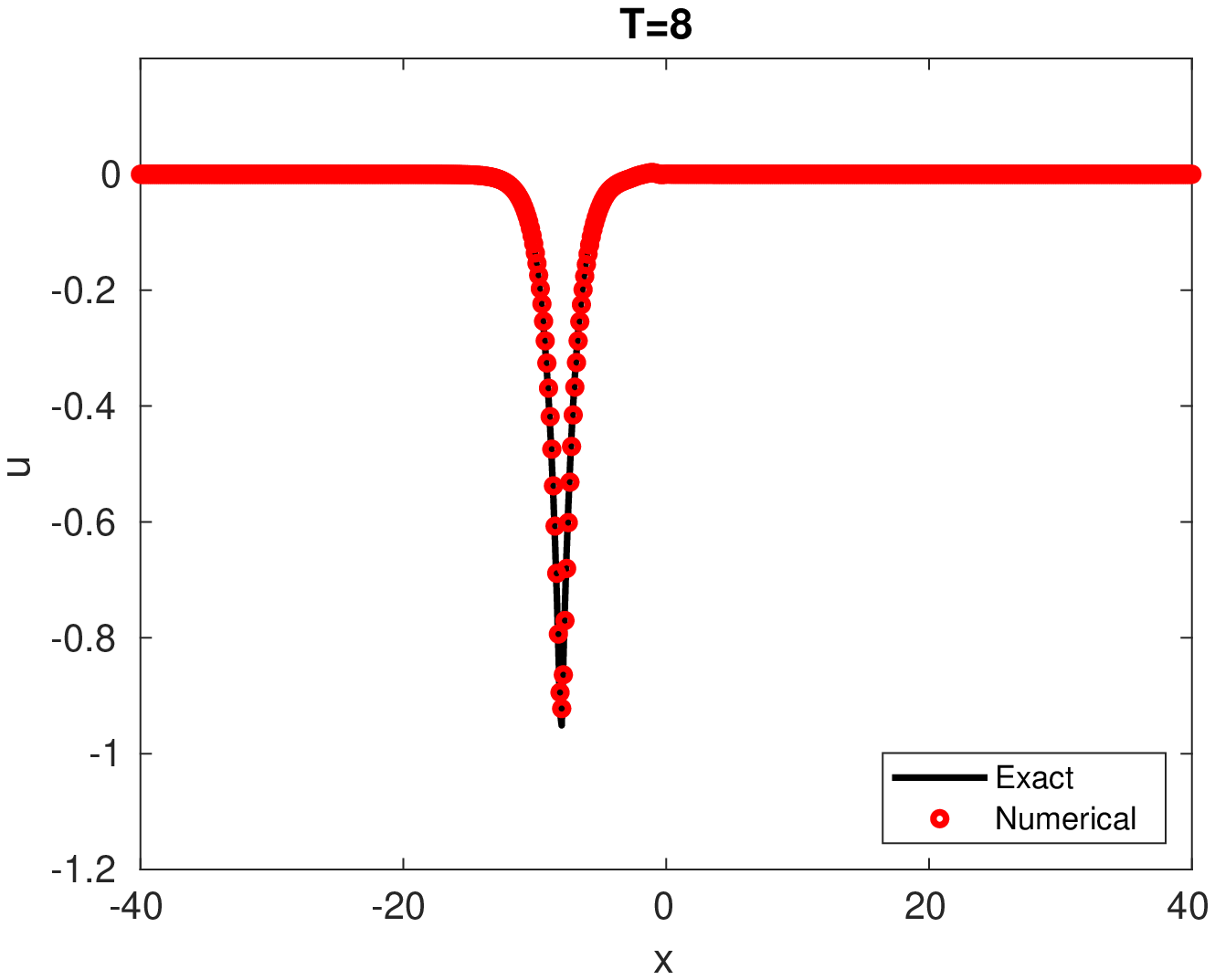}
\includegraphics[width=0.45\textwidth]{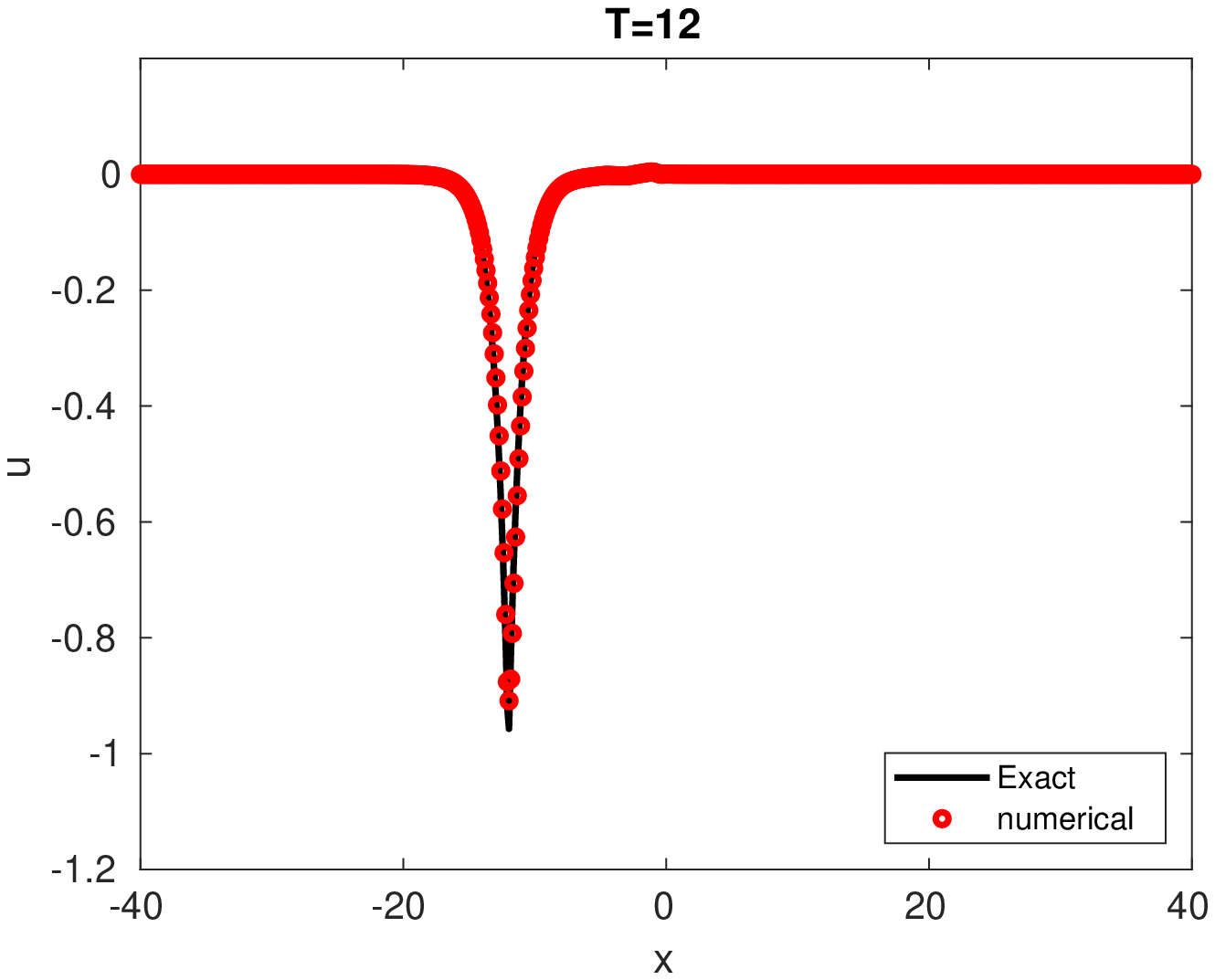} 
\includegraphics[width=0.45\textwidth]{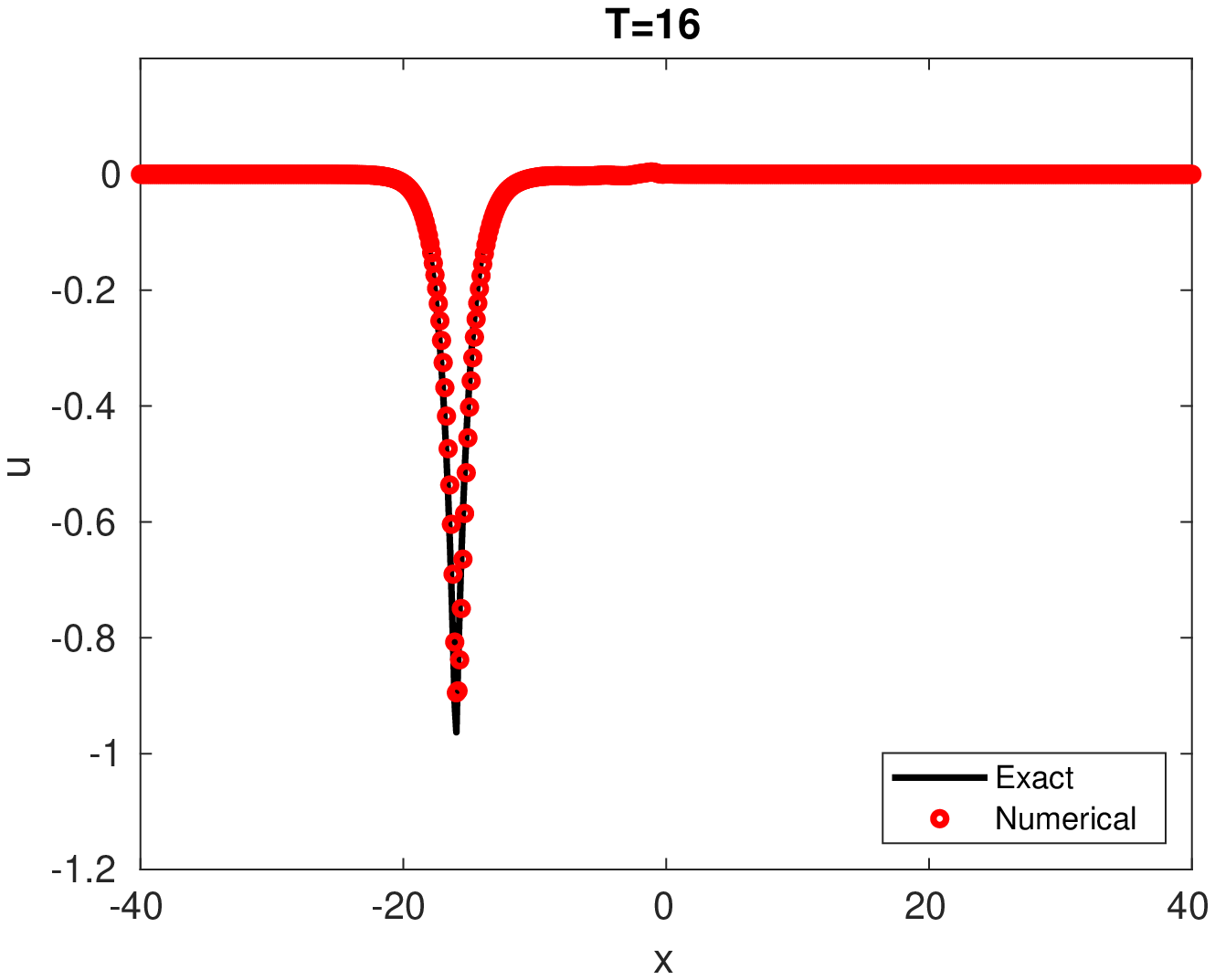}	
\caption{The single anti-peakon solution of the DP equation in Example \ref{Ex:peakon_one}. $N = 640$. WENO5.}
\label{Fig:anti-peakon_NWENO5}
\end{figure} 
%
	
%
\end{exa}

\begin{exa}
\label{Ex:twopeakon}
{\bf Two-peakon interaction and two-anti-peakon interaction}\\
In this example, we consider the two-peakon interaction \cite{Lundmark_JNS2007, Feng.Liu_JCP2009} of the DP equation  with the initial condition
\begin{equation*}
u(x,0)= c_{1}e^{-|x-x_{1}|}+c_{2}e^{-|x-x_{2}|},
\end{equation*}
and the two-anti-peakon interaction with the initial condition
\begin{align*}
u(x,0) = -c_{1}e^{-|x-x_{1}|}-c_{2}e^{-|x-x_{2}|}.
\end{align*}
In these interactions, the peakon should preserve its shape and velocity before and after encountering a nonlinear interaction with the other peakon. We set the parameters  as $c_{1} = 2$, $c_{2} = 1$, $x_{1} = -13.792,\ x_{2} = -4$, and the computational domain as $[-40, 40]$. In Figure \ref{Fig:peakon2_WENO5}
and \ref{Fig:anti-peakon2_WENO5},
we show  the two-peakon interaction and the two-anti-peakon interaction at $T =0, 4, 8, 12$ with $N = 1280$. {\color{red} The reference solutions are obtained by the  fifth order classical WENO scheme \cite{xia_weighted_2017} with $2560$ meshes.} We can see clearly that the moving peakons of both cases are well resolved.

%
	
\begin{figure}[!htbp]
\centering%
\includegraphics[width=0.45\textwidth]{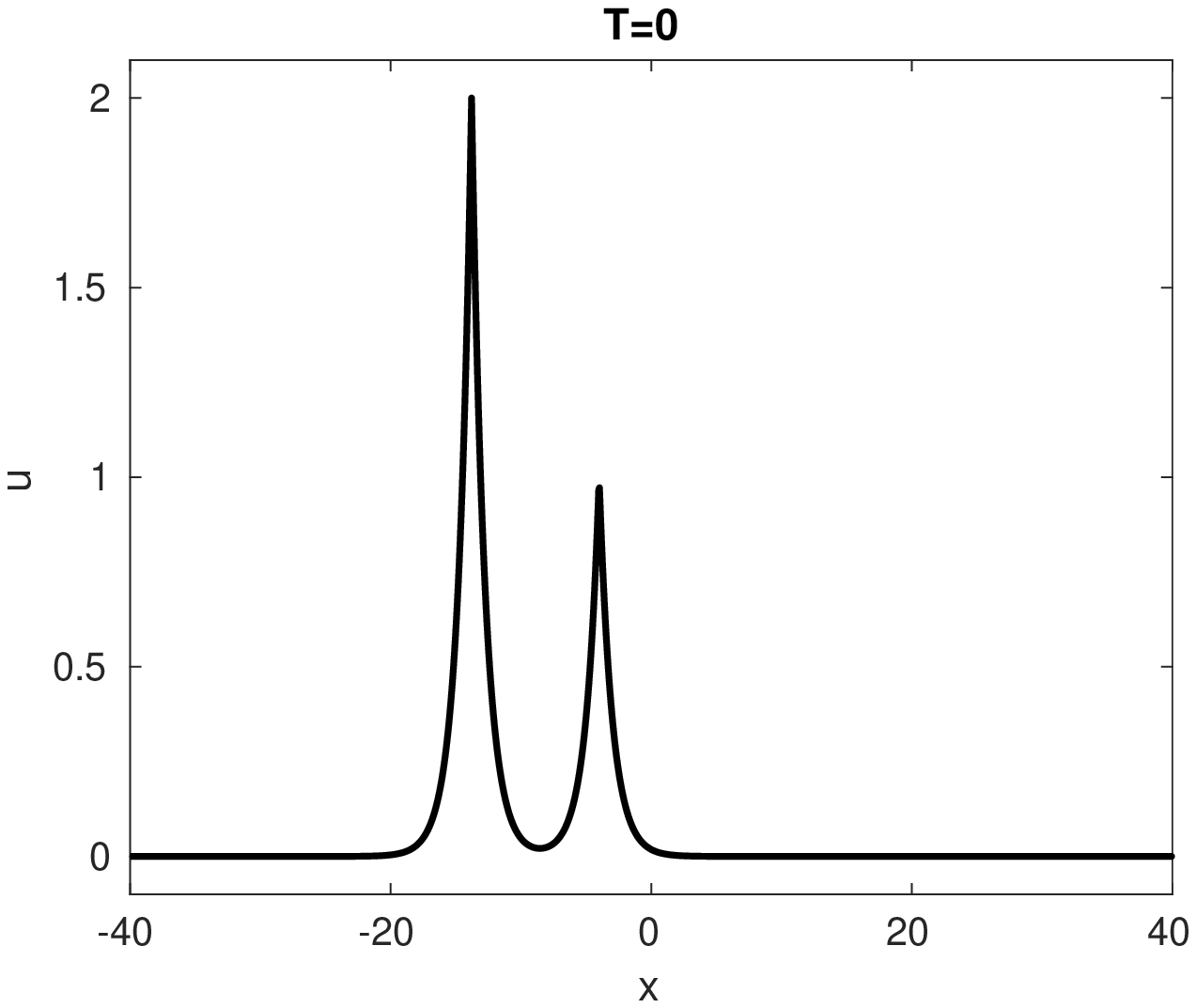} 
\includegraphics[width=0.45\textwidth]{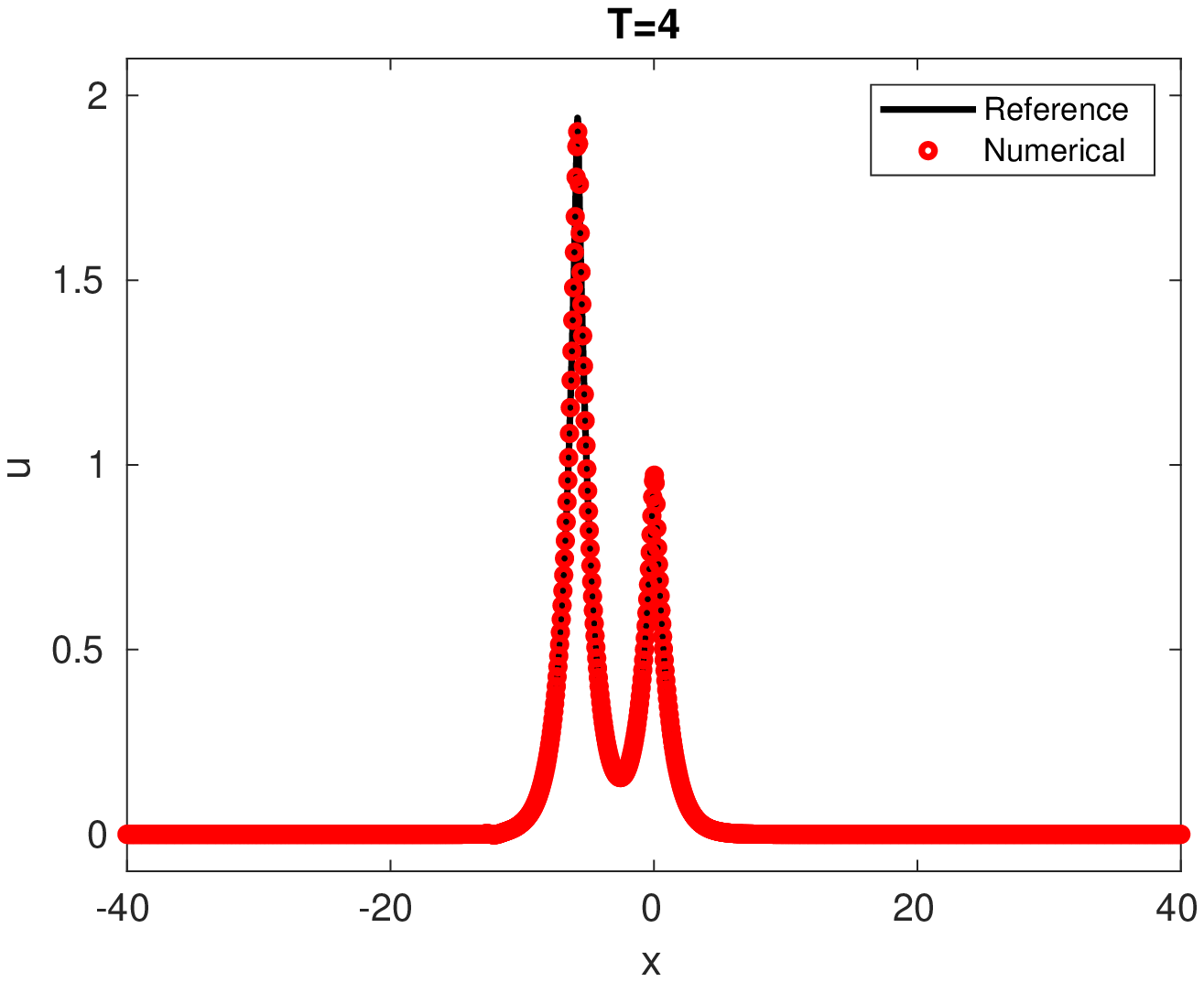}
\includegraphics[width=0.45\textwidth]{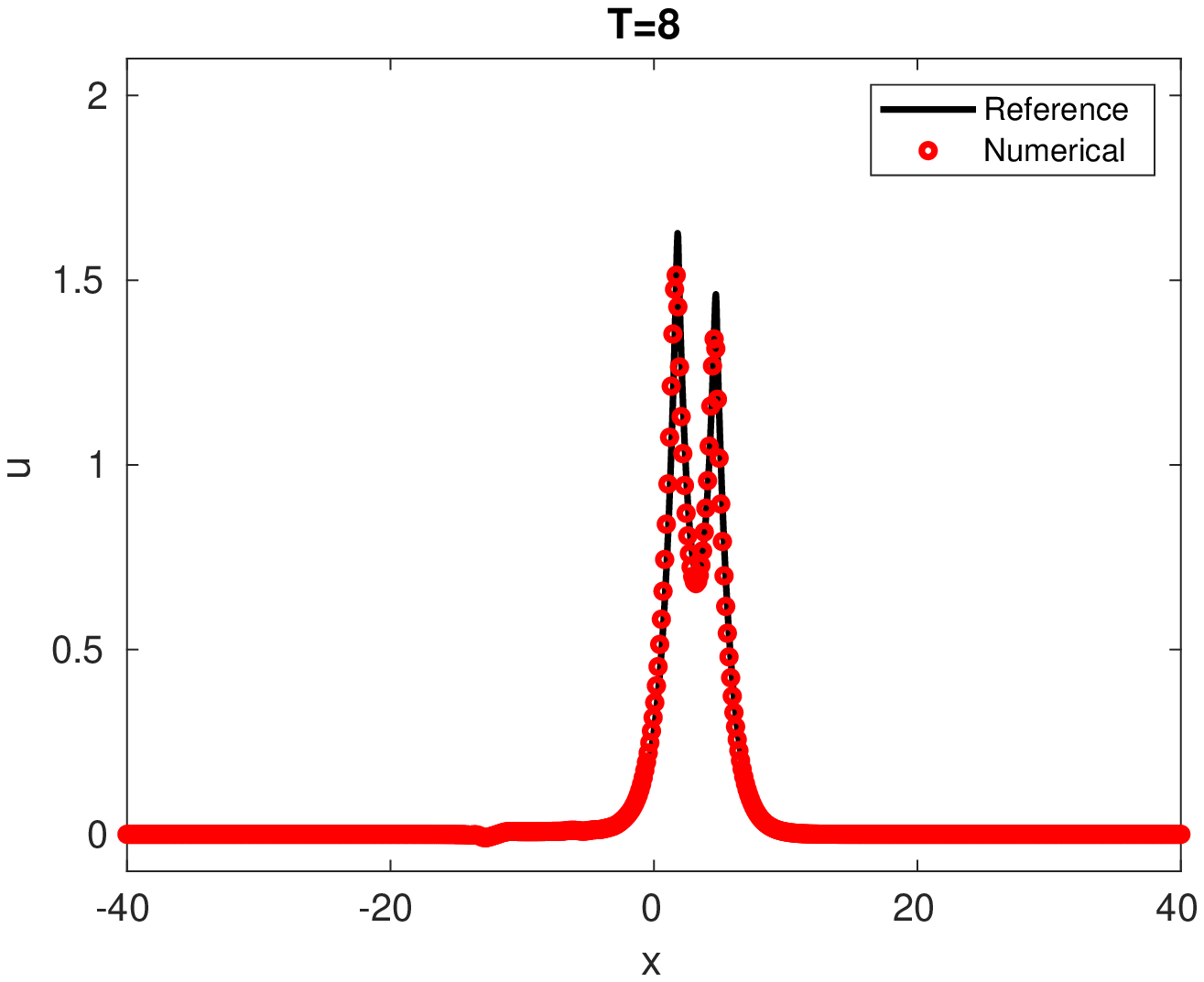} 
\includegraphics[width=0.45\textwidth]{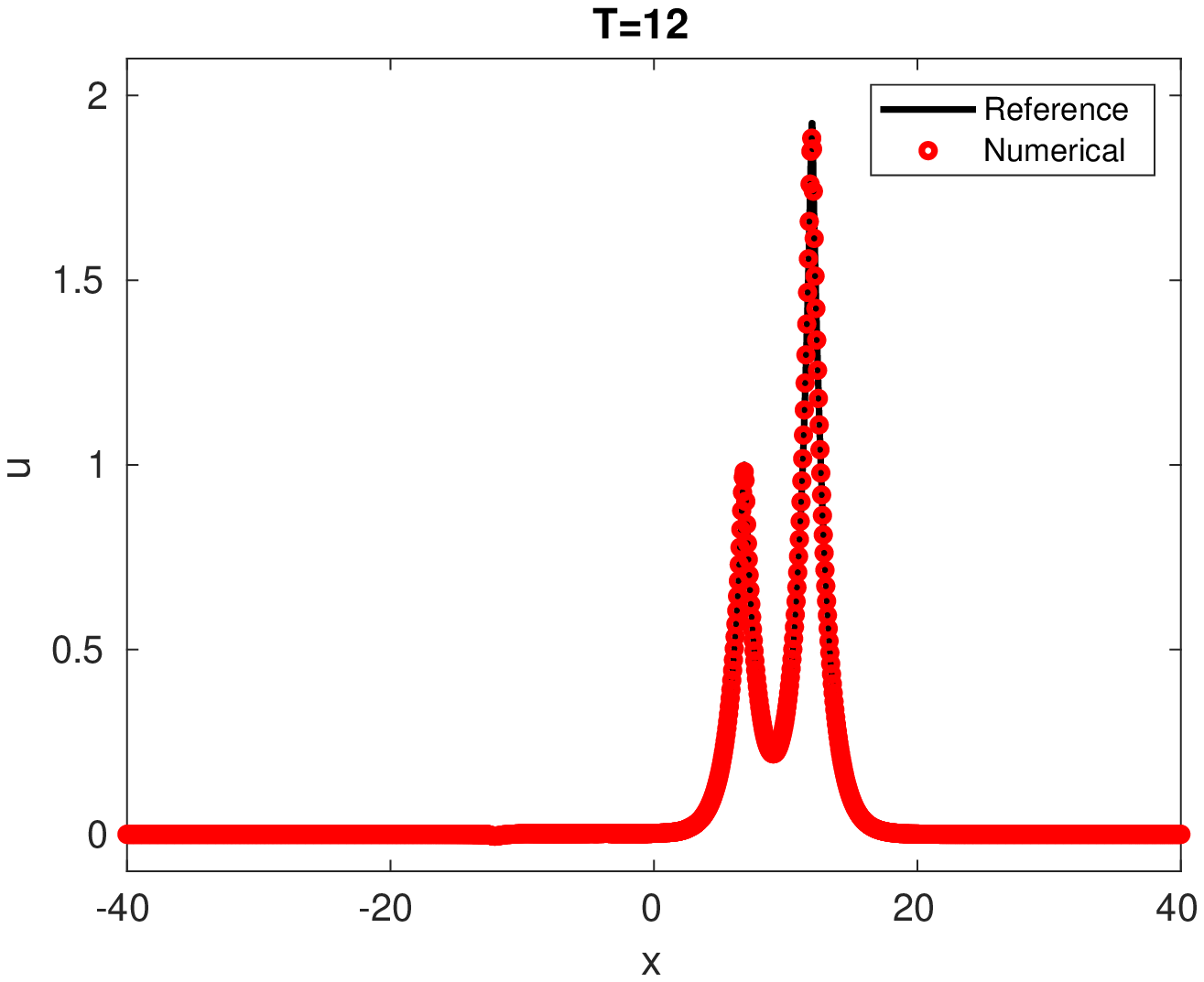}
\caption{{\color{red}Two-peakon interaction of the DP equation in Example \ref{Ex:twopeakon}.  $N = 1280$. MR-WENO5.}}
\label{Fig:peakon2_WENO5}
\end{figure} 
%
%
	
%
%
	
\begin{figure}[!htbp]
\centering%
\includegraphics[width=0.45\textwidth]{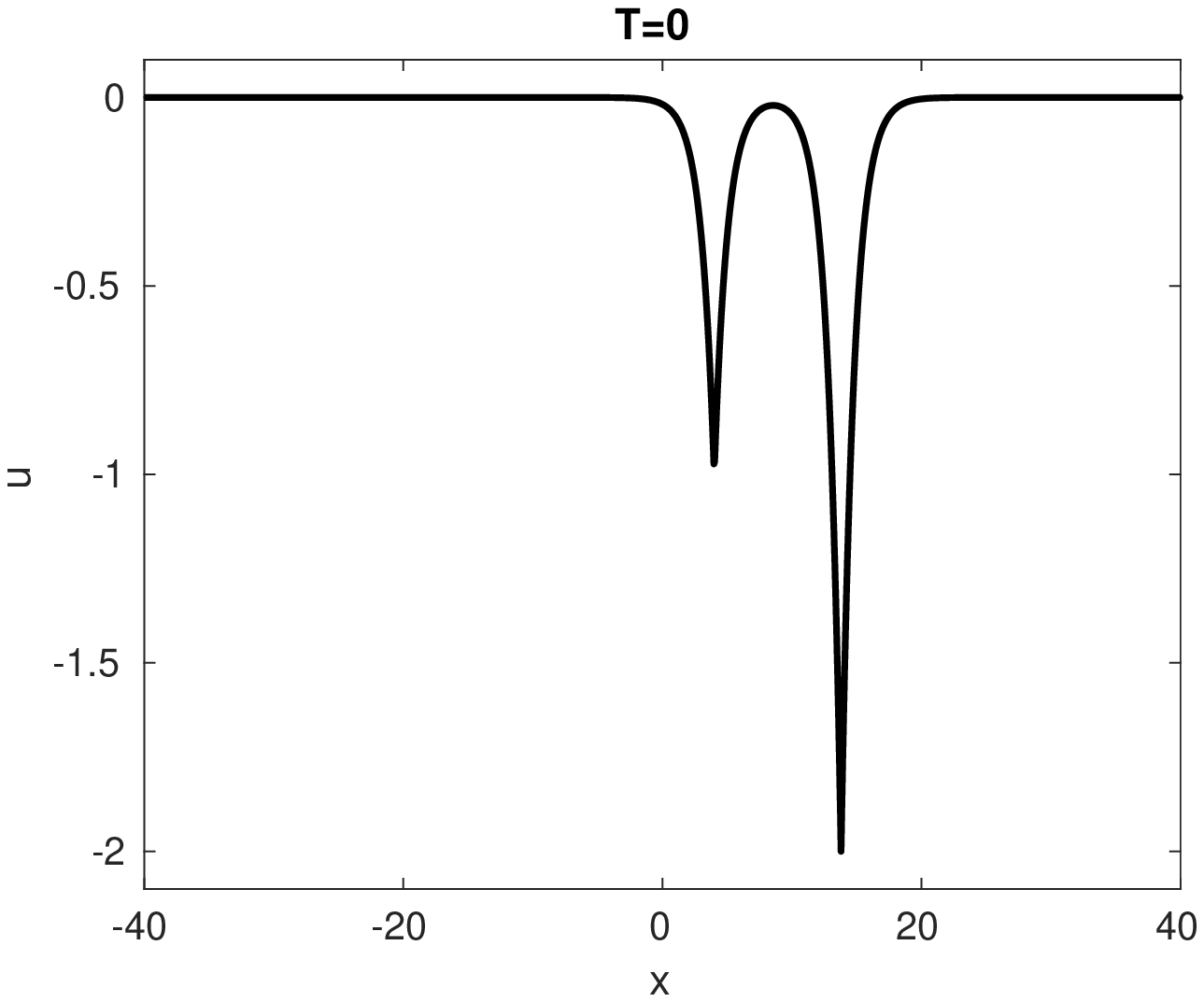} 
\includegraphics[width=0.45\textwidth]{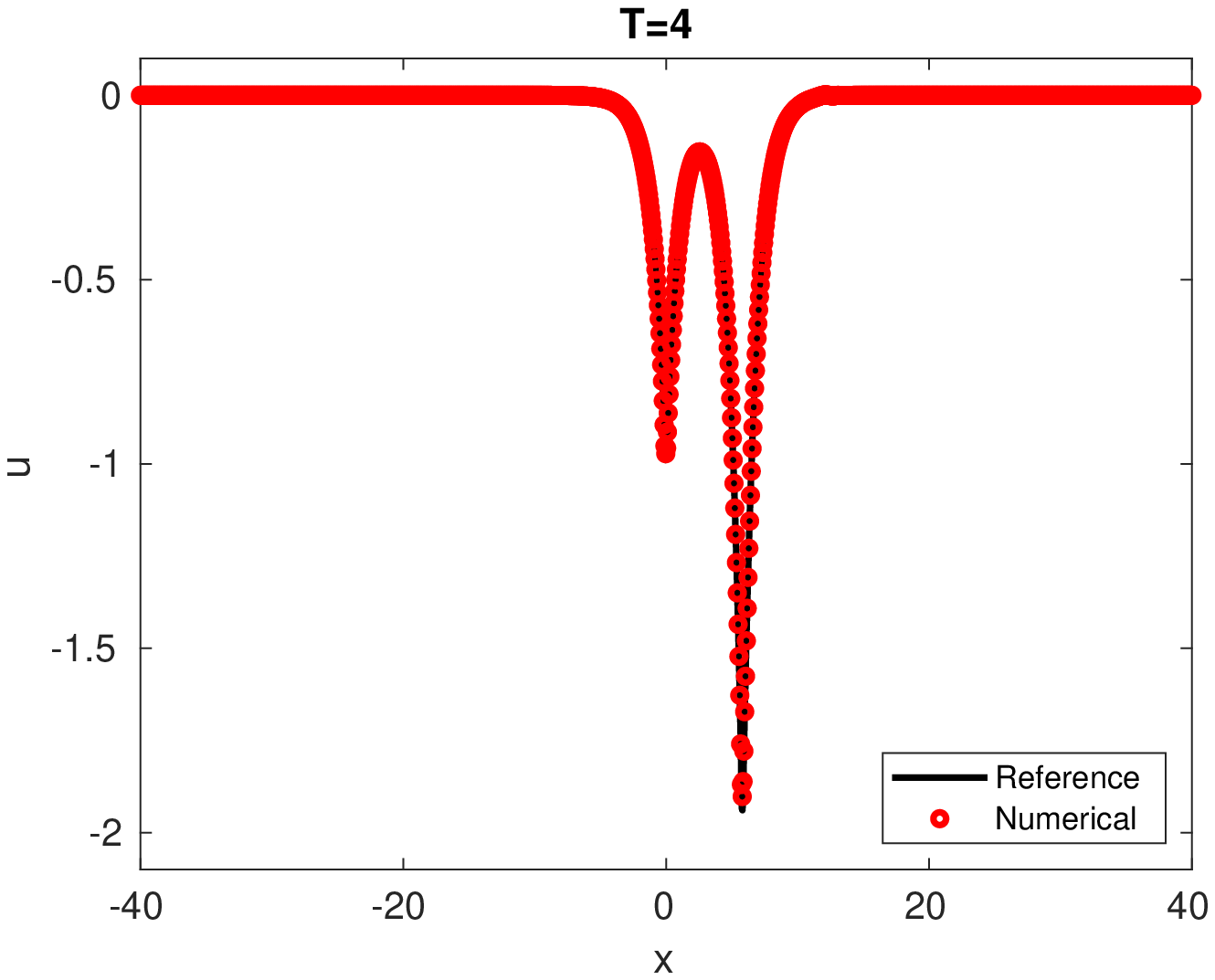}
\includegraphics[width=0.45\textwidth]{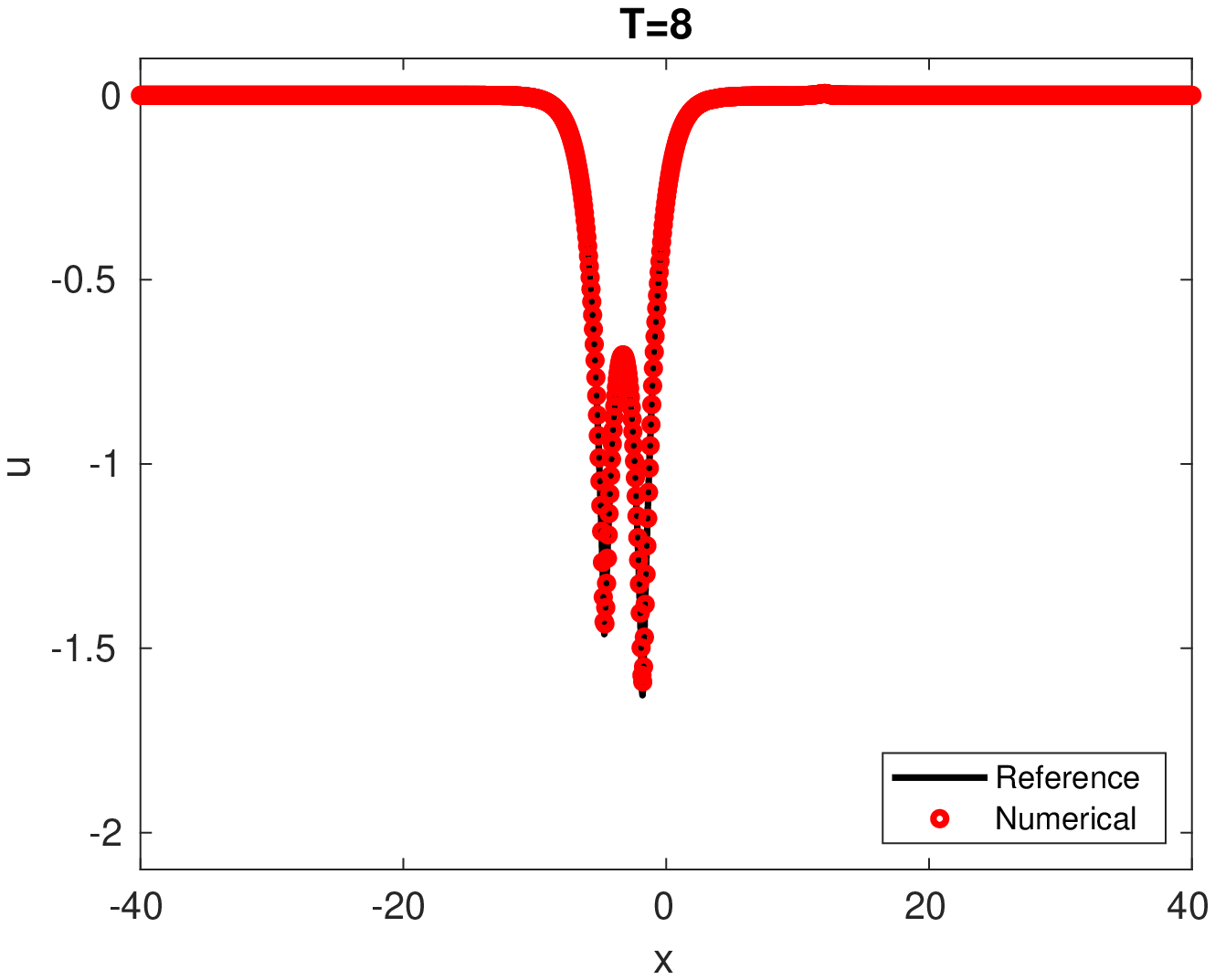} 
\includegraphics[width=0.45\textwidth]{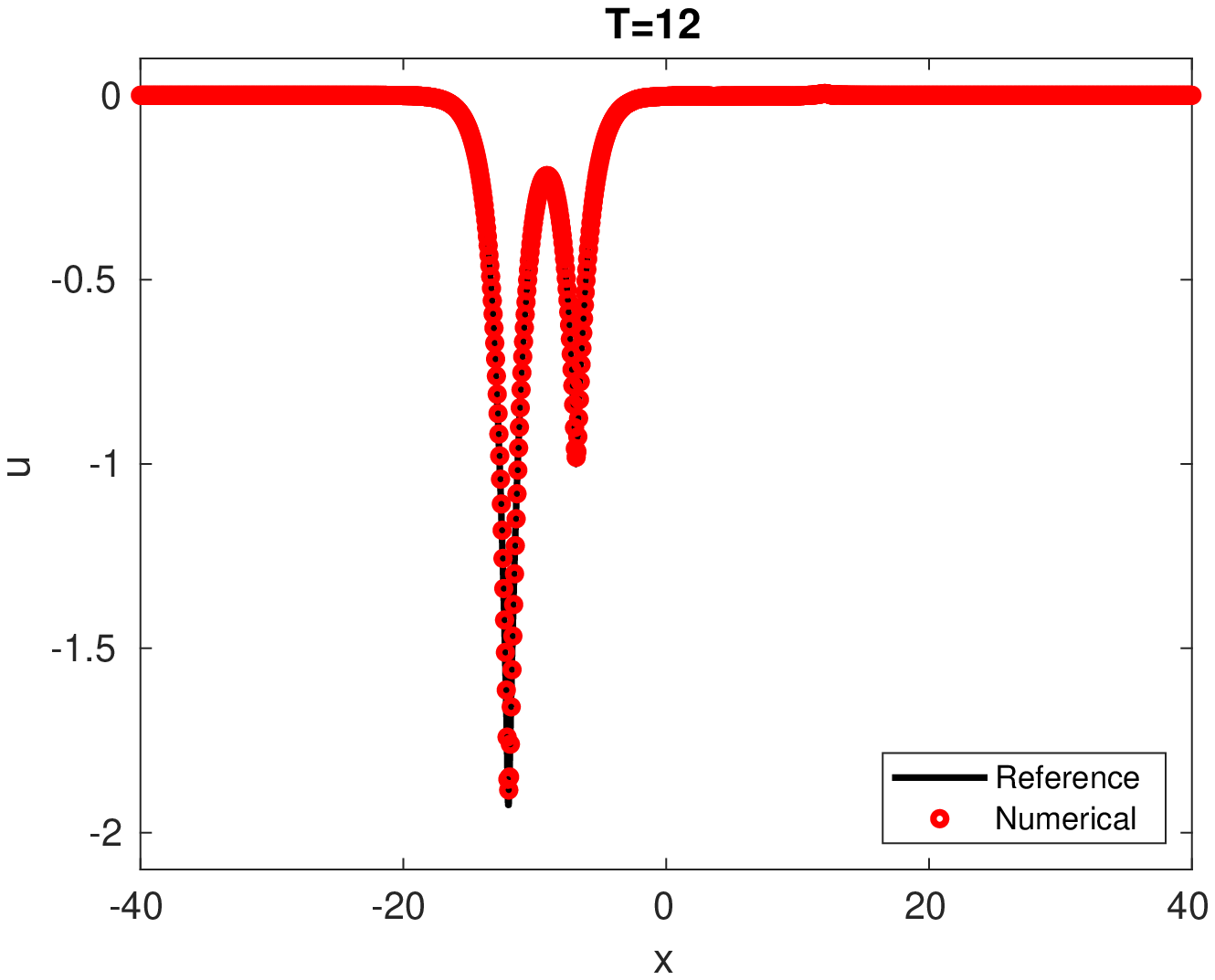}
\caption{{\color{red}Two-anti-peakon interaction of the DP equation in Example \ref{Ex:twopeakon}.  $N = 1280$. MR-WENO5.}}
\label{Fig:anti-peakon2_WENO5}
\end{figure} 
%
%
\end{exa}

\begin{exa}
\label{Ex:shockpeakon}
{\bf Shock peakon solution}\\
In this example, we consider the DP equation  with the shock peakon solution
\begin{equation}
u(x, t) = -\frac{1}{t+1}\text{sign}(x)e^{-|x|}.
\end{equation}
The computational domain is taken as $[-25, 25]$. Figure \ref{Fig:shock-peakon_NWENO5}
shows the numerical solutions obtained by WENO5, MR-WENO5 and MR-WENO7 at $T = 3$ and $6$  with $N = 640$. We observe that  there is no numerical oscillation near the discontinuity at $x=0$ and the shock interface at $T=3$ is very sharp for all three schemes.  At $T=6$, the shock interface obtained by the seventh order scheme MR-WENO7 is better resolved than the fifth order schemes, e.g. WENO5, MR-WENO5 shown in Figure \ref{Fig:shock-peakon_NWENO5} and the classical WENO method in \cite{xia_weighted_2017}.
\begin{figure}
\centering%
\includegraphics[width=0.45\textwidth]{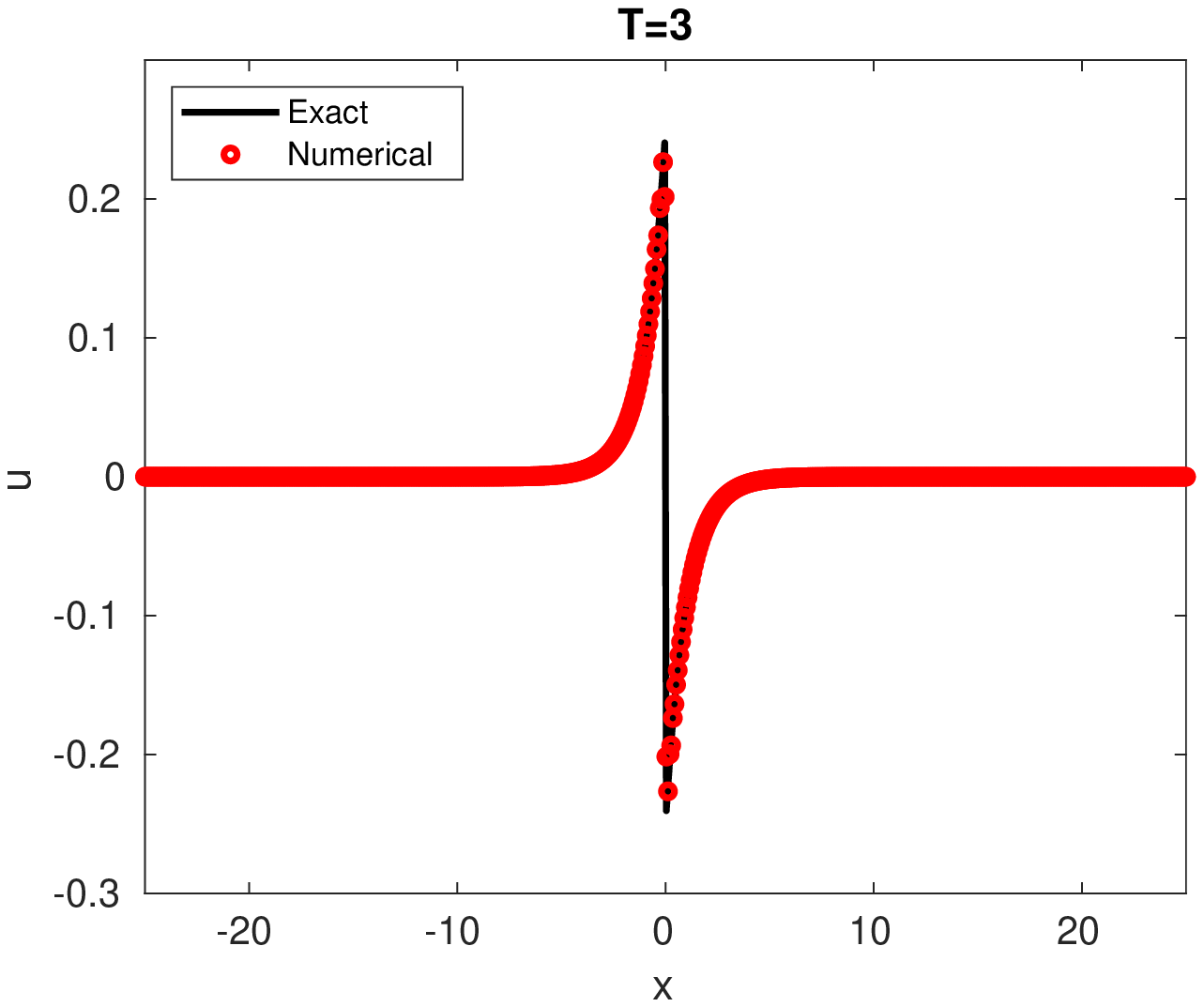} 
\includegraphics[width=0.45\textwidth]{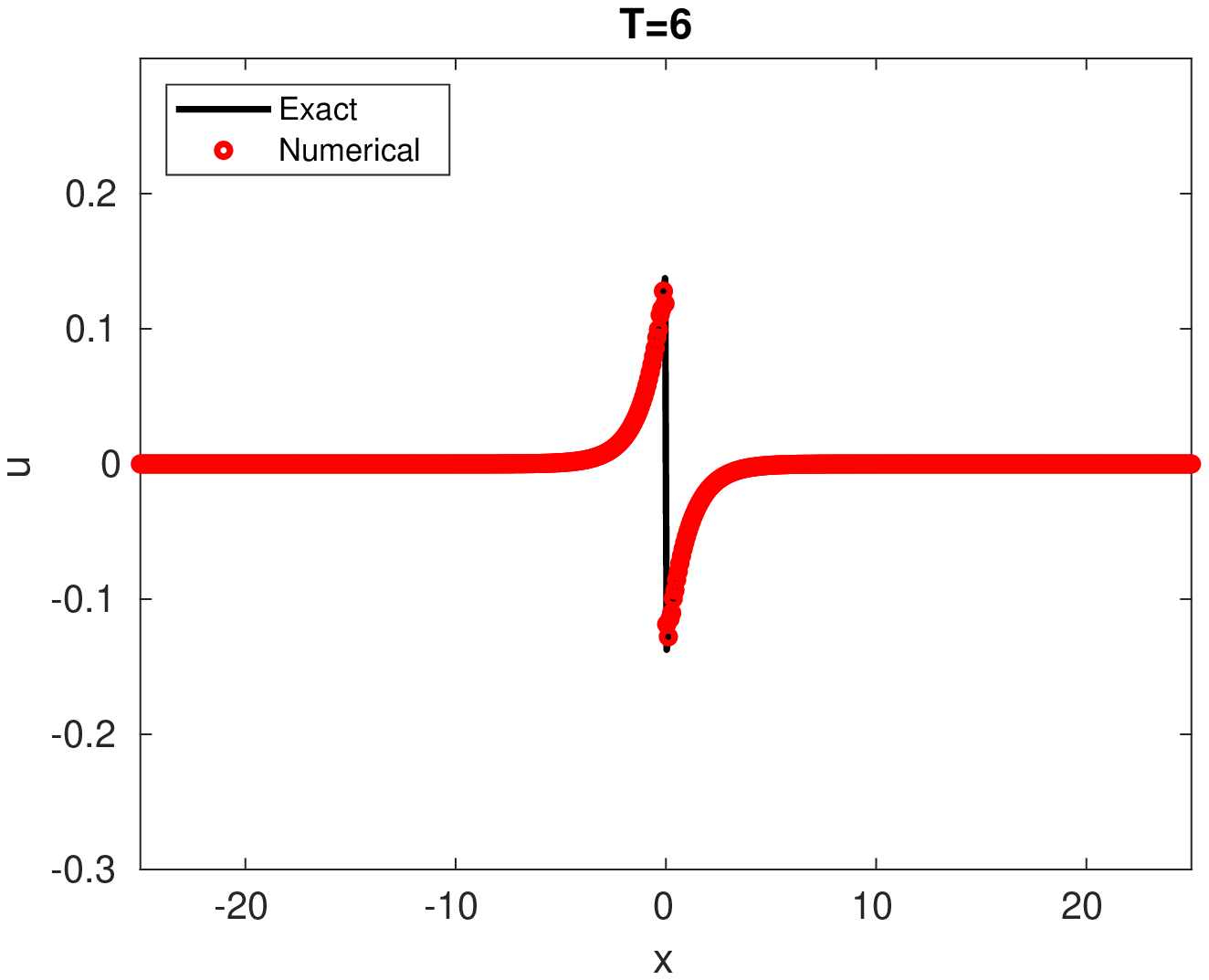}
\includegraphics[width=0.45\textwidth]{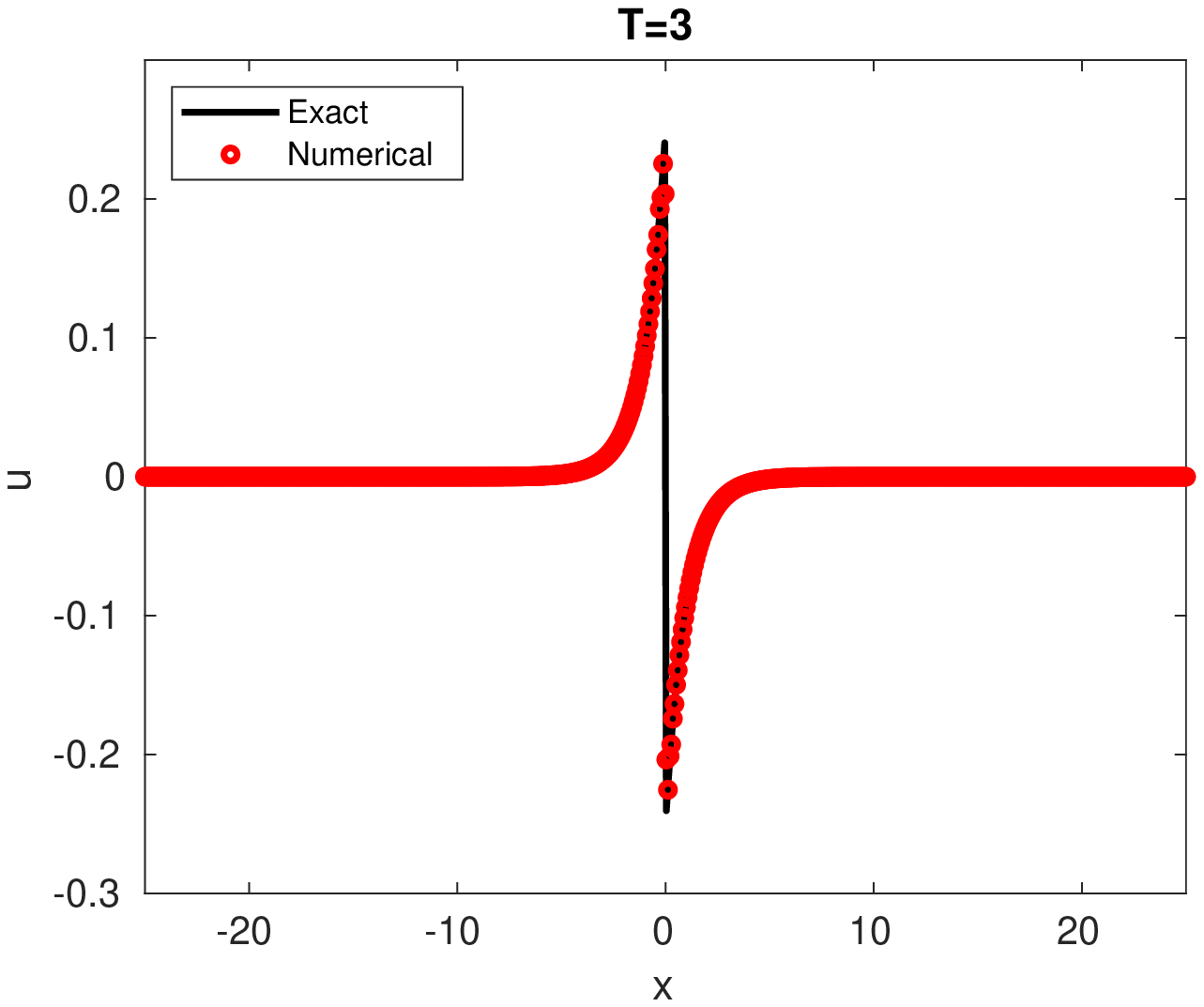} 
\includegraphics[width=0.45\textwidth]{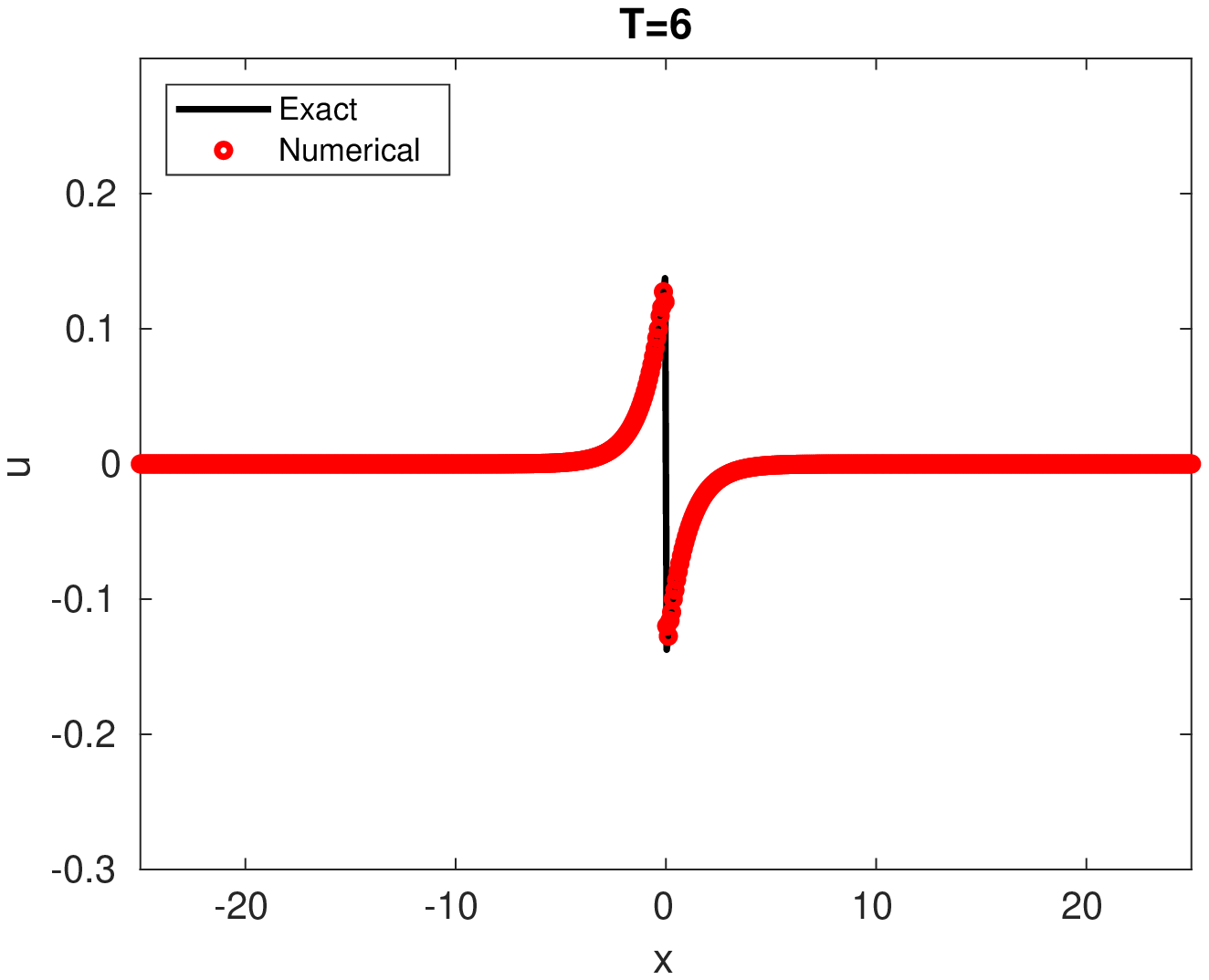}
\includegraphics[width=0.45\textwidth]{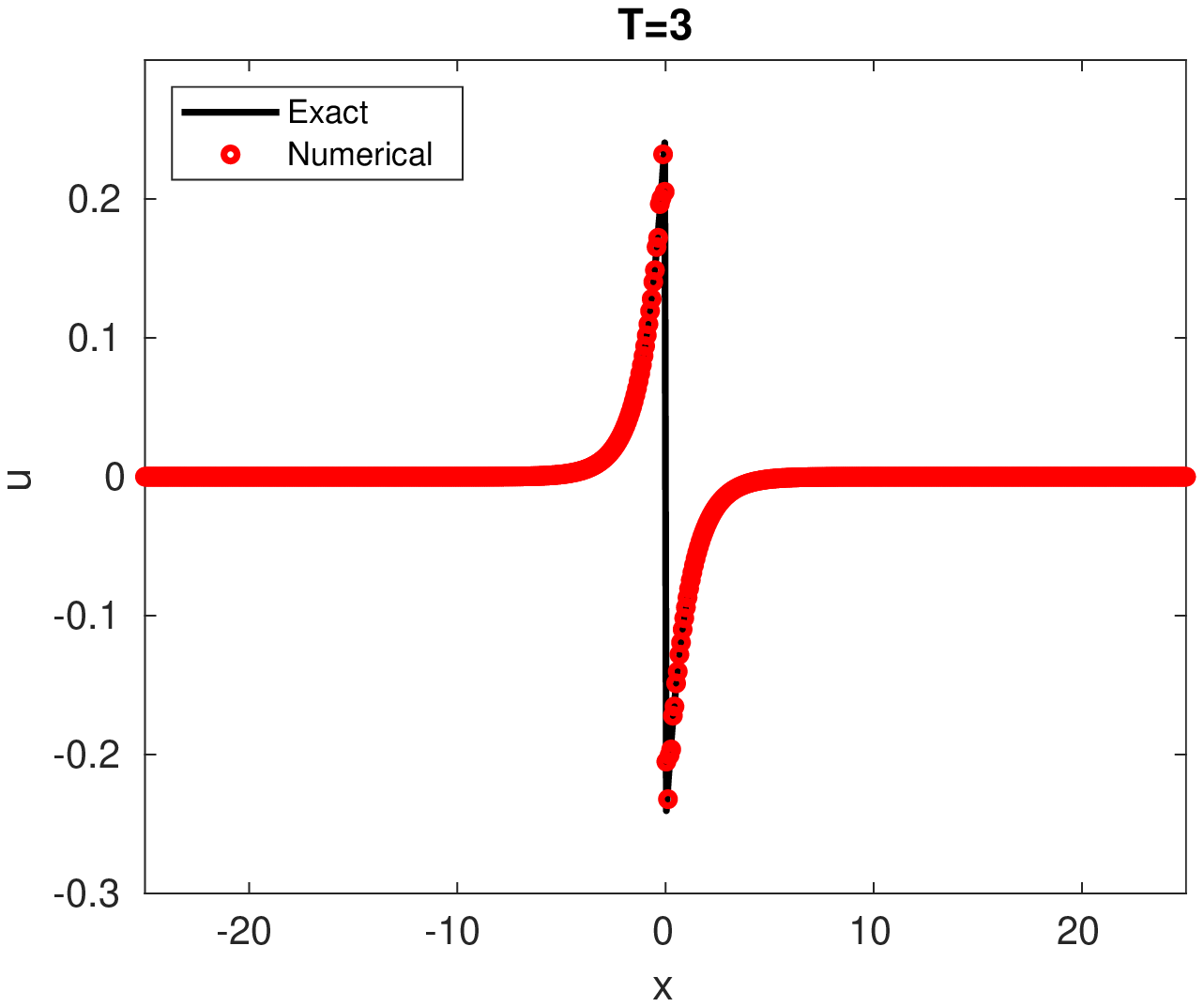} 
\includegraphics[width=0.45\textwidth]{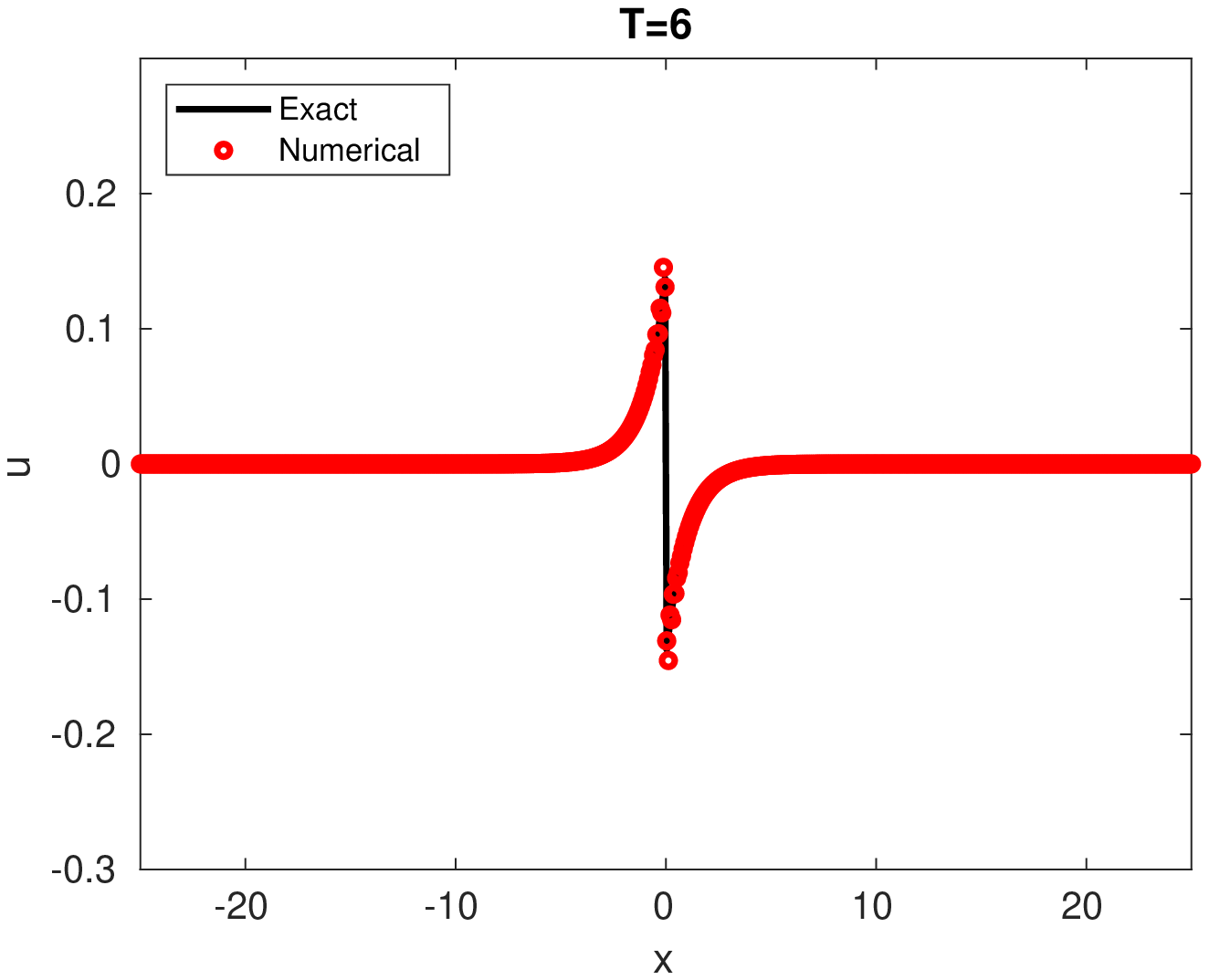}
\caption{The shock peakon solution of the DP equation in Example \ref{Ex:shockpeakon}. 
Top: WENO5; Middle: MR-WENO5; Bottom: MR-WENO7. $N = 640$.}
\label{Fig:shock-peakon_NWENO5}
\end{figure} 
	
%
%
%
%
%
	
%
%
%
%
\end{exa}

\begin{exa}
\label{Ex:peakon+anti-peakon}
{\bf Peakon and anti-peakon interaction}\\	
In this example, we consider the peakon and anti-peakon interaction \cite{Lundmark_JNS2007,Feng.Liu_JCP2009} of the DP equation with the initial condition
\begin{equation*}
u(x, 0) = e^{-|x+5|}-e^{-|x-5|}.
\end{equation*}
The computational domain is set as $[-20, 20]$. In this case,  a shock peakon is formed at $t\approx 5$, see e.g.  \cite{Lundmark_JNS2007,Feng.Liu_JCP2009, xu_local_2011} for more details.
We plot the numerical solutions at $T=0,\ 4,\ 5$ and $T=7$ with $N = 640$  in Figure \ref{Fig:peakon_anti-peakon_WENO7}. {\color{red} The reference solutions are obtained by the fifth order classical WENO scheme \cite{xia_weighted_2017} with 2560 meshes.} Again, there is no numerical oscillation during the peakon and anti-peakon interaction, and the shock interface is well resolved.
%
%
%
\begin{figure}[!htbp]
\centering%
\includegraphics[width=0.45\textwidth]{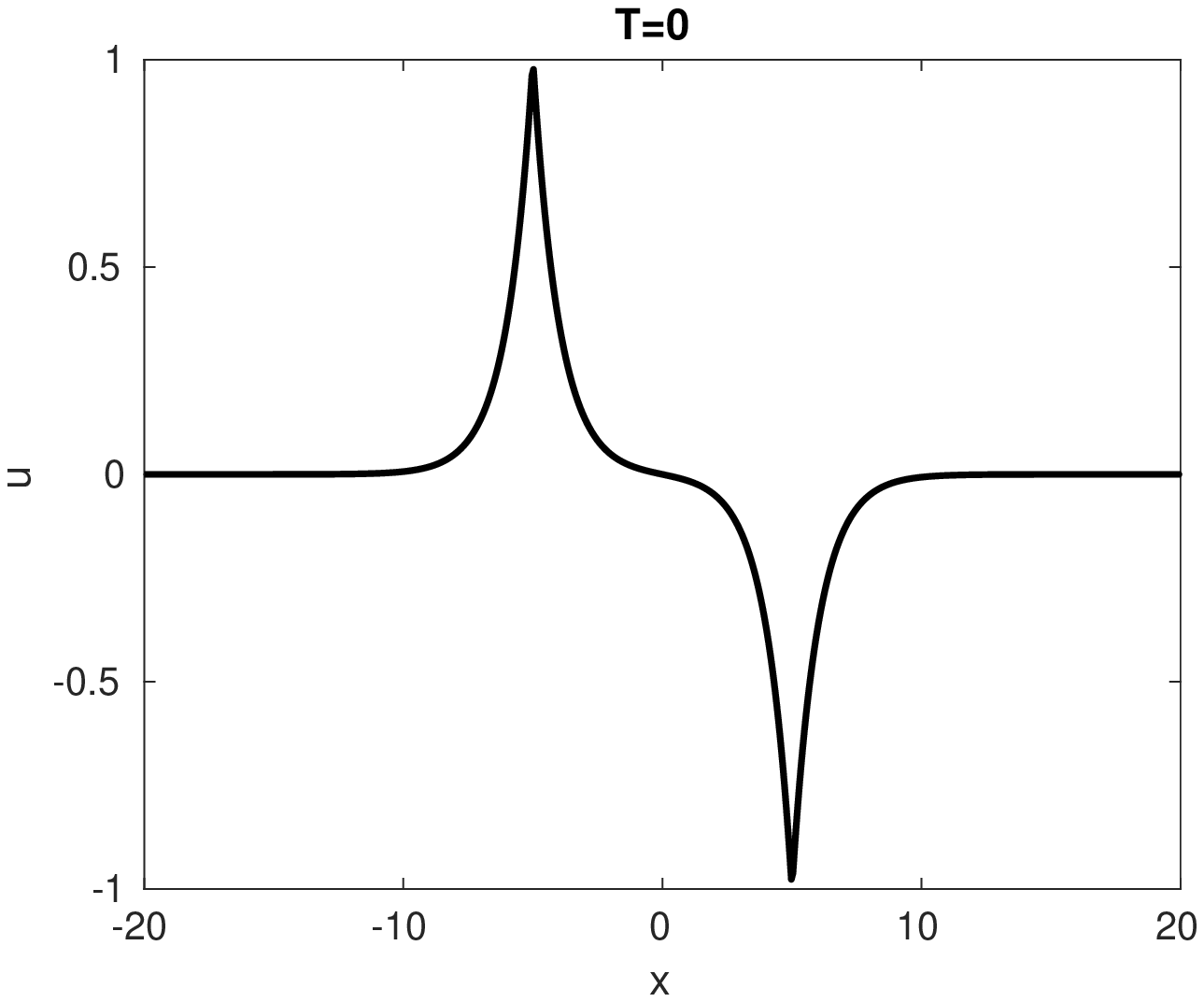} 
\includegraphics[width=0.45\textwidth]{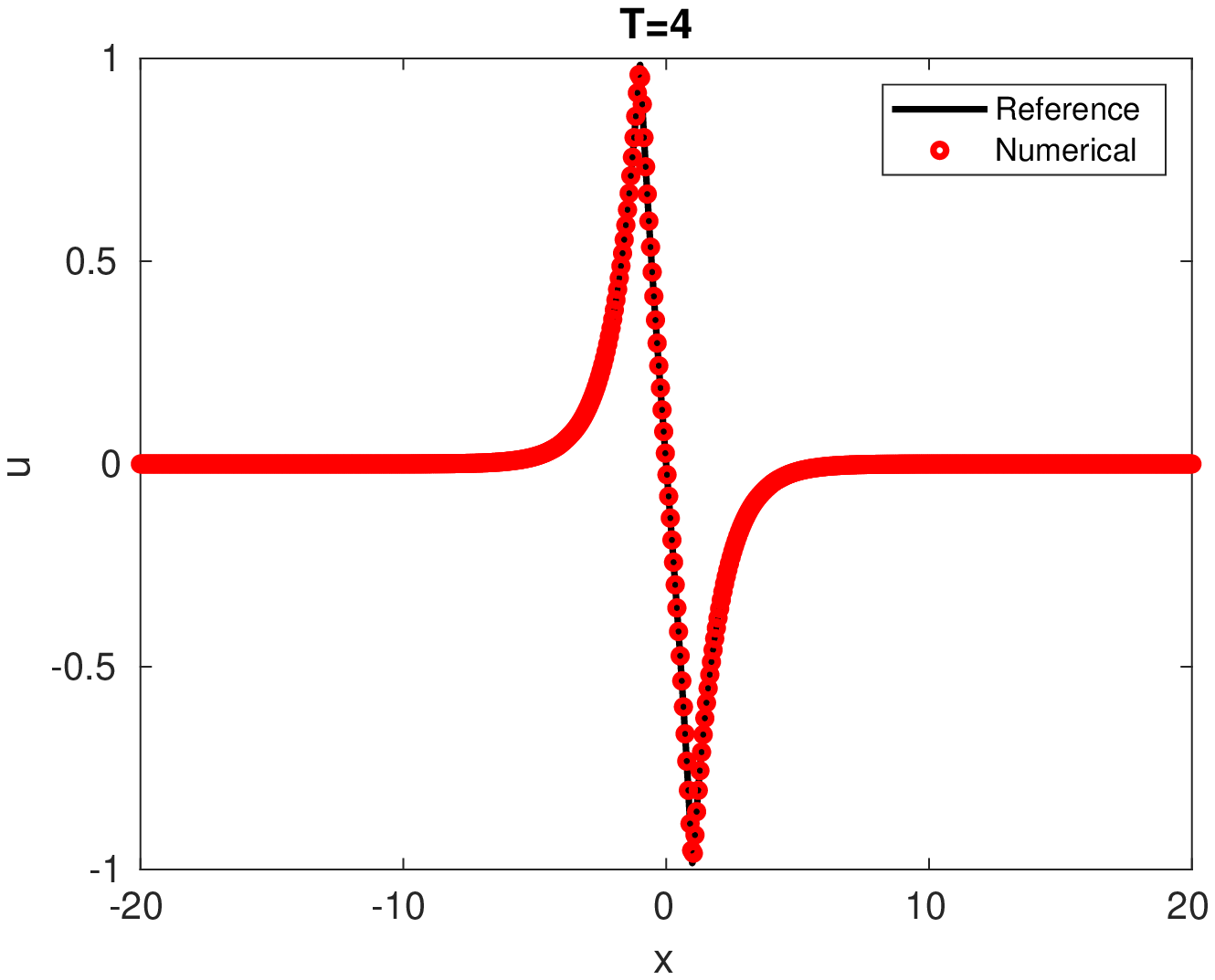}
\includegraphics[width=0.45\textwidth]{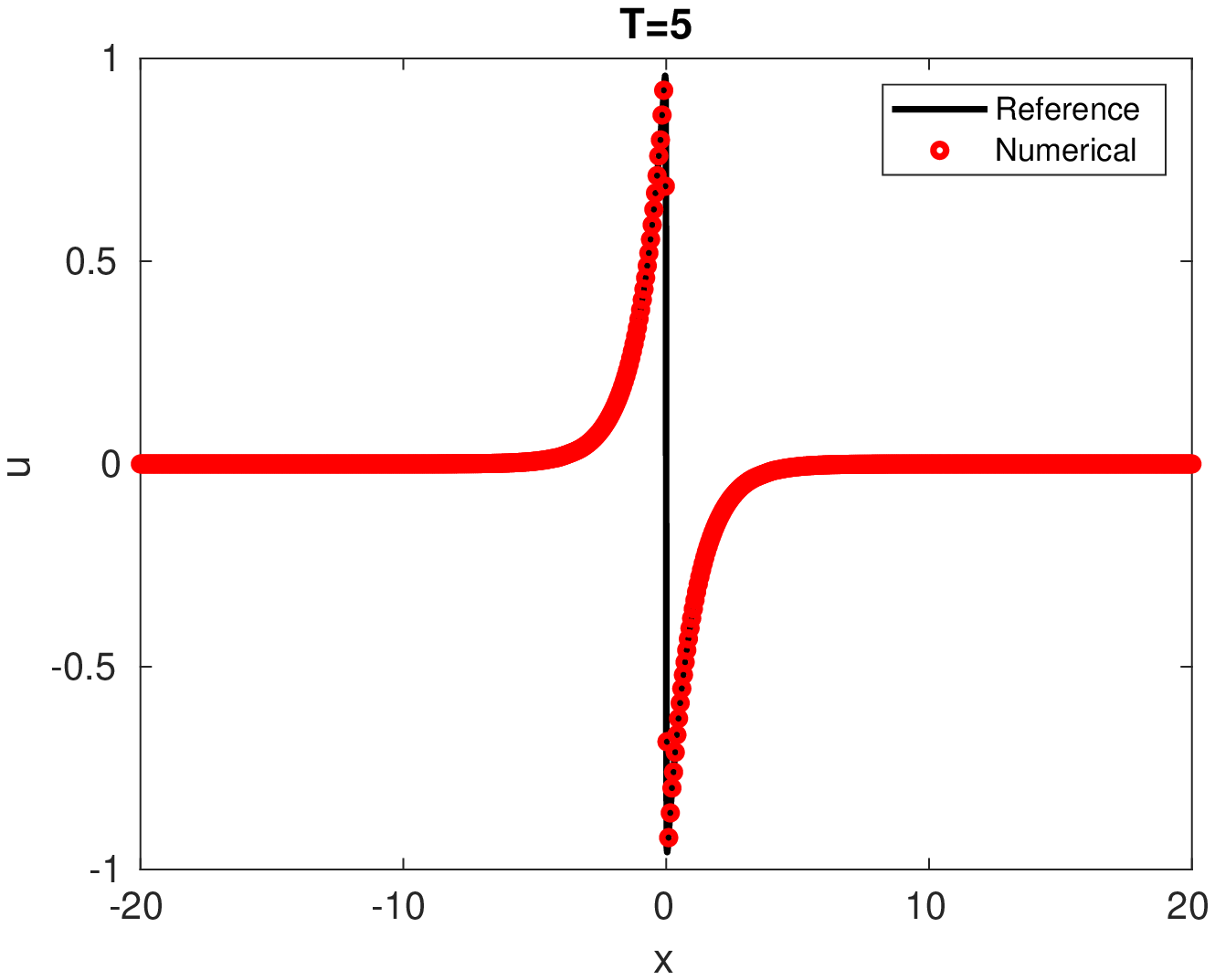}
\includegraphics[width=0.45\textwidth]{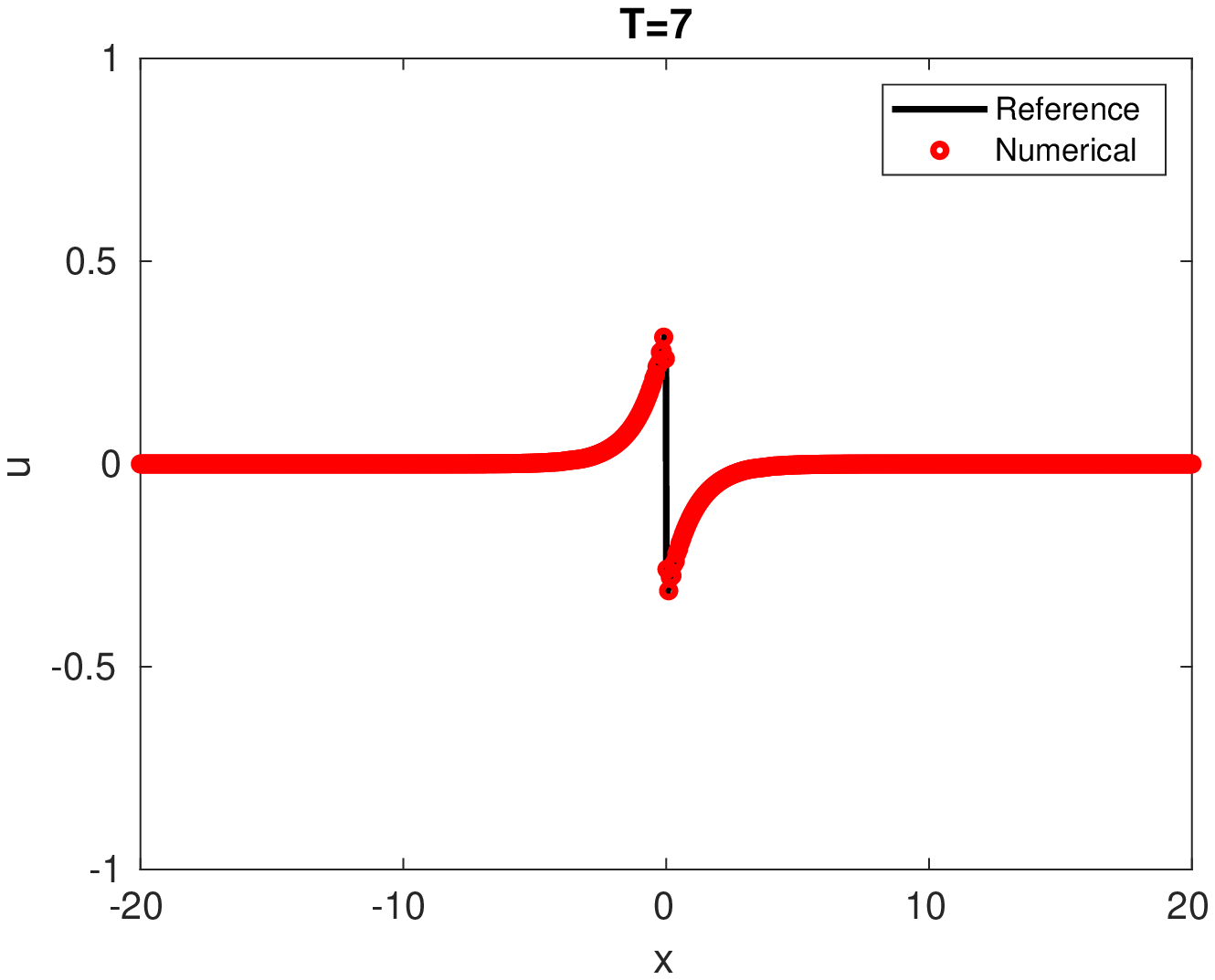}

\caption{{\color{red}The peakon  and anti-peakon interaction of the DP equation in Example \ref{Ex:peakon+anti-peakon}.  $N = 640$. MR-WENO7.}}
\label{Fig:peakon_anti-peakon_WENO7}
\end{figure} 
\end{exa}

\begin{exa}\label{Ex:triple}
{\bf Triple interaction}\\
In this example, we consider a triple interaction among a peakon, an anti-peakon and a stationary shock peakon  of the DP equation with the initial condition
\begin{equation*}
u(x, 0) = e^{-|x+5|}+\text{sign}(x)e^{-|x|}-e^{-|x-5|}.
\end{equation*}
This example was theoretically studied in \cite{Lundmark_JNS2007} and numerically in \cite{Coclite.Karlsen.Risebro_JNA2008,Feng.Liu_JCP2009,xu_local_2011,xia_weighted_2017}, etc.
The exact solution is a triple collision among a peakon, an anti-peakon and a shock peakon  when $T\approx 5.32$.  The second shock peakon is formed when $T>5.32$. 
The simulations are carried out in the domain  $[-20, 20]$ with $N=640$ up to $T=7$. We show the numerical solutions at $T = 0,\ 2,\ 5.32$ and $7$ in Figure \ref{Fig:triple_NWENO5}. {\color{red} The reference solutions are obtained by the fifth order classical WENO scheme with 2560 meshes.} Again, there is no numerical oscillation during the triple interaction, and the shock interface is well resolved.
\begin{figure}[!htbp]
\centering%
\includegraphics[width=0.45\textwidth]{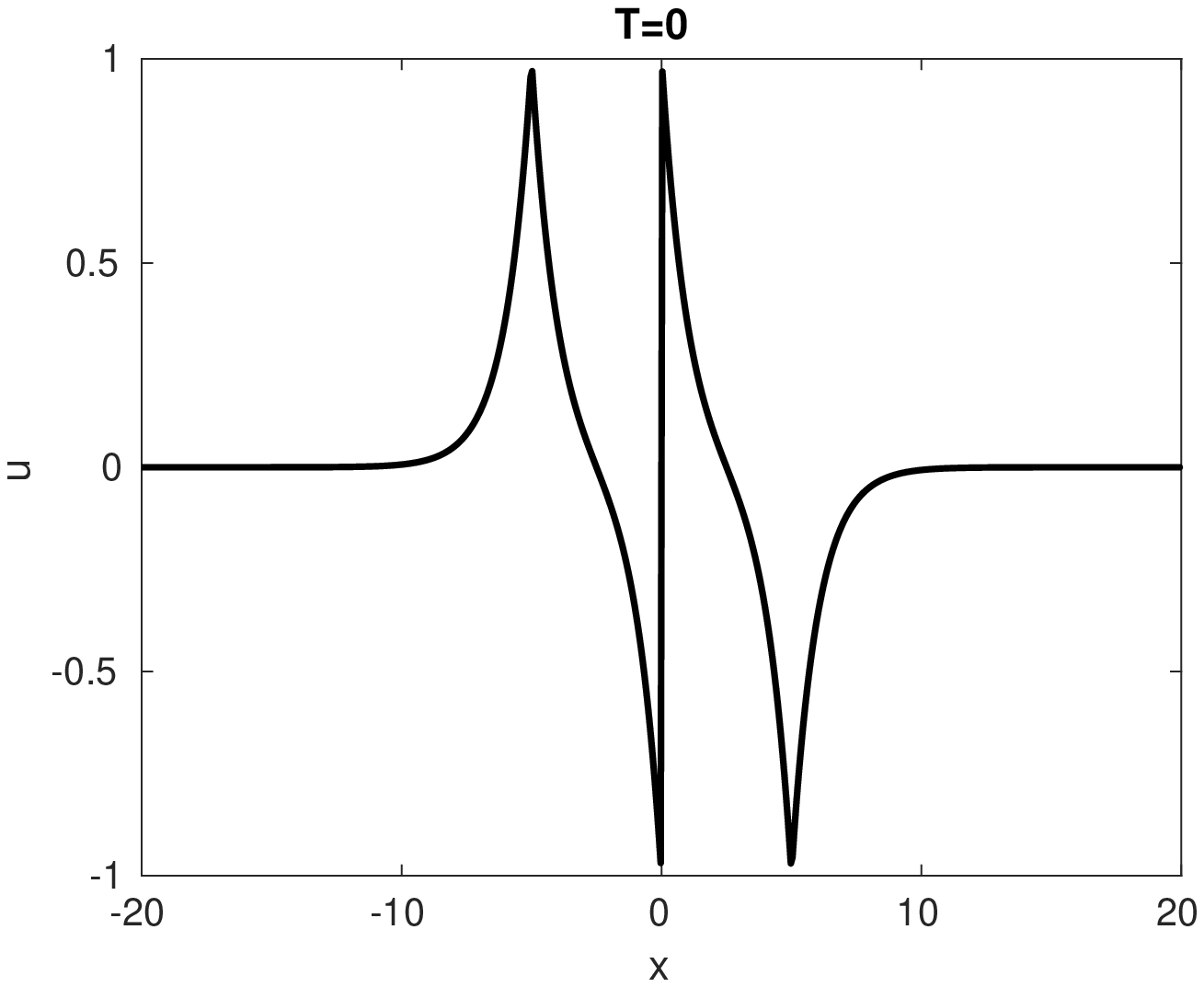} 
\includegraphics[width=0.45\textwidth]{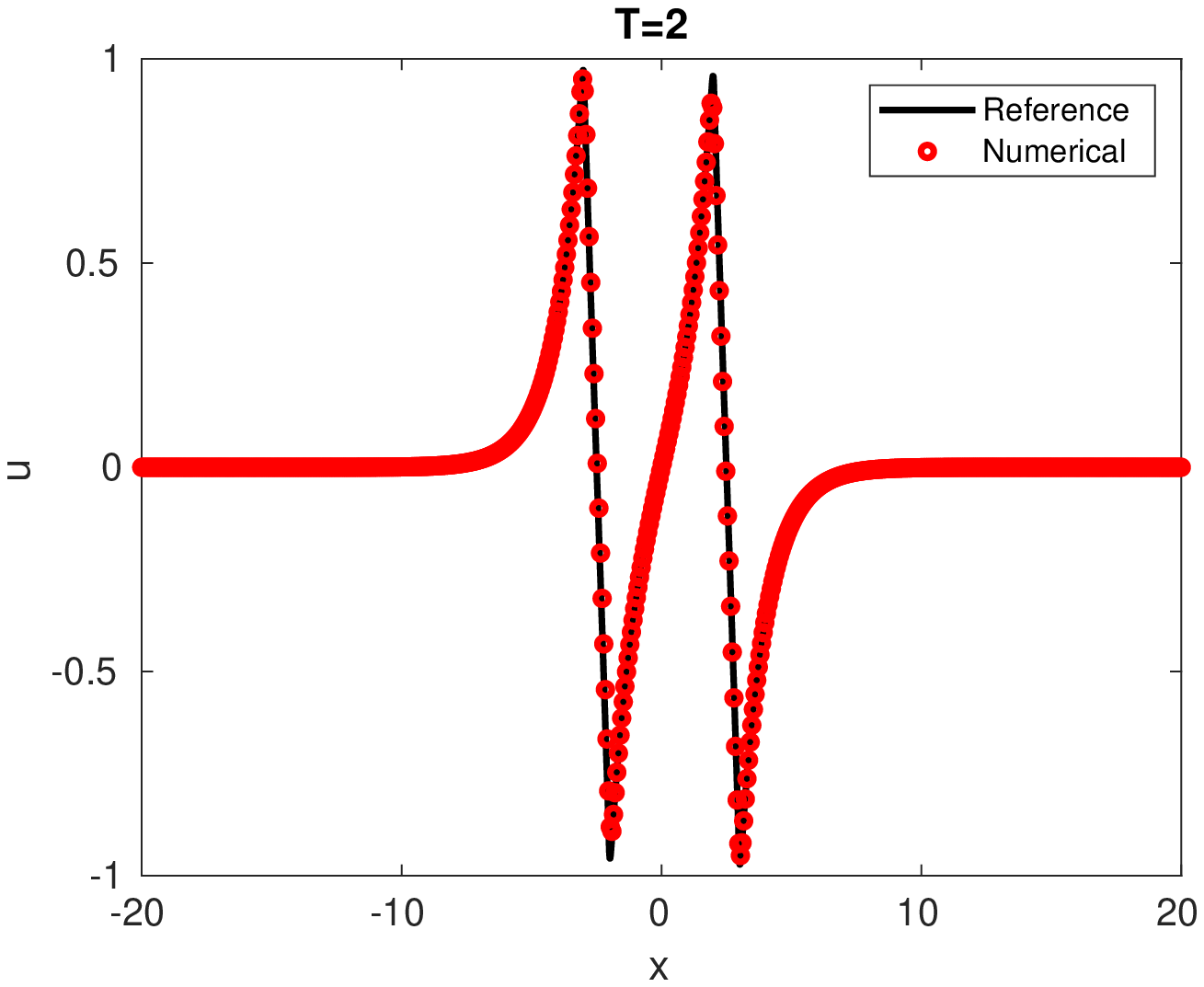} 
\includegraphics[width=0.45\textwidth]{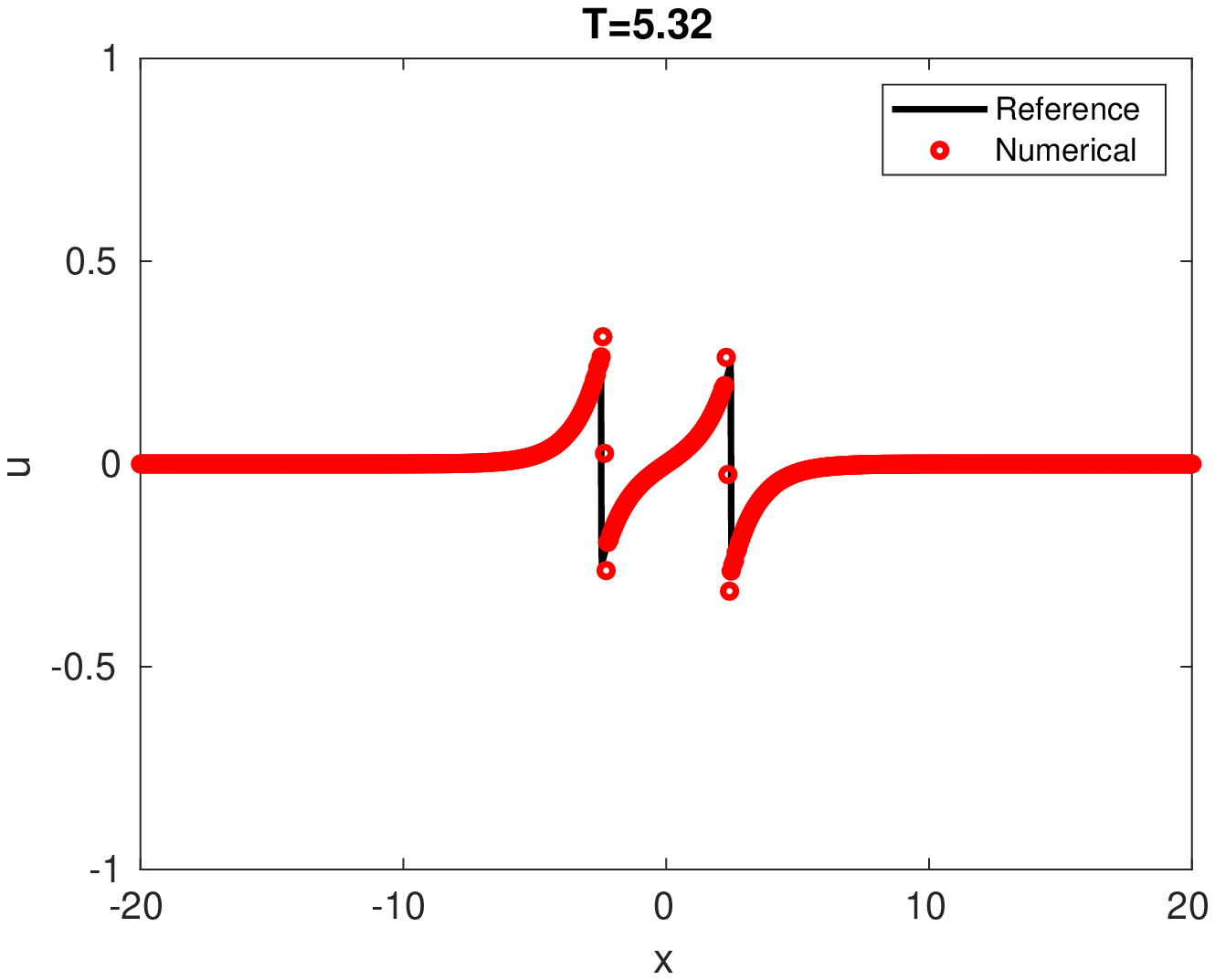} 
\includegraphics[width=0.45\textwidth]{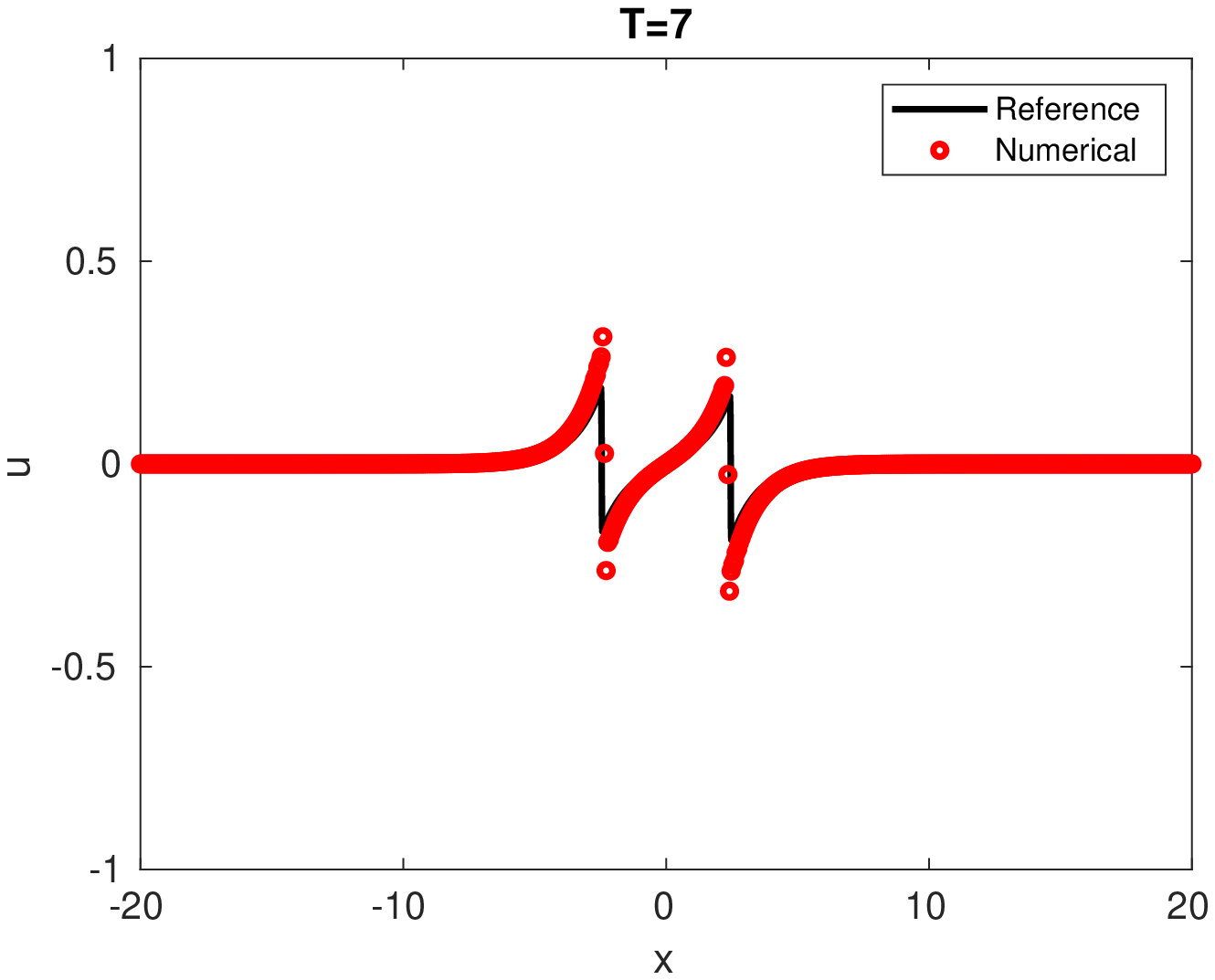} 
\caption{{\color{red}The triple interaction  of the DP in Example \ref{Ex:triple}. $N = 640$. WENO5.}}
\label{Fig:triple_NWENO5}
\end{figure} 
	
%
%
%
\end{exa}


\begin{exa}
\label{Ex: wavebreak}
{\bf Wave-breaking phenomena}\\
In this example, we consider the wave breaking phenomena of the DP equation, which was theoretically studied in \cite{liu_global_2006}. Briefly speaking, assume the initial condition 
$u(x,0)\in H^s(\mathbb{R}), s>\frac{1}{2}$, and there exists $x_0 \in \mathbb{R}$ such that the so-called momentum density, defined as $m_0(x)=u(x,0)-u_{xx}(x,0)$ changes  the sign from positive to negative at $x=x_0$, then the corresponding solution breaks in finite time $T_c<\infty$, i.e. the wave profile remains bounded but its slope becomes infinity  at time $T_c$. The shock waves usually appear afterwards. To verify these theoretical results, we consider two initial conditions, given by
\begin{align}
\label{init:wavebreak1}
u(x, 0) &= e^{0.5x^{2}}\sin(\pi x),\\
u(x, 0) &= \text{sech}^{2}(0.1(x+50)). \label{init:wavebreak2}
\end{align}
Figure \ref{Fig:shock_WENO5} shows the numerical solutions with the initial condition \eqref{init:wavebreak1} at  $T = 0, \ 0.18,\ 0.5$ and $1.1$ with $N = 640$  in the domain $[-2, 2]$.  {\color{red}The reference solutions are obtained by the fifth order classical WENO scheme \cite{xia_weighted_2017} with 2560 meshes.} The solution is smooth when $T<0.18$ and a shock is formed afterwards.
Figure \ref{Fig:wavebk_WENO7} show the numerical solutions with the initial condition \eqref{init:wavebreak2} at $T = 0,\ 10, \ 20$ and $30$ with $N = 2560$ in the domain $[-100, 100]$.  {\color{red} The reference solutions are obtained by the fifth order classical WENO scheme \cite{xia_weighted_2017} with 5120 meshes.}
The results with both initial settings agree well with those in \cite{Coclite.Karlsen.Risebro_JNA2008, Feng.Liu_JCP2009, xu_local_2011,xia_weighted_2017}.
\begin{figure}[!htbp]
\centering%
\includegraphics[width=0.45\textwidth]{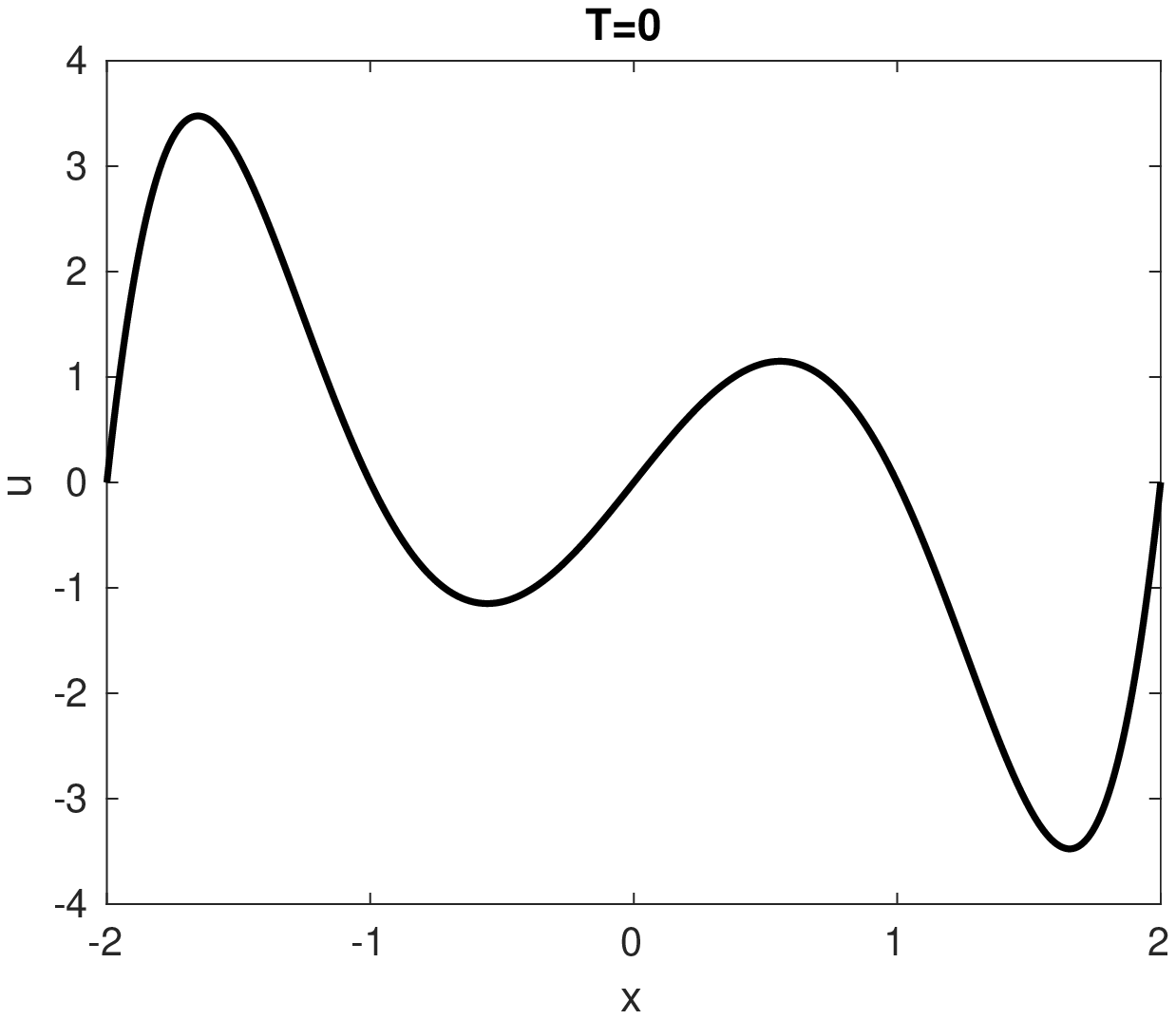} 
\includegraphics[width=0.45\textwidth]{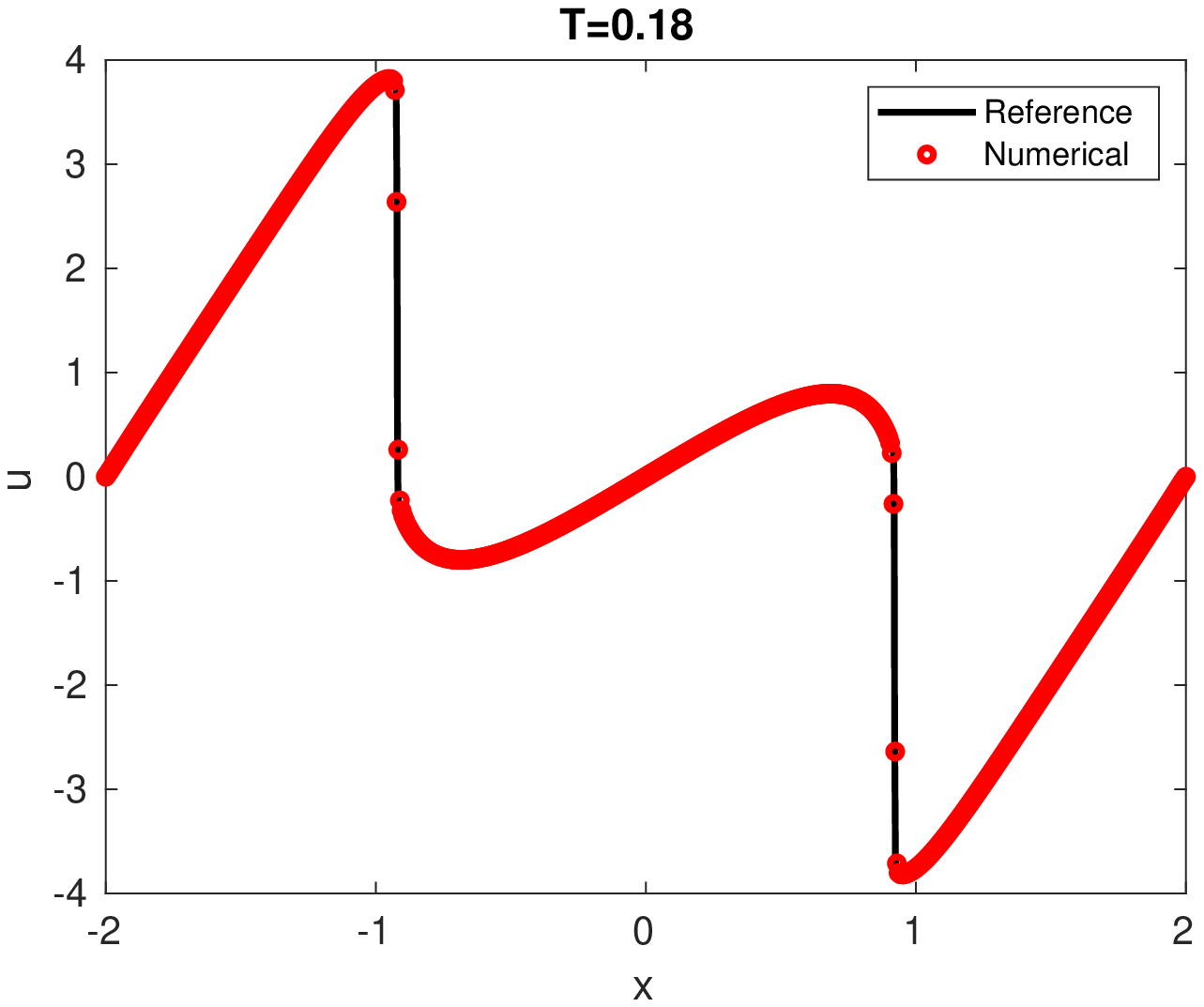}
\includegraphics[width=0.45\textwidth]{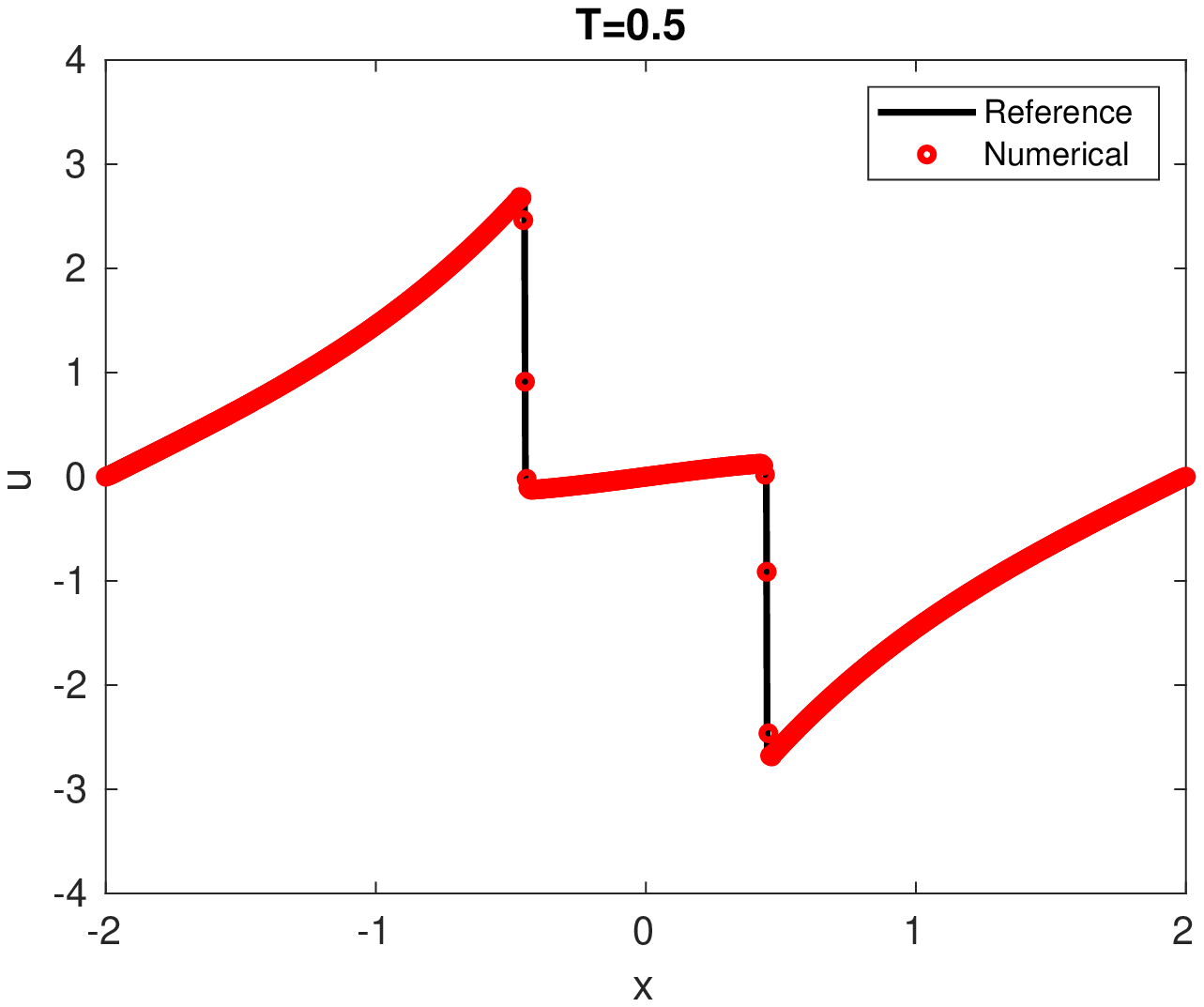} 
\includegraphics[width=0.45\textwidth]{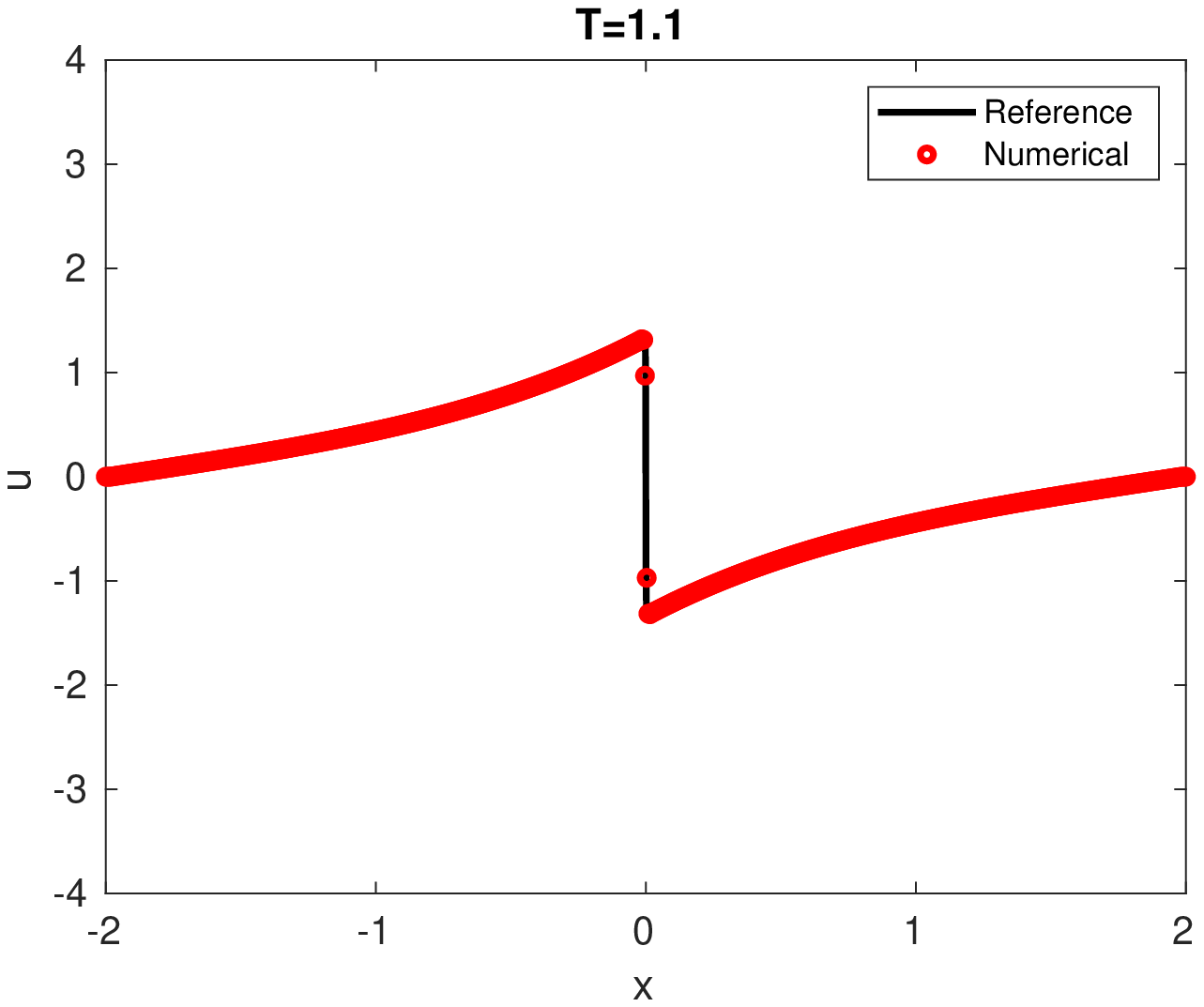}
\caption{{\color{red}Wave breaking of the DP equation in Example \ref{Ex: wavebreak} with initial condition \eqref{init:wavebreak1}. $N = 640$. MR-WENO5.}}
\label{Fig:shock_WENO5}
\end{figure} 

\begin{figure}[!htbp]
\centering%
\includegraphics[width=0.45\textwidth]{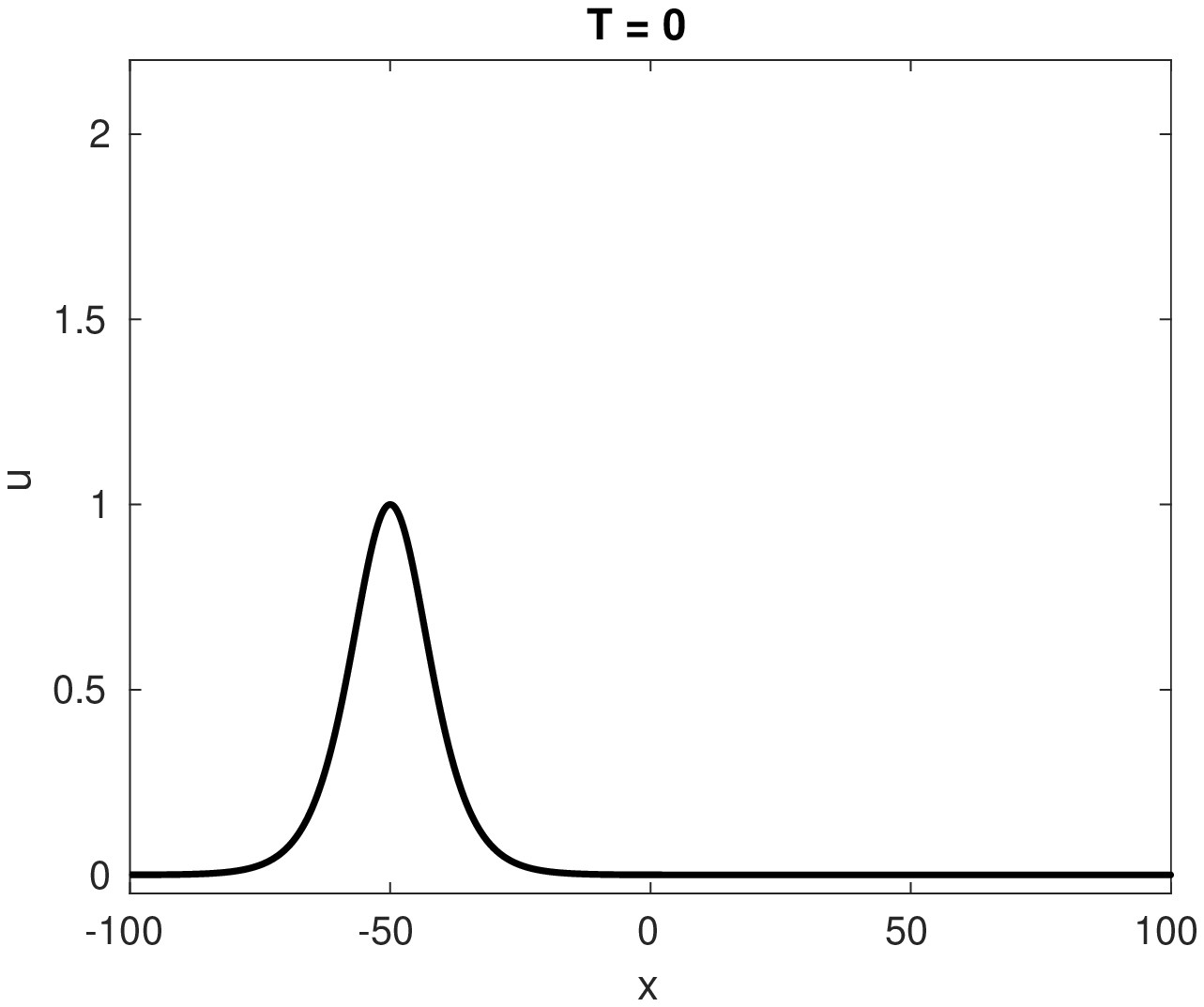}
\includegraphics[width=0.45\textwidth]{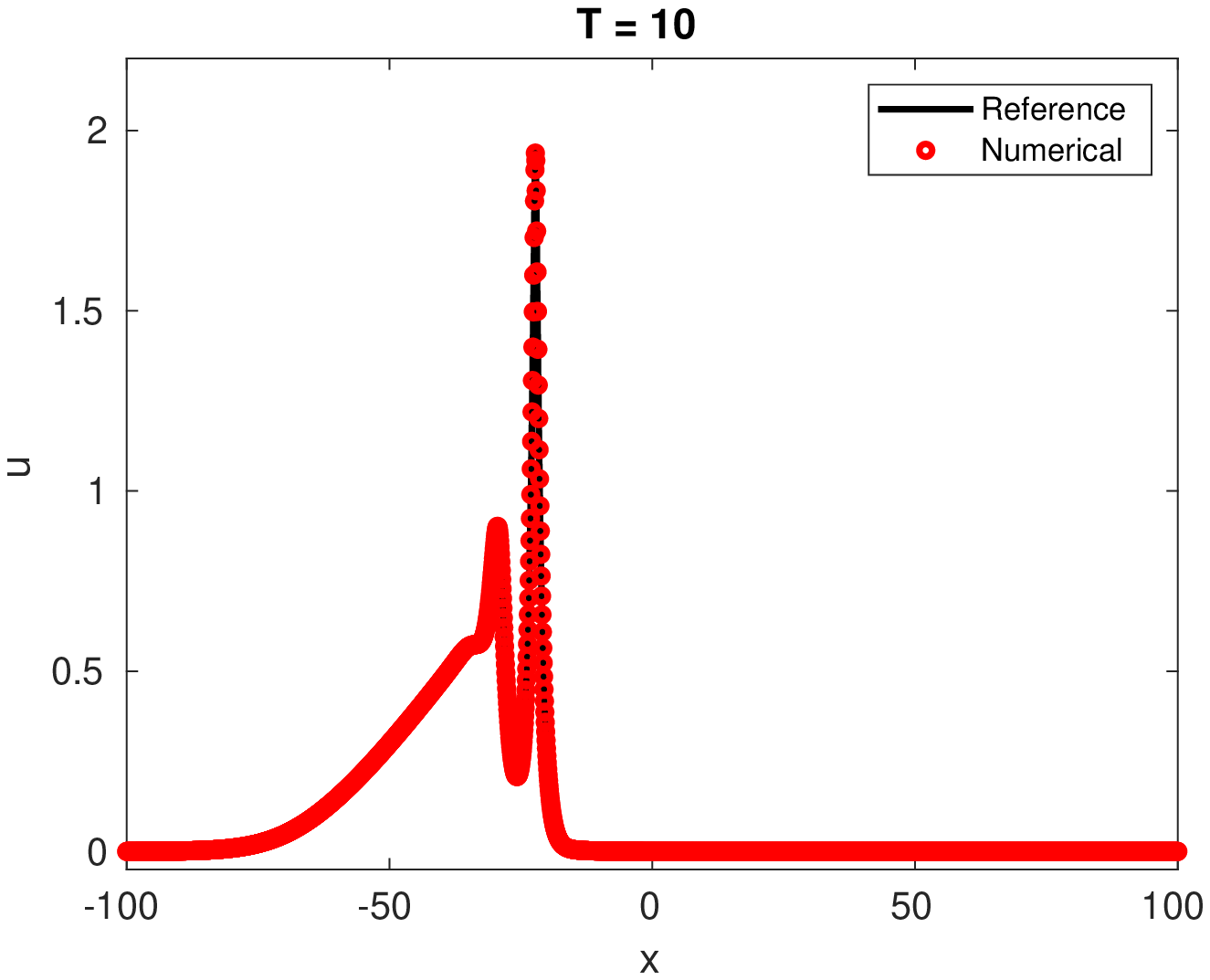}
\includegraphics[width=0.45\textwidth]{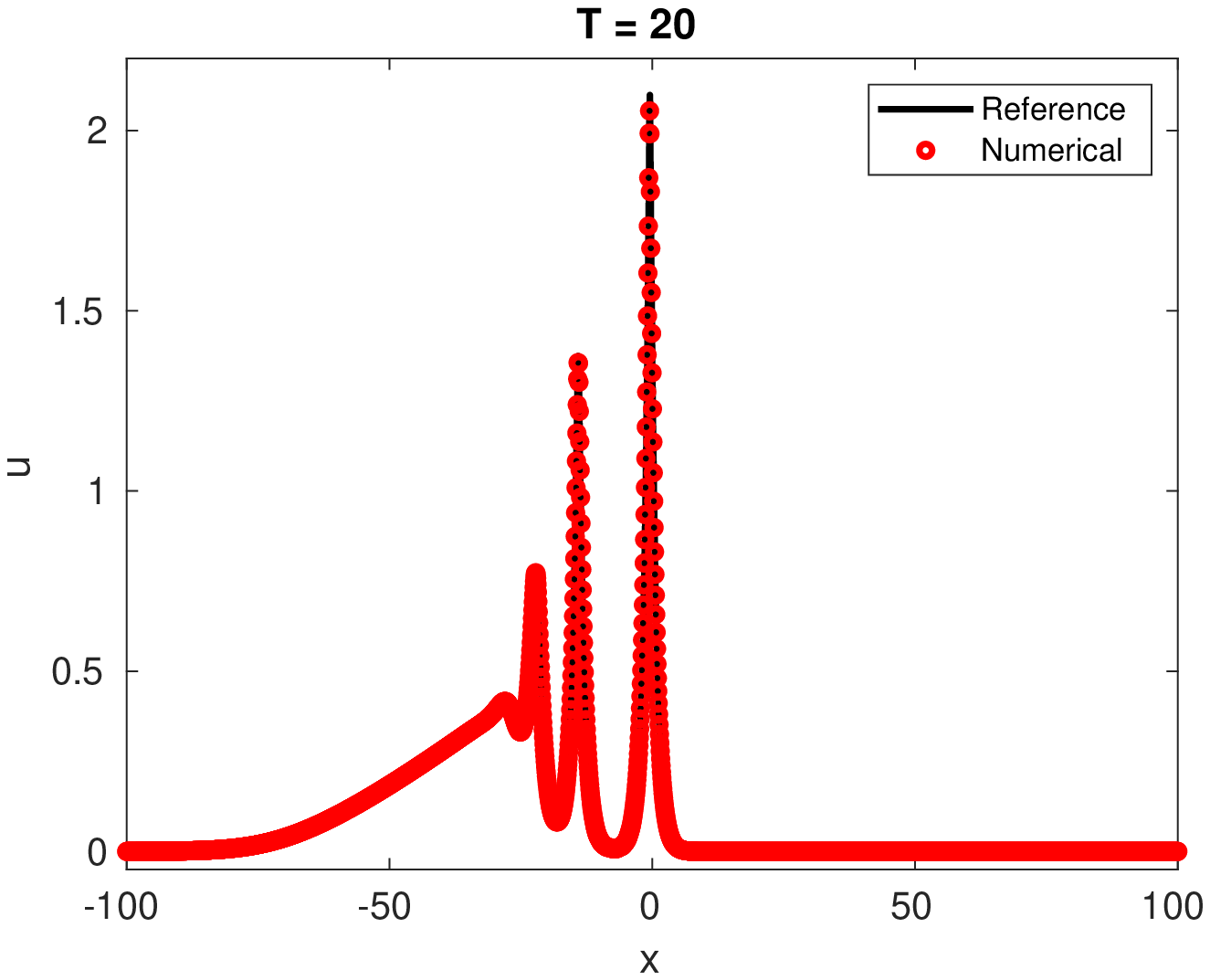}
\includegraphics[width=0.45\textwidth]{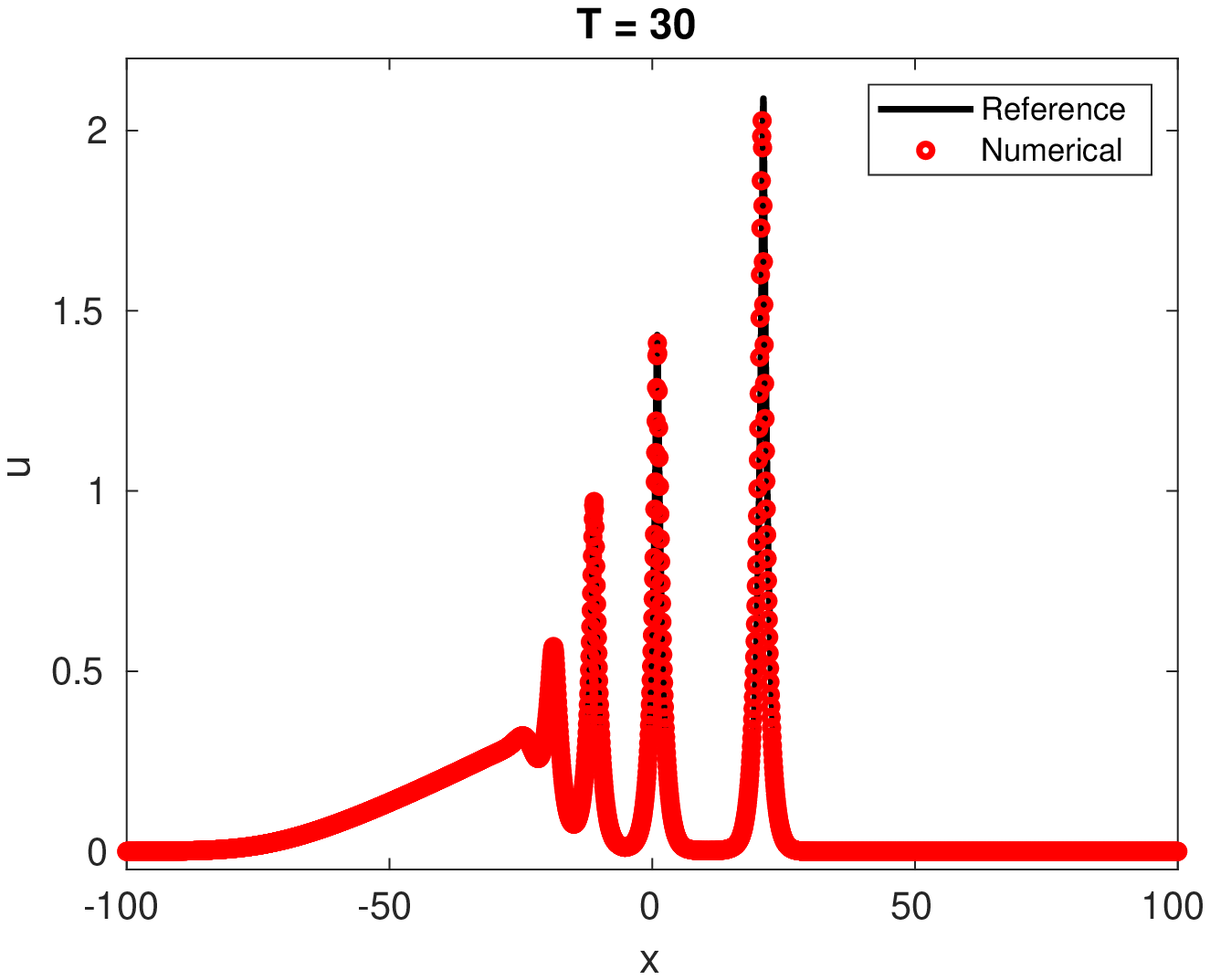}
\caption{{\color{red}Wave breaking of the DP equation in Example \ref{Ex: wavebreak} with initial condition \eqref{init:wavebreak2}. $N = 2560$. MR-WENO7.}}
\label{Fig:wavebk_WENO7}
\end{figure}
\end{exa}

\subsection{Numerical results for the $\mu$DP equation}
\label{numerical:mudp}
In this section, we present the numerical results to demonstrate the performance of WENO5, MR-WENO5 and MR-WENO7 for the \mdp equation with different initial conditions.

\begin{exa}
\label{Ex:mu-DP_twave_sm}
{\bf Accuracy test for smooth periodic waves}\\
In this example, we consider the $\mu$DP equation with  initial condition  $u(x, 0) = \phi(x)$, where $\phi$ is a solution of the following ODE
\begin{equation}\label{Eqn:CH_sm}
(\phi')^{2} = \frac{2\mu_{0}(M-\phi)(\phi-m)}{c-\phi},
\end{equation}
with  the constants $M$, $m$ and $c$ satisfying the condition $m<c<M$.
We remark here that this initial condition is taken as the one used for the $\mu$CH equation with smooth periodic waves in the form of $u(x,t)=\phi(x-ct)$, due to the lack of examples 
for the $\mu$DP equation with smooth traveling waves.
We set $M = 1.5$, $m=0.5$, and $c = 2$. This setting leads to a smooth periodic traveling wave of  the $\mu$CH equation with period $T_p=2.73321849515629$. $\mu_0$ can be further obtained with $\mu_0=2.55499933801271.$  An initial condition for $\phi$ with
$\phi(0.796433828683979) = 1$ can also be computed by setting $\theta=\pi/2$ in (6.10) in \cite{Lenells.Misiolek_CMP2010}. More details can be found in \cite{Lenells.Misiolek_CMP2010, xu_local_2011}. 
Then we can get a high-precision numerical solution of \eqref{Eqn:CH_sm} in the domain $[-T_p/2,T_p/2]$ by  a fourth-order RK method with $2^{21}=2097152$ cells. We compute the error between the numerical solution with $N$ cells and with $2N$ cells to test the order of accuracy.
We list the $L^{1}$ and $L^{\infty}$ errors and  orders of accuracy with WENO5, MR-WENO5 and MR-WENO7 at $T=0.1$ in Table \ref{Tab:mu-DP_twave_sm}. 
We can see that all three methods achieve the desired order of accuracy, i.e. fifth-order accuracy for WENO5 and MR-WENO5 and seventh-order accuracy for MR-WENO7.
\begin{table}[!htbp]
\centering
\caption{The \mdp equation with sufficiently smooth solution in Example \ref{Ex:mu-DP_twave_sm} at $T=0.1$.} 
\label{Tab:mu-DP_twave_sm}
\smallskip		
\begin{tabular}{ccccccccc}
\toprule
& \multicolumn{4}{c}{WENO5} 	& \multicolumn{4}{c}{MR-WENO5}  \\ \cmidrule(lr){2-5}\cmidrule(lr){6-9}
$N$& $L^{1}$ \mbox{error} & \mbox{Order} & $L^{\infty}$ \mbox{error} & \mbox{Order} & $L^{1}$ \mbox{error} & \mbox{Order} & $L^{\infty}$ \mbox{error} & \mbox{Order}\\ \hline
32  & 2.56E-05 &      & 2.77E-04 &       & 1.34E-04  &       & 1.24E-03  &       \\
64  & 7.08E-07 & 5.17 & 7.19E-06 & 5.27  & 6.58E-07  & 7.67  & 6.69E-06  & 7.54  \\ 
128 & 2.07E-08 & 5.09 & 3.24E-07 & 4.47  & 2.00E-08  & 5.03  & 2.42E-07  & 4.79  \\ 
256 & 6.26E-10 & 5.05 & 6.20E-09 & 5.70  & 6.26E-10  & 5.00  & 6.22E-09  & 5.28  \\
512 & 1.92E-11 & 5.02 & 1.96E-10 & 4.98  & 1.92E-11  & 5.02  & 1.96E-10  & 4.99  \\
\midrule
&\multicolumn{8}{c}{MR-WENO7}  \\ \cmidrule(lr){2-9}
$N$ & &$L^{1}$ \mbox{error} && \mbox{Order} & &$L^{\infty}$ \mbox{error} && \mbox{Order} \\ \hline
32  & & 5.48E-06  & &       & & 3.67E-05 & &   \\ 
64  & & 5.32E-08  & & 6.69  & & 5.92E-07 & & 5.95   \\ 
128 & & 4.52E-10  & & 6.88  & & 5.29E-09 & & 6.81   \\ 
256 & & 3.85E-12  & & 6.88  & & 4.16E-11 & & 6.99 \\
\bottomrule    
\end{tabular}
\end{table}
\end{exa}

\begin{exa}
\label{Ex:mudp_peakon}
{\bf Peakon solutions}\\
In this example, we consider the \mdp equation with $M$-peakon solutions \cite{Lenells.Misiolek_CMP2010} in the form of
\begin{equation}\label{mu-DP_peakon}
u(x, t) = \sum\limits^{M}_{i=1}\psi_{i}(t)g(x-\varphi_i(t)),
\end{equation}
where $g(x)$ is the Green's function given by
\begin{equation}
\label{eq:g}
g(x) = \frac{1}{2}x(x-1)+\frac{13}{12},~~ x\in[0, 1),
\end{equation}
and is extended periodically to $\mathbb{R}$, namely
\begin{equation}
g(x-x') = \frac{(x-x')^{2}}{2}-\frac{|x-x'|}{2}+\frac{13}{12}, ~~x\in[x',x'+1),
\end{equation}
where $x'$ denotes a translation of one periodic interval. The time-dependent variables $\psi_{i}(t)$ and $\varphi_i(t)$ satisfy the following ODE
\begin{align}
\label{eq:variable}
\frac{\mathrm{d}\varphi_i}{dt}	 =\sum\limits^{M}_{j=1}\psi_{j}g(\varphi_i-\varphi_j), \quad\quad
\frac{\mathrm{d}\psi_i}{dt}=-2\sum\limits^{M}_{j=1}\psi_{i}\psi_{j}g'(\varphi_i-\varphi_j),
\end{align}
where $g'(x)$ is the derivative of $g(x)$  in \eqref{eq:g} with the value 0 assigned
to the otherwise undetermined derivative. That is,
\begin{equation}\label{green_der}
g'(x) := \begin{cases}
0, & x= 0, \\
x-\frac{1}{2}, &0<x<1.
\end{cases}
\end{equation}
%

Now we simulate the $\mu$DP equation at  $T =0,\  1, \ 5$ and $10$ with $N=160$, under the following initial condition settings:
\begin{itemize}
\item One peakon
\begin{equation}
\psi_{1}(0)=0.333, ~~\varphi_1(0) = -0.5;
\end{equation}
\item Two peakons
\begin{equation}
\begin{split}
\psi_{1}(0) &=0.1, ~~\varphi_1(0) = 0.4, \\
\psi_{2}(0) &=0.08, ~~\varphi_2(0) = 0.1.
\end{split}
\end{equation}
\end{itemize}
	
The solutions of the \mdp equation with one peakon are shown in Figure \ref{Fig:mu-DP_weno5-zq_p1} and with two peakons are shown in Figure \ref{Fig:mu-DP_mrweno5_p2}.
We can see clearly that the moving peakon profile are well resolved. There is no numerical oscillation near the wave crest. The results match very well with the exact solution and agree well with \cite{zhang_local_2019, zhao_high_2020}.
\begin{figure}[!htbp]
\centering%
\includegraphics[width=0.45\textwidth]{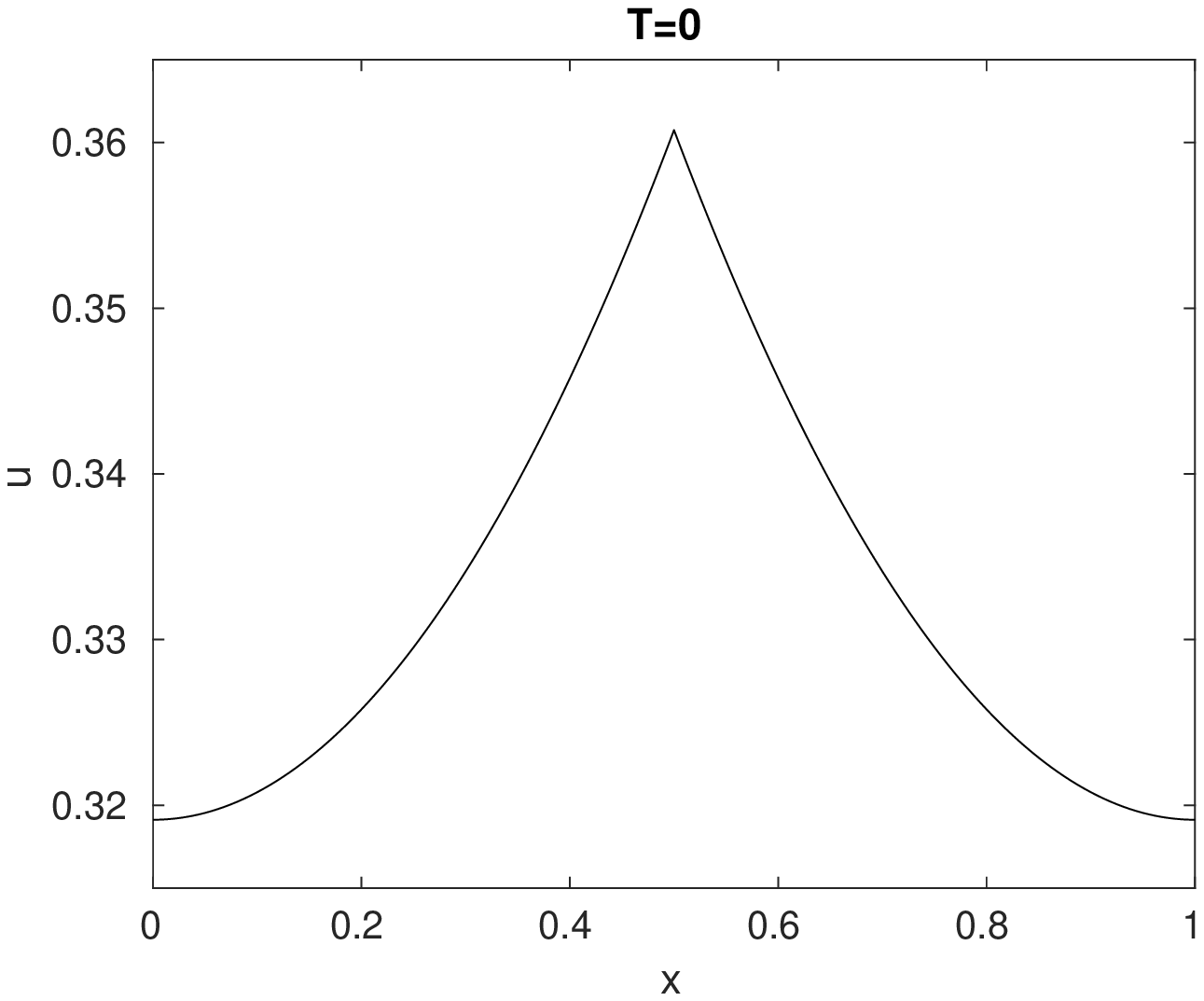}
\includegraphics[width=0.45\textwidth]{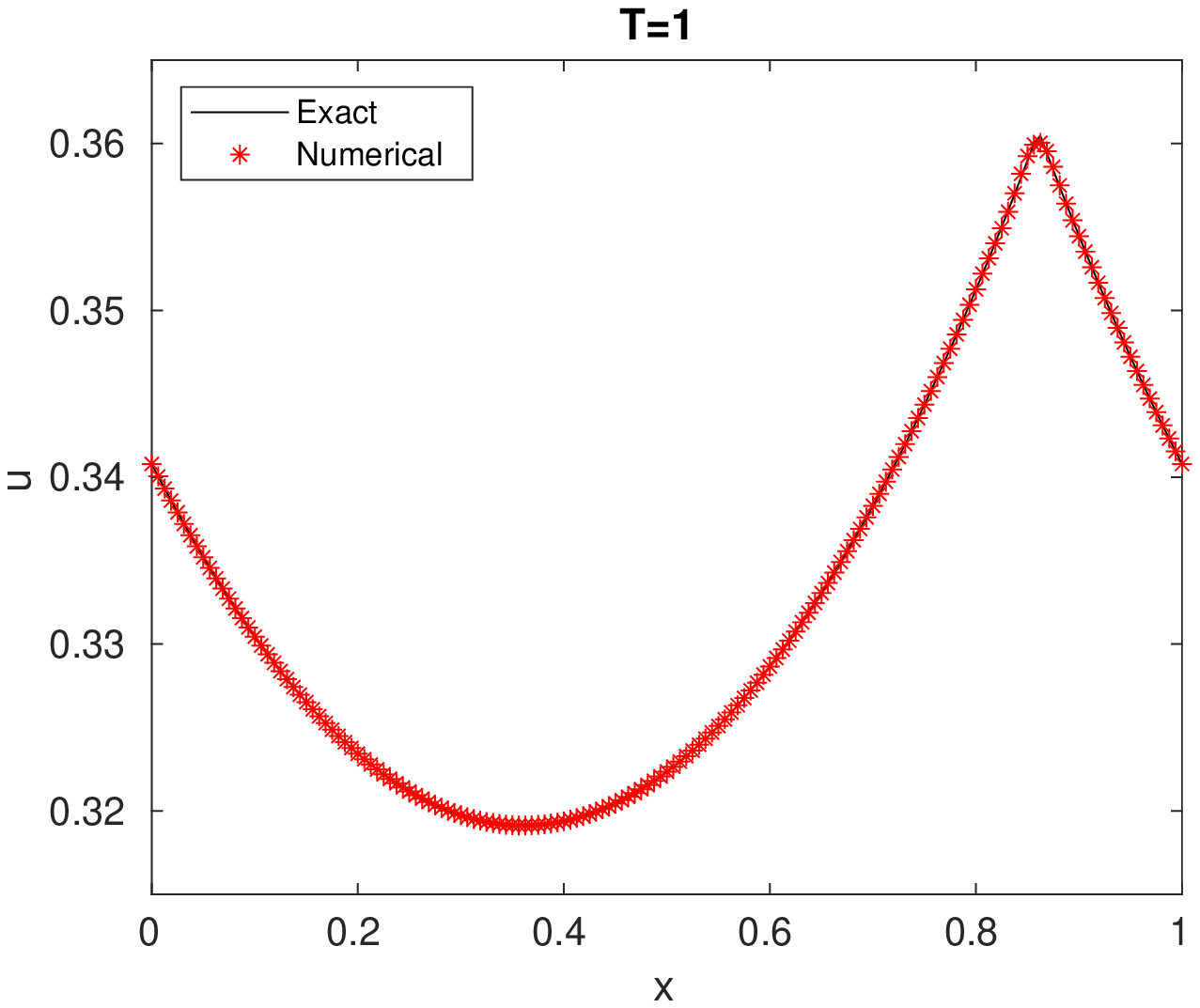}
\includegraphics[width=0.45\textwidth]{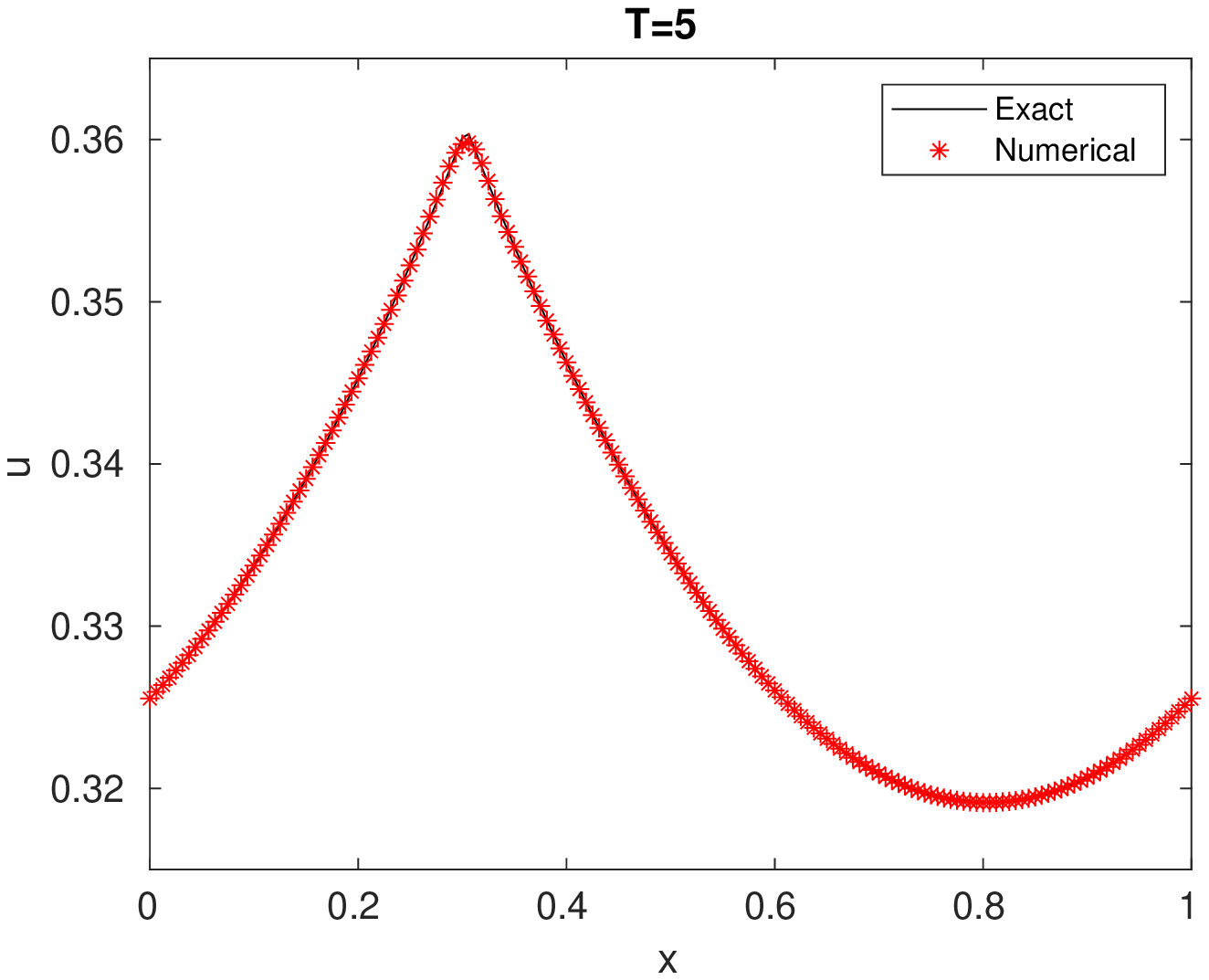} 
\includegraphics[width=0.45\textwidth]{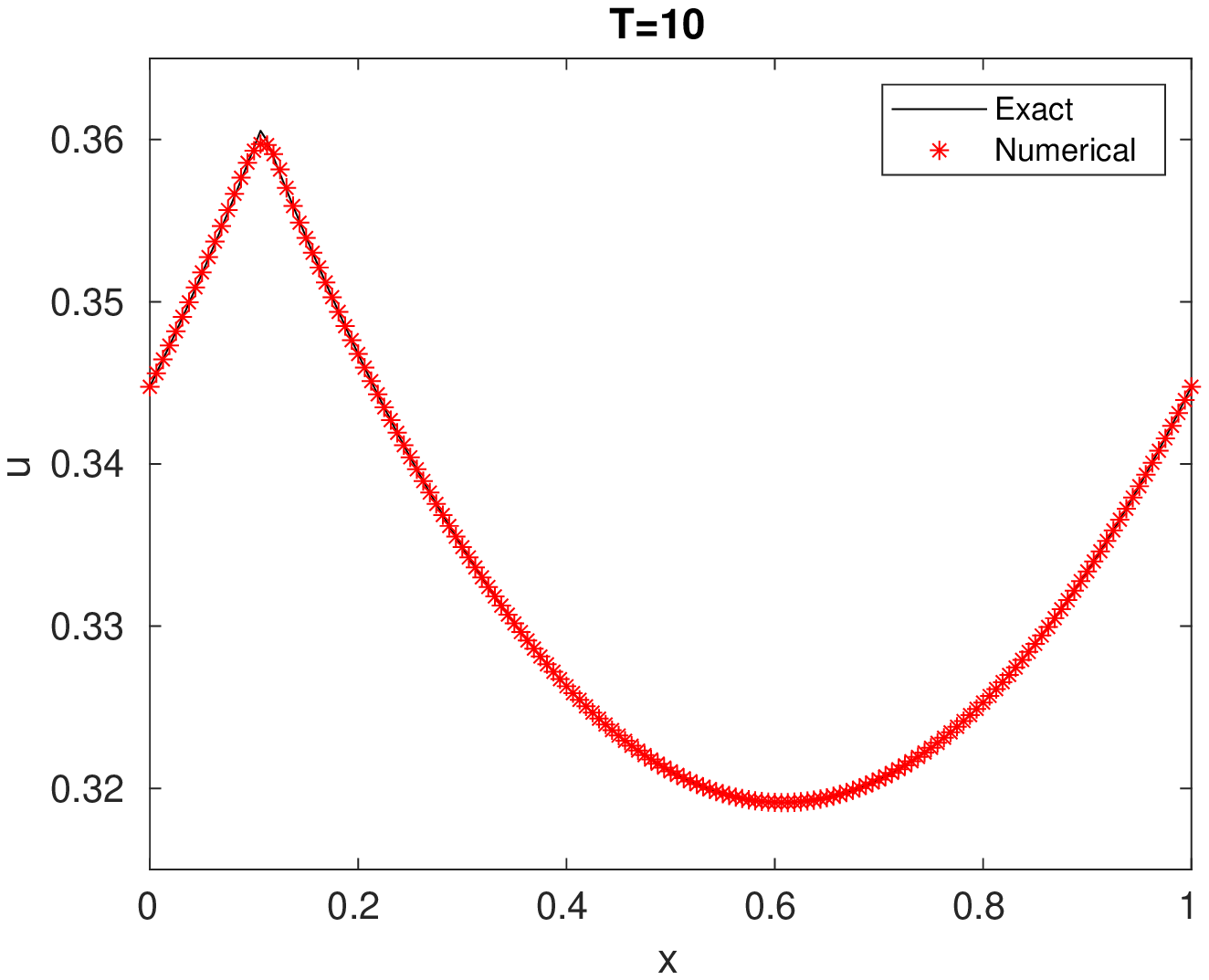}
\caption{One-peakon solution of the $\mu$DP equation in Example \ref{Ex:mudp_peakon}. $N = 160$. WENO5.}
\label{Fig:mu-DP_weno5-zq_p1}
\end{figure}

\begin{figure}[!htbp]
\centering%
\includegraphics[width=0.45\textwidth]{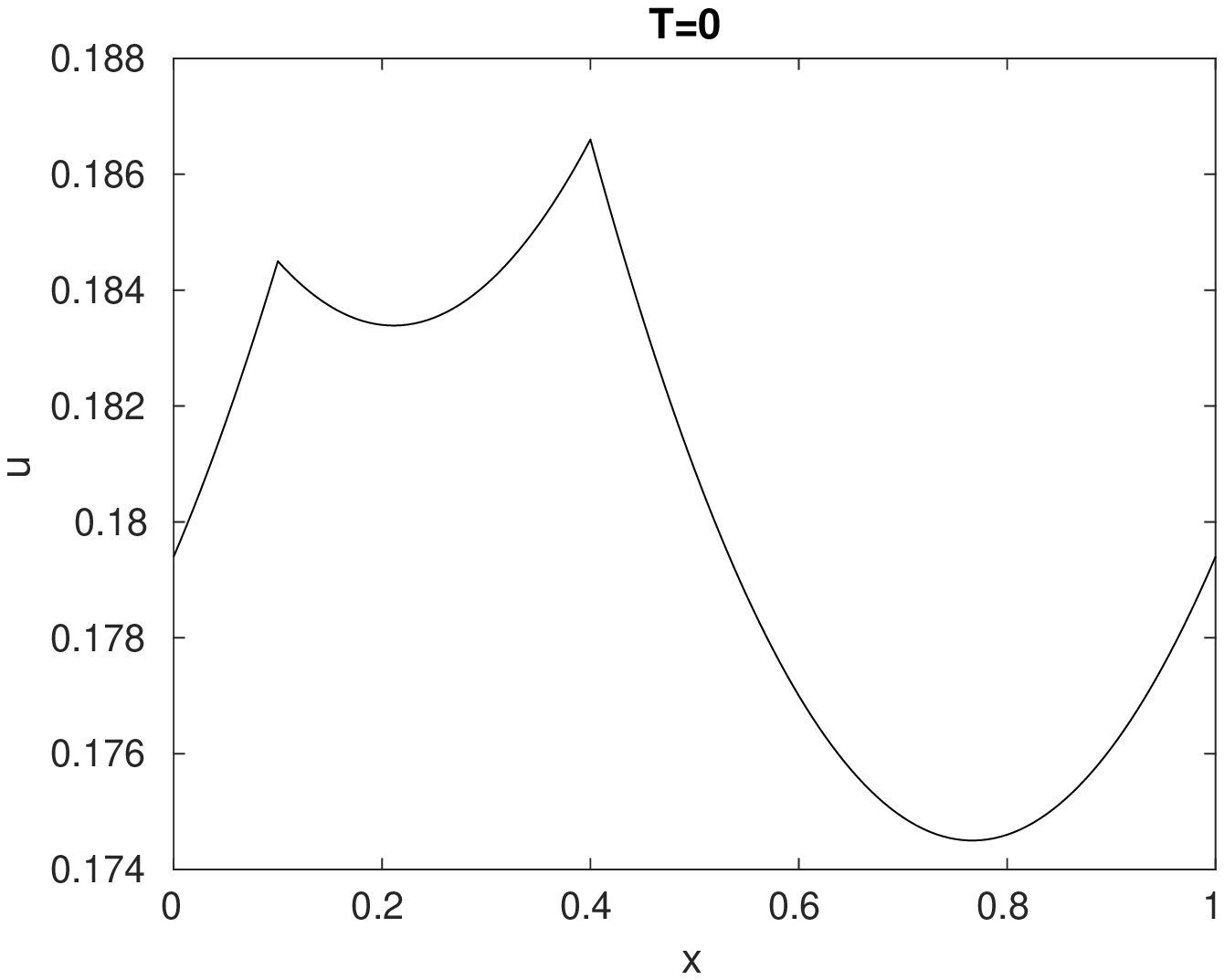} 
\includegraphics[width=0.45\textwidth]{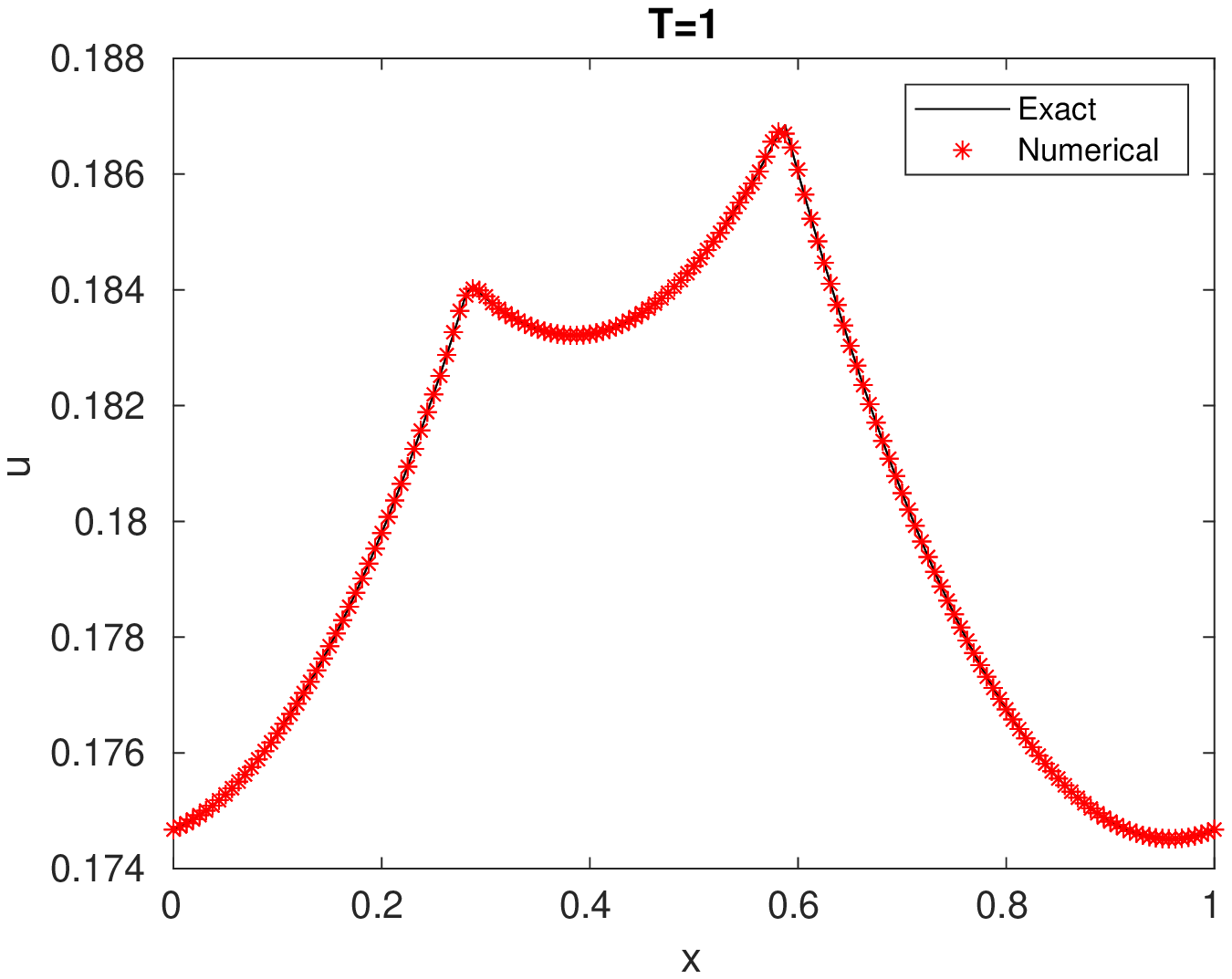}
\includegraphics[width=0.45\textwidth]{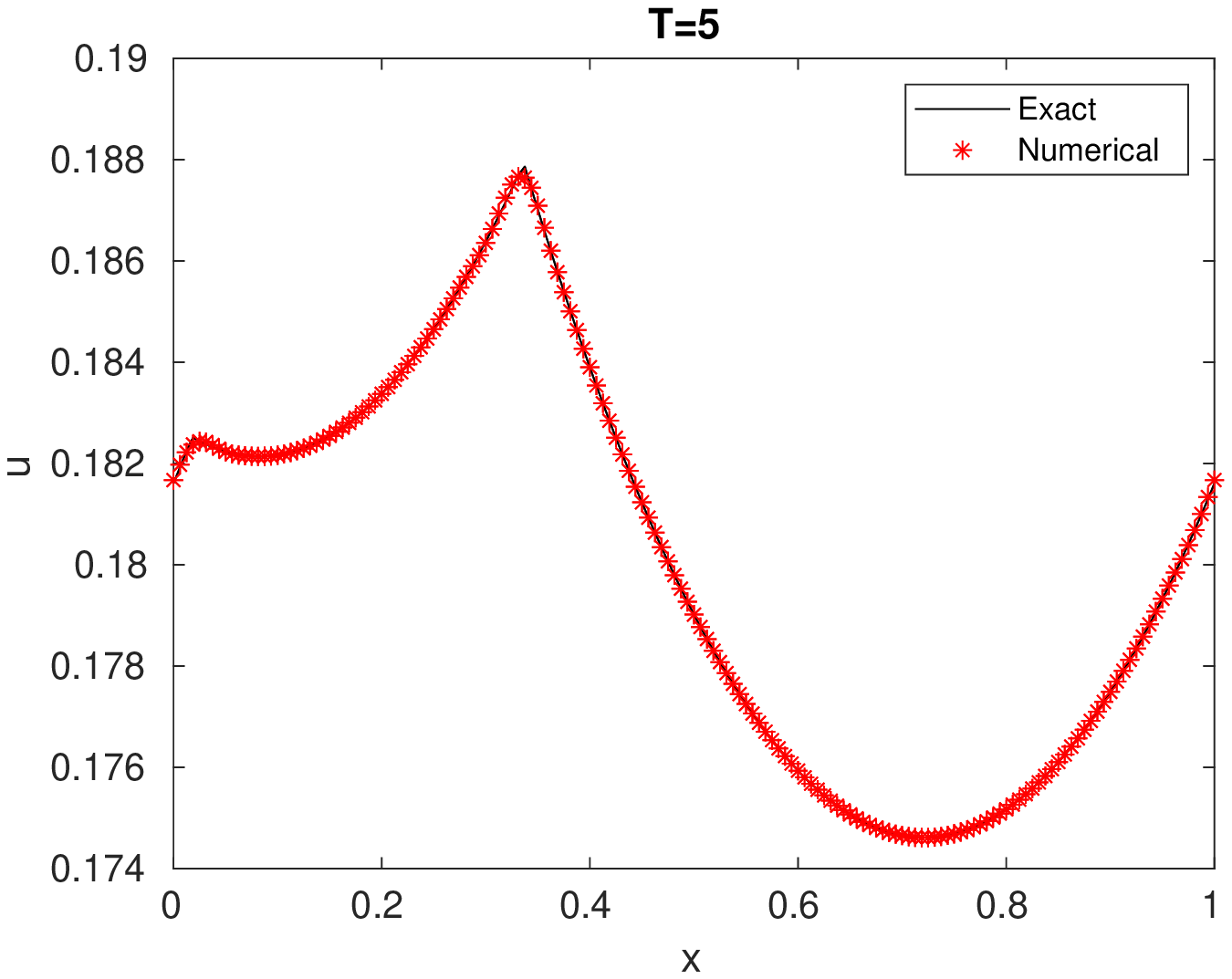} 
\includegraphics[width=0.45\textwidth]{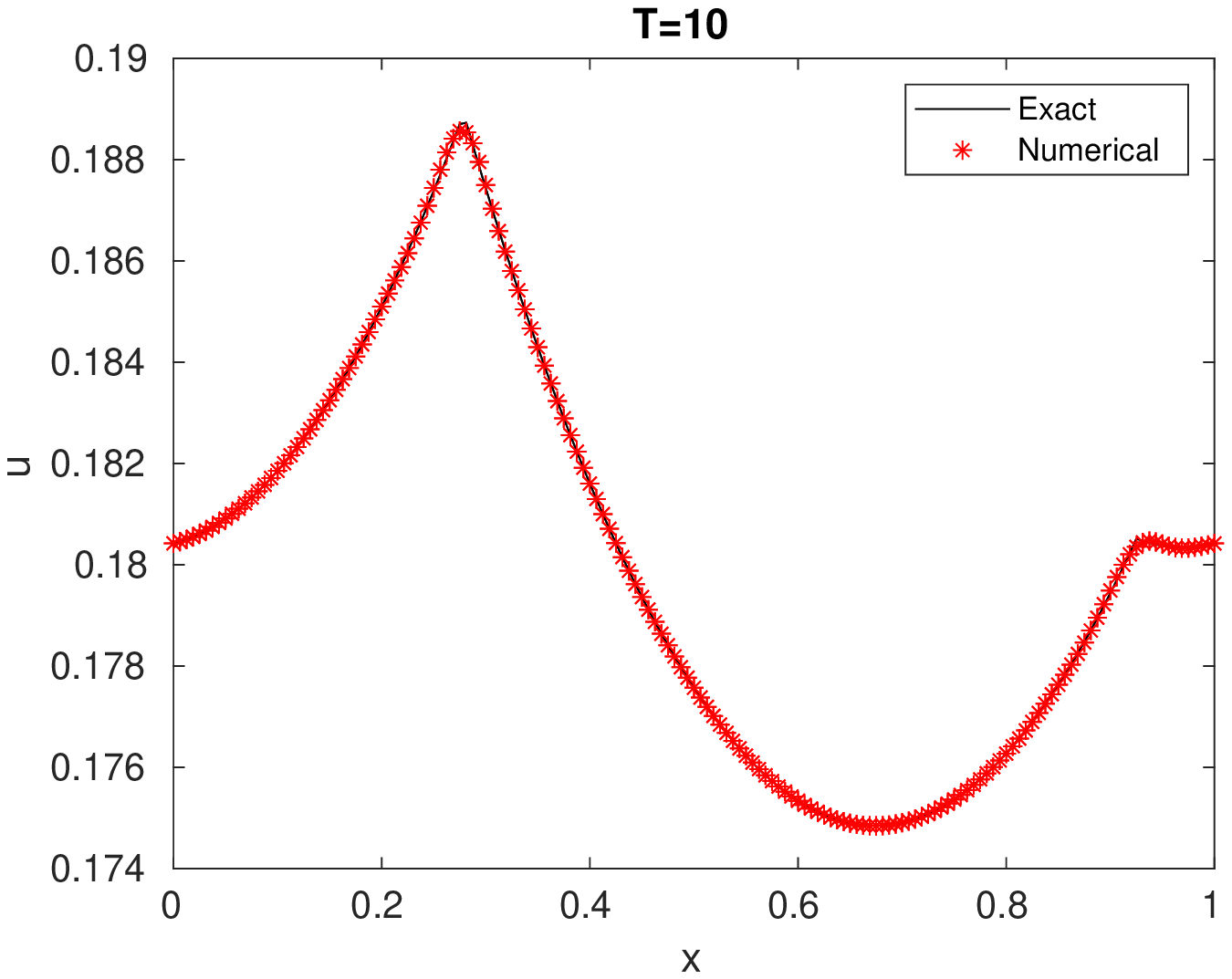}
\caption{Two-peakon solution of the $\mu$DP equation in Example \ref{Ex:mudp_peakon}. $N = 160$. MR-WENO5.}
\label{Fig:mu-DP_mrweno5_p2}
\end{figure}

%
%
\end{exa}

\begin{exa}
\label{Ex:shock_mudp}	
{\bf Shock solutions}\\
In this example, we consider the \mdp equation with $M$-shock solutions \cite{Lenells.Misiolek_CMP2010} in the form of
\begin{equation}\label{mu-DP_shock}
u = \sum\limits^{M}_{i=1}\left(\psi_{i}g(x-\varphi_i)+s_{i}g'(x-\varphi_i)\right),
\end{equation}
with $g(x)$  defined by \eqref{eq:g} and $g'(x)$  defined by \eqref{green_der}. 
The time-dependent variables $\psi_{i}(t),\ \varphi_i(t)$  and $s_i(t)$ satisfy the following ODE
\begin{align}
\label{eq:variable2}
\begin{split}
\frac{\mathrm{d}\varphi_i}{dt}	 &=\sum\limits^{M}_{j=1}\left(\psi_{j}g(\varphi_i-\varphi_j)+s_ig'(\varphi_i-\varphi_j)\right), \\
\frac{\mathrm{d}\psi_i}{dt}&=2\sum\limits^{M}_{j=1}\left(\psi_j-\psi_{i}\psi_{j}g'(\varphi_i-\varphi_j)\right),\\
\frac{\mathrm{d}s_i}{dt}&=-\sum\limits^{M}_{j=1}s_{i}\psi_{j}g'(\varphi_i-\varphi_j).
\end{split}
\end{align}
Now we simulate the $\mu$DP equation at $T =0,\  1, \ 3$ and $5$ with $N=320$, under the following initial condition settings:
\begin{itemize}
\item One shock
\begin{equation}
\psi_{1}(0) = 0.333,~~\varphi_1(0)=0.1, ~~s_{1} = 0.1;
\end{equation}
\item Two shocks
\begin{equation}
\begin{split}
	\psi_{1}(0) &= 0.3,~~ \varphi_{1}(0) = 0.2, ~~s_{1}(0) = 0.4, \\
	\psi_{2}(0) &= 0.1, ~~\varphi_{2}(0) = 0.5, ~~s_{2}(0) = 0.2.
\end{split}
\end{equation}
\end{itemize}	

Due to the complicated feature of the \mdp equation, the solution with shocks are  sensitive to the choice of linear weights. To get better non-oscillatory solutions, we set linear weights as
$\gamma_{1} = 0.4, \;\gamma_{2} = 0.3, \; \gamma_{3} = 0.3$,
for WENO5, $\gamma_{2, 1} = 1/11,   \gamma_{2, 2} = 10/11,  
\gamma_{3, 1} = 0.666, \gamma_{3, 2} = 0.001, \gamma_{3, 3} = 0.333$
for MR-WENO5, and $\gamma_{2, 1} = 1/11, \gamma_{2, 2} = 10/11,  
\gamma_{3, 1} = 1/111, \gamma_{3, 2} = 10/111, \gamma_{3, 3} = 100/111,   
\gamma_{4, 1} = 0.665,  \gamma_{4, 2} = 0.001,  \gamma_{4, 3} = 0.001,  \gamma_{4, 4} = 0.333$,
for MR-WENO7.
\begin{figure}[!htbp]
	\centering%
	\mbox{
		\includegraphics[width=0.45\textwidth]{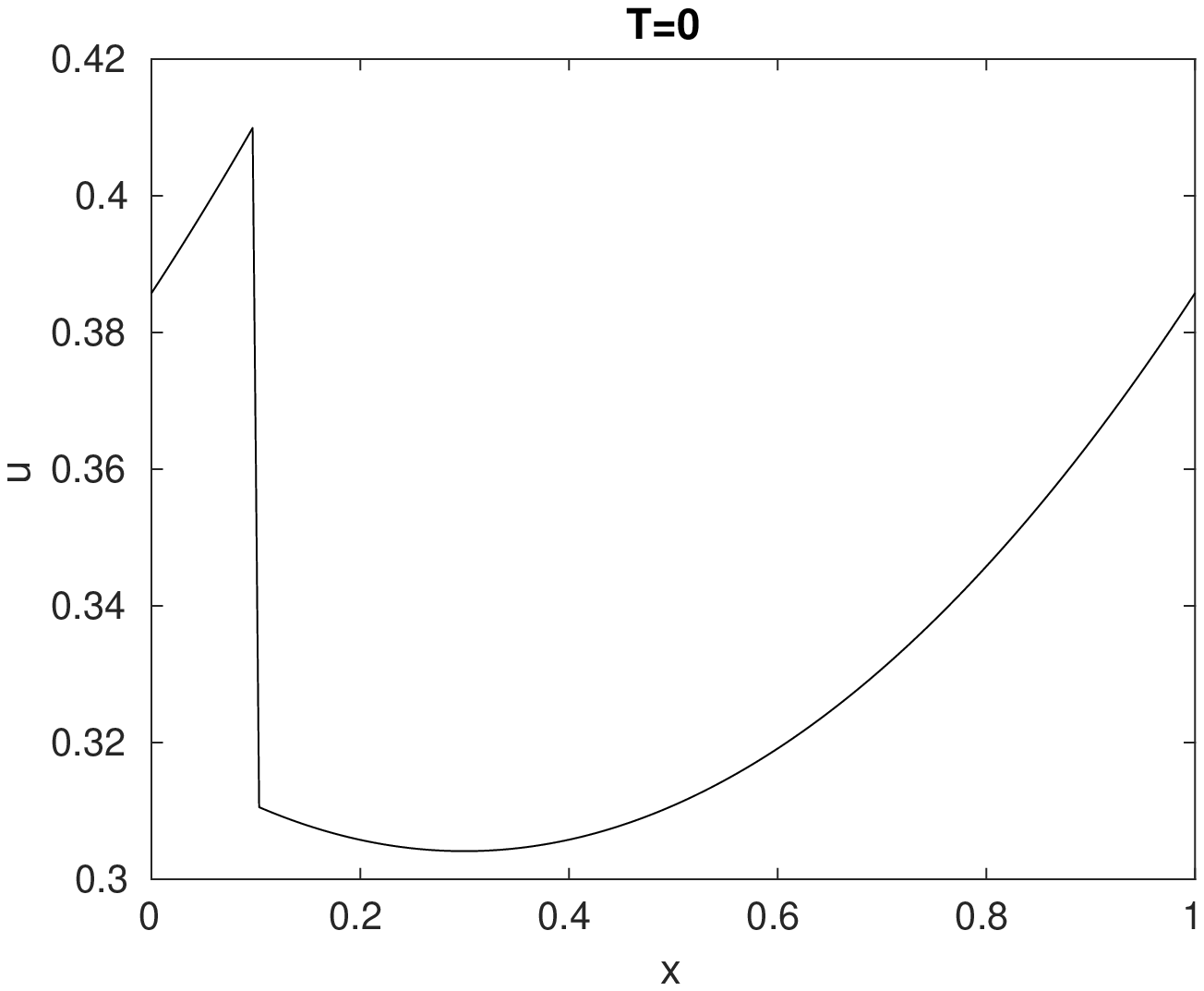} 
		\quad 
		\includegraphics[width=0.45\textwidth]{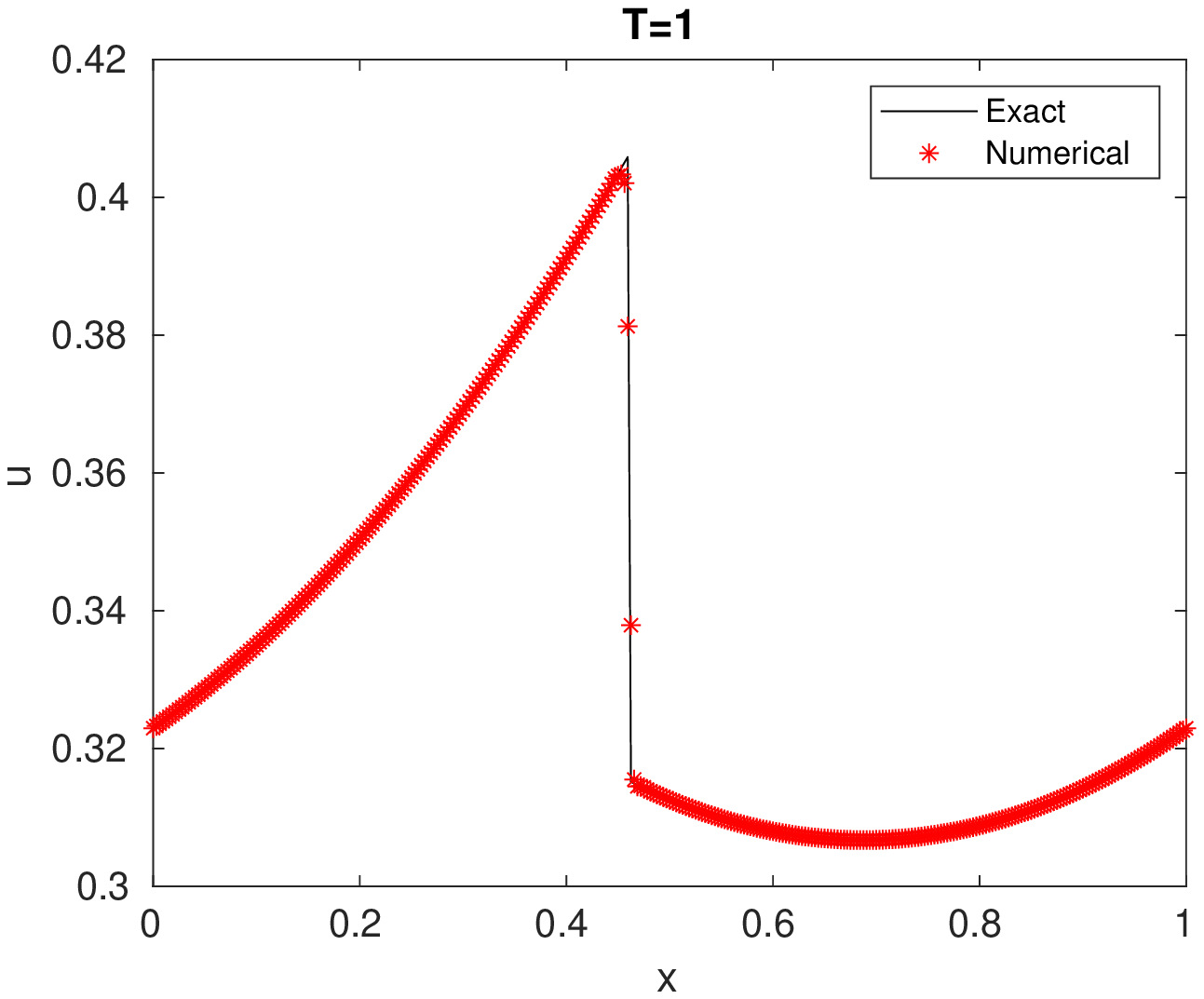}}
	\quad
	\mbox{
		\includegraphics[width=0.45\textwidth]{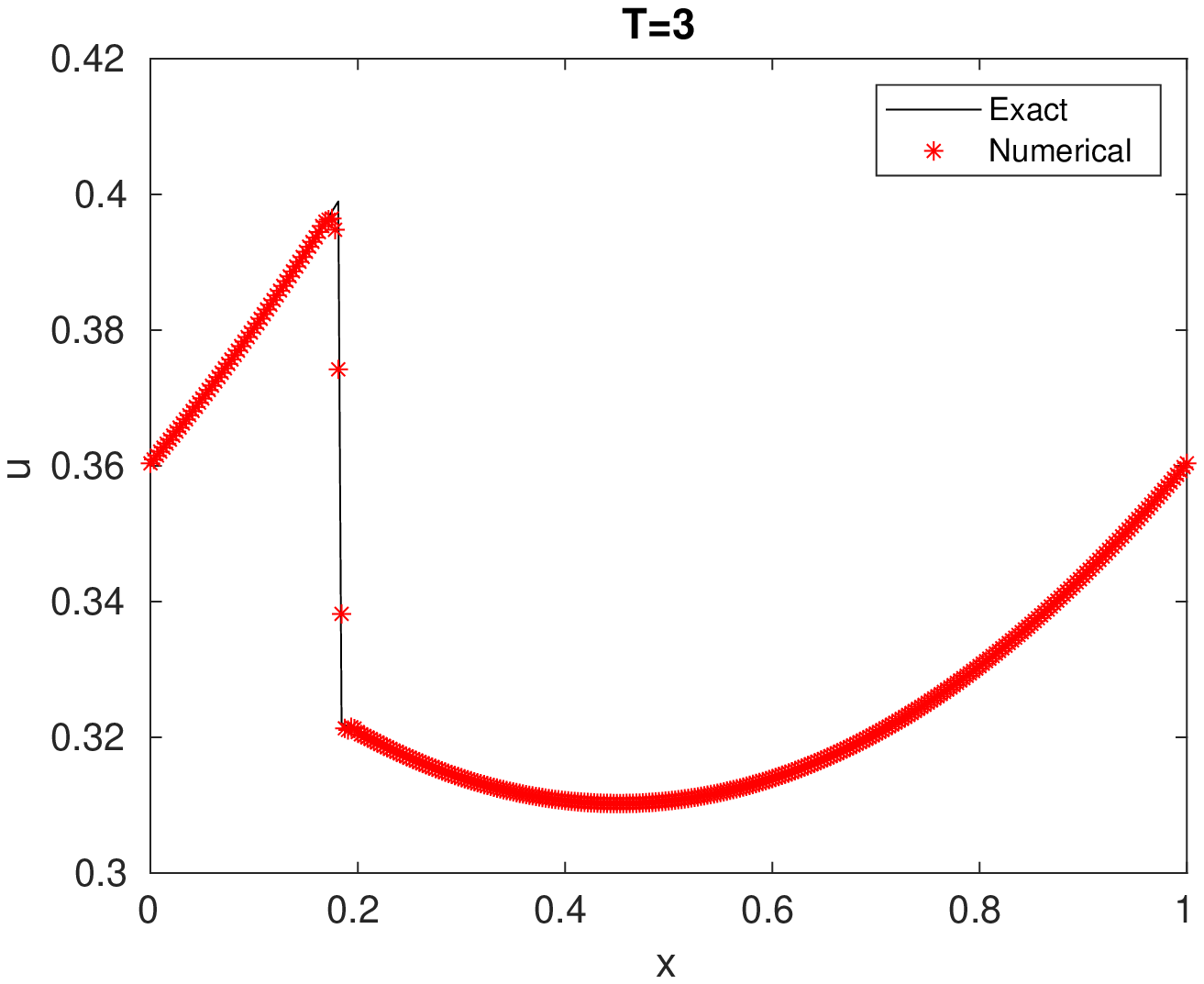} 
		\quad 
		\includegraphics[width=0.45\textwidth]{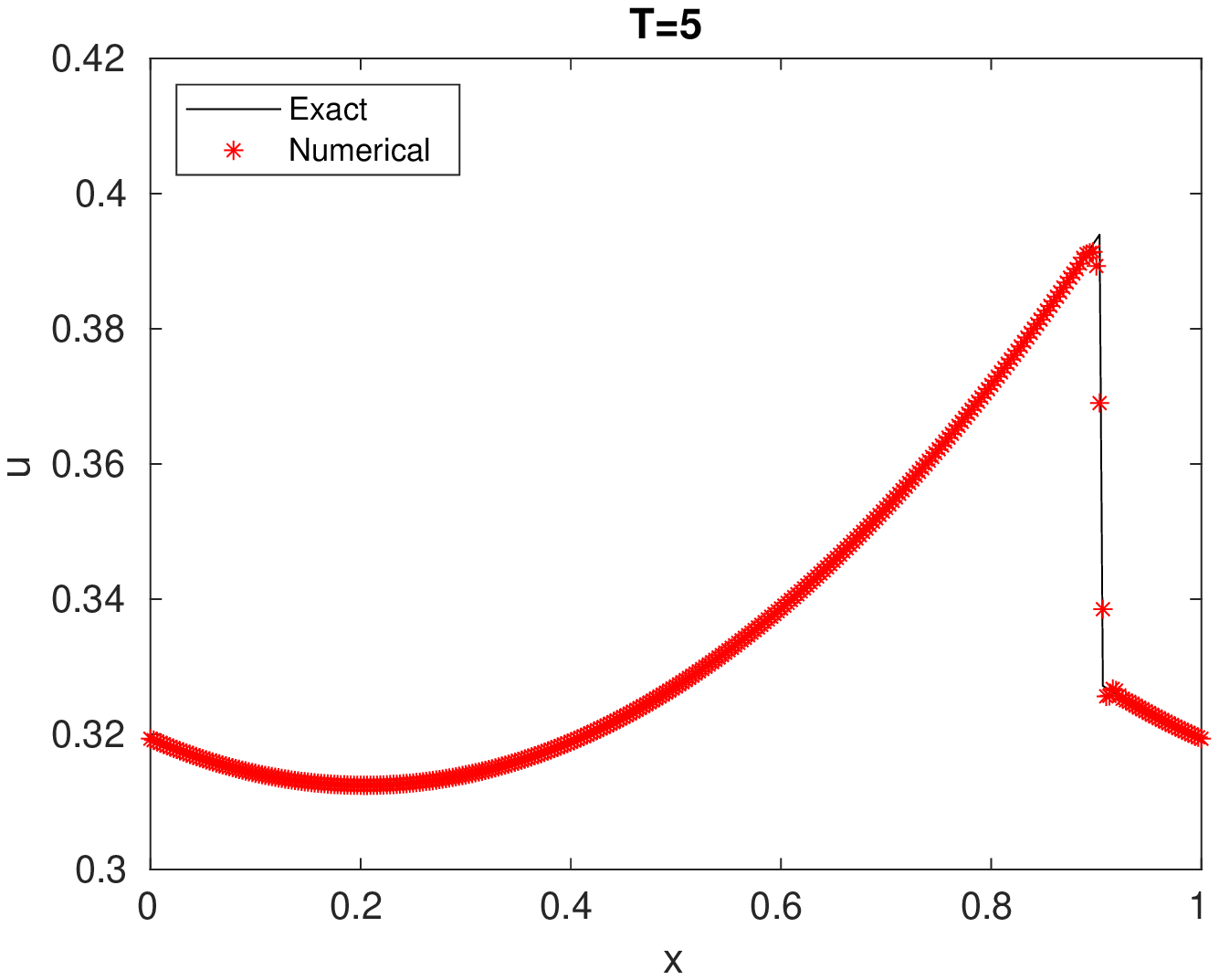}}
	
	\caption{One-shock solution of the $\mu$DP equation in Example \ref{Ex:shock_mudp}. $N = 320$. WENO5.}
	\label{Fig:mu-DP_weno5-zq_s1}
\end{figure} 

In Figure \ref{Fig:mu-DP_weno5-zq_s1}, \ref{Fig:mu-DP_mrweno5_s1} and \ref{Fig:mu-DP_mrweno7_s1}, we show the single shock solutions of the \mdp equation obtained by WENO5, MRWENO-5 and MR-WENO7, respectively. 
In Figure \ref{Fig:mu-DP_weno5-zq_s2}, \ref{Fig:mu-DP_mrweno5_s2} and \ref{Fig:mu-DP_mrweno7_s2}, we show two-shock solutions of the \mdp equation obtained by WENO5, MRWENO-5 and MR-WENO7, respectively. 
We observe sharp, non-oscillatory shock transitions for all three schemes. Our results agree well with those in
\cite{zhang_local_2019, zhao_high_2020}.
%
%
%
%

\begin{figure}[!htbp]
\centering%
\mbox{
\includegraphics[width=0.45\textwidth]{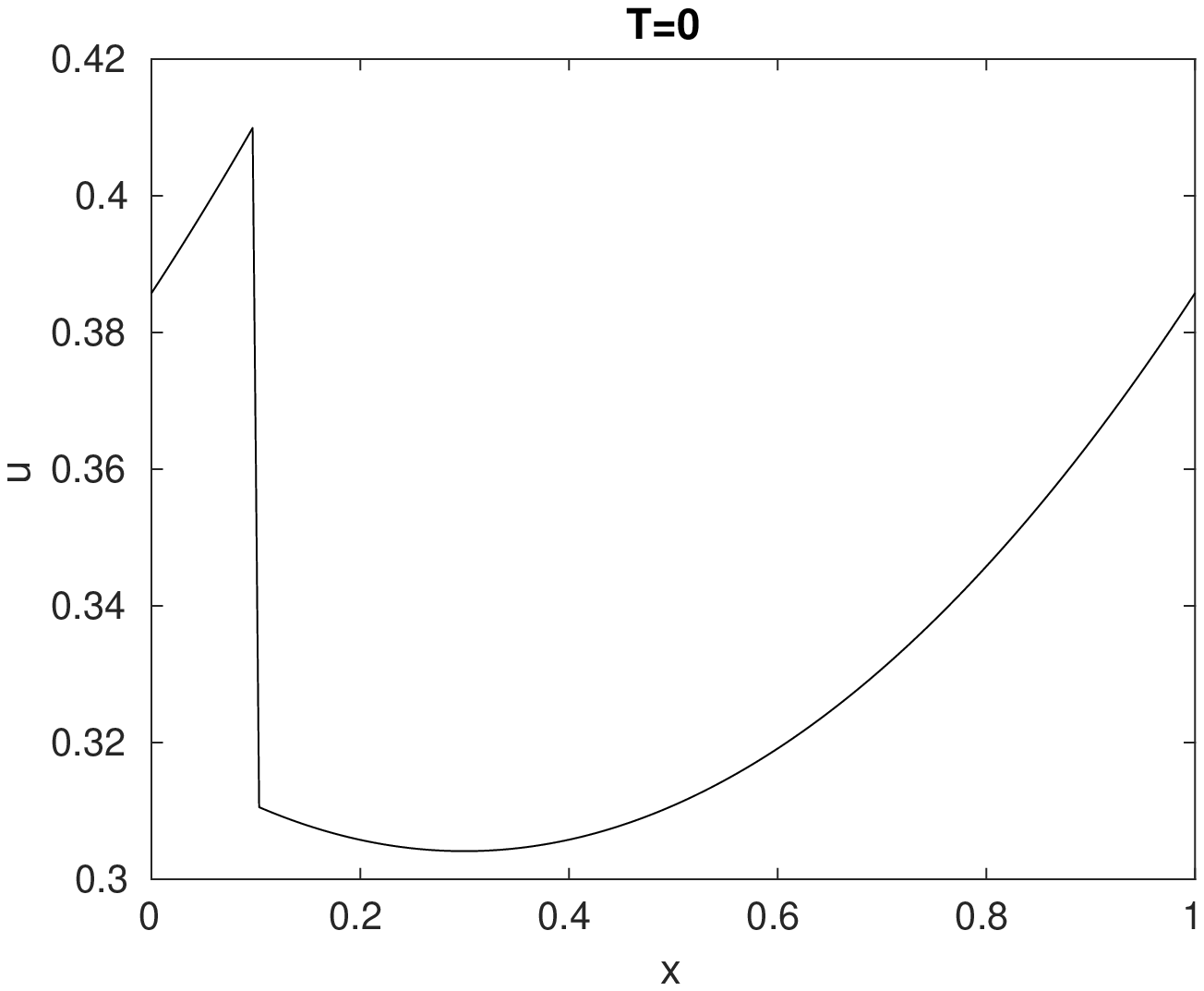} 
\quad 
\includegraphics[width=0.45\textwidth]{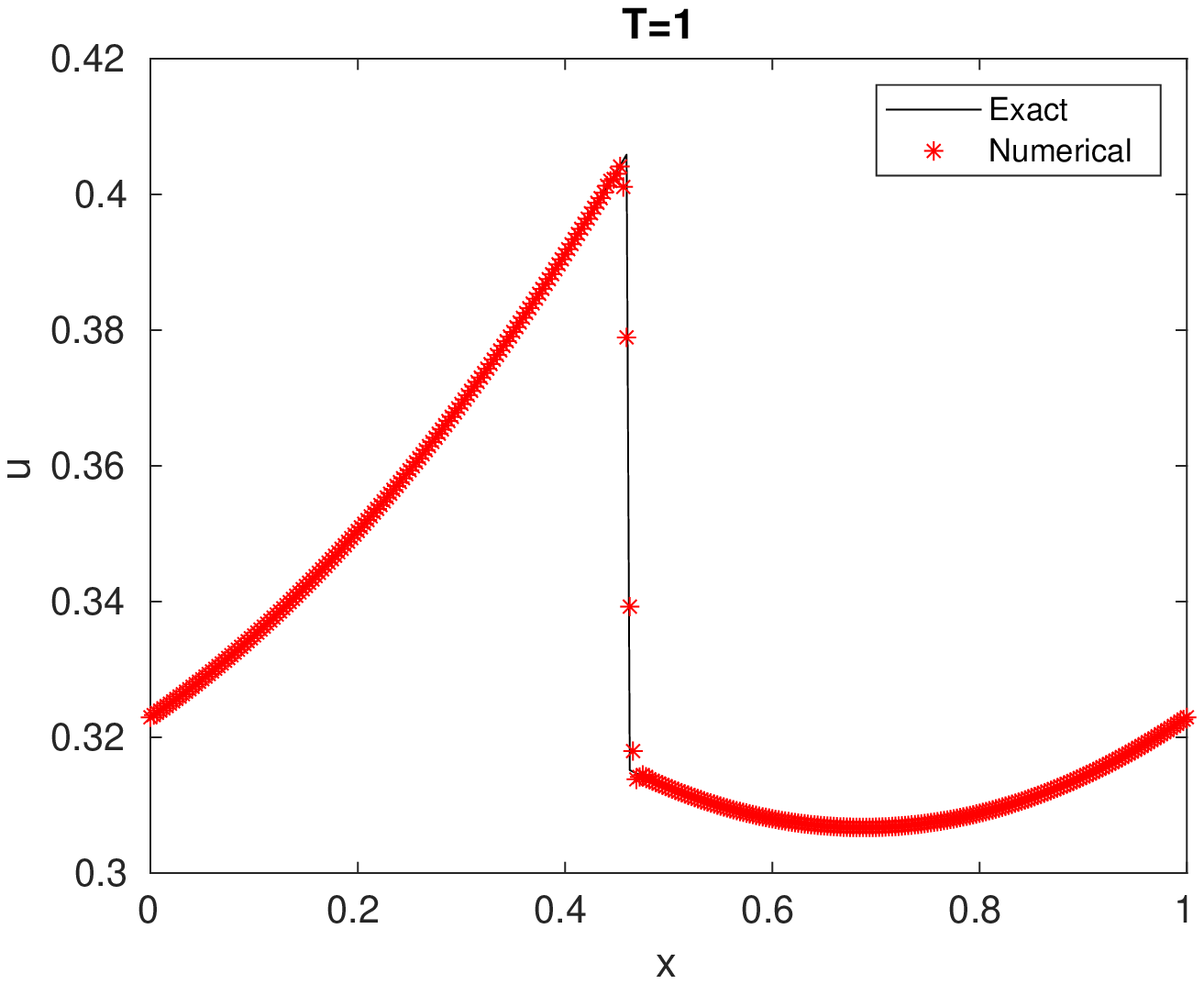}}
\quad
\mbox{
\includegraphics[width=0.45\textwidth]{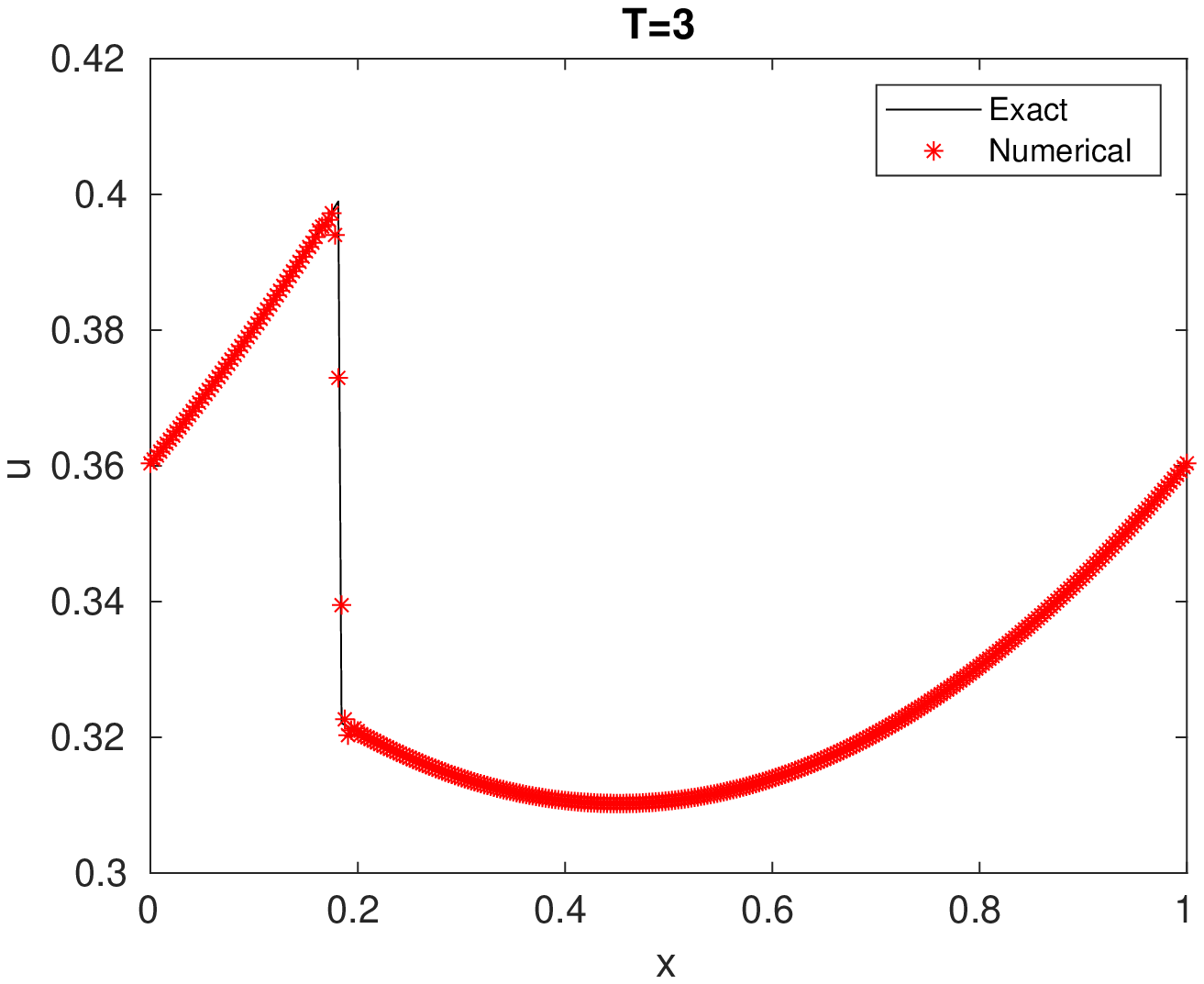} 
\quad 
\includegraphics[width=0.45\textwidth]{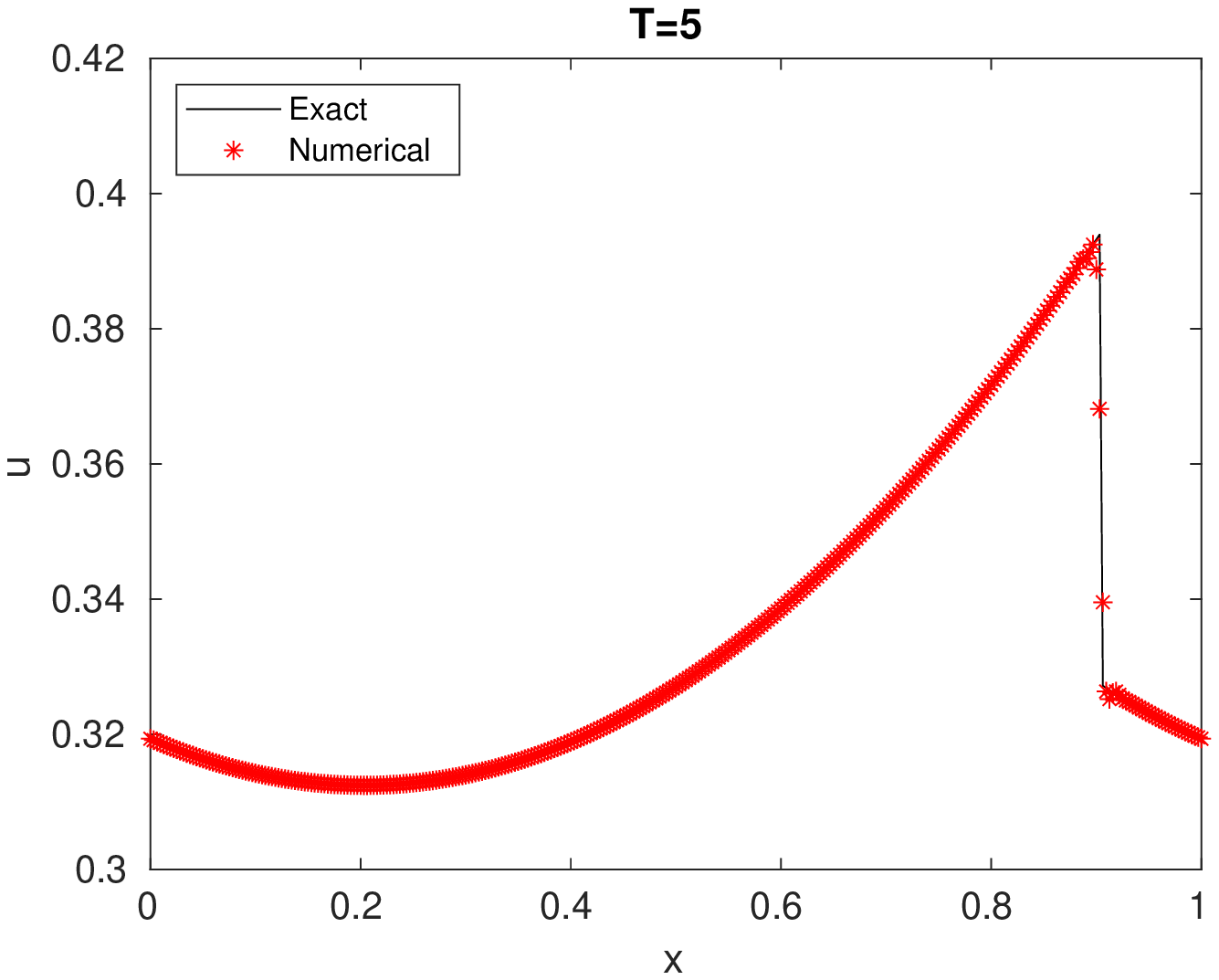}}

\caption{One-shock solution of the $\mu$DP equation in Example \ref{Ex:shock_mudp}. $N = 320$. MR-WENO5.}
\label{Fig:mu-DP_mrweno5_s1}
\end{figure}

\begin{figure}[!htbp]
\centering%
\mbox{
\includegraphics[width=0.45\textwidth]{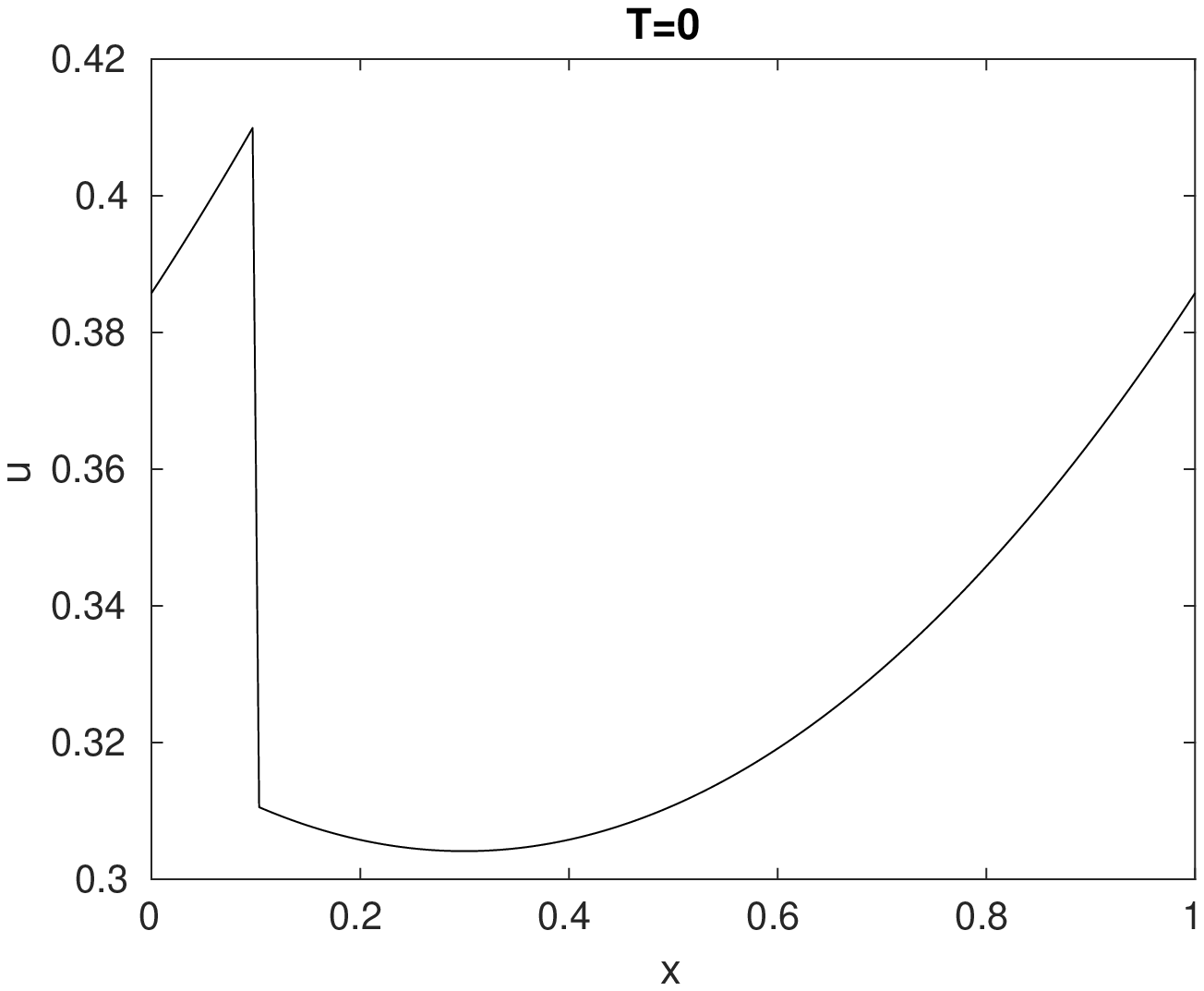} 
\quad 
\includegraphics[width=0.45\textwidth]{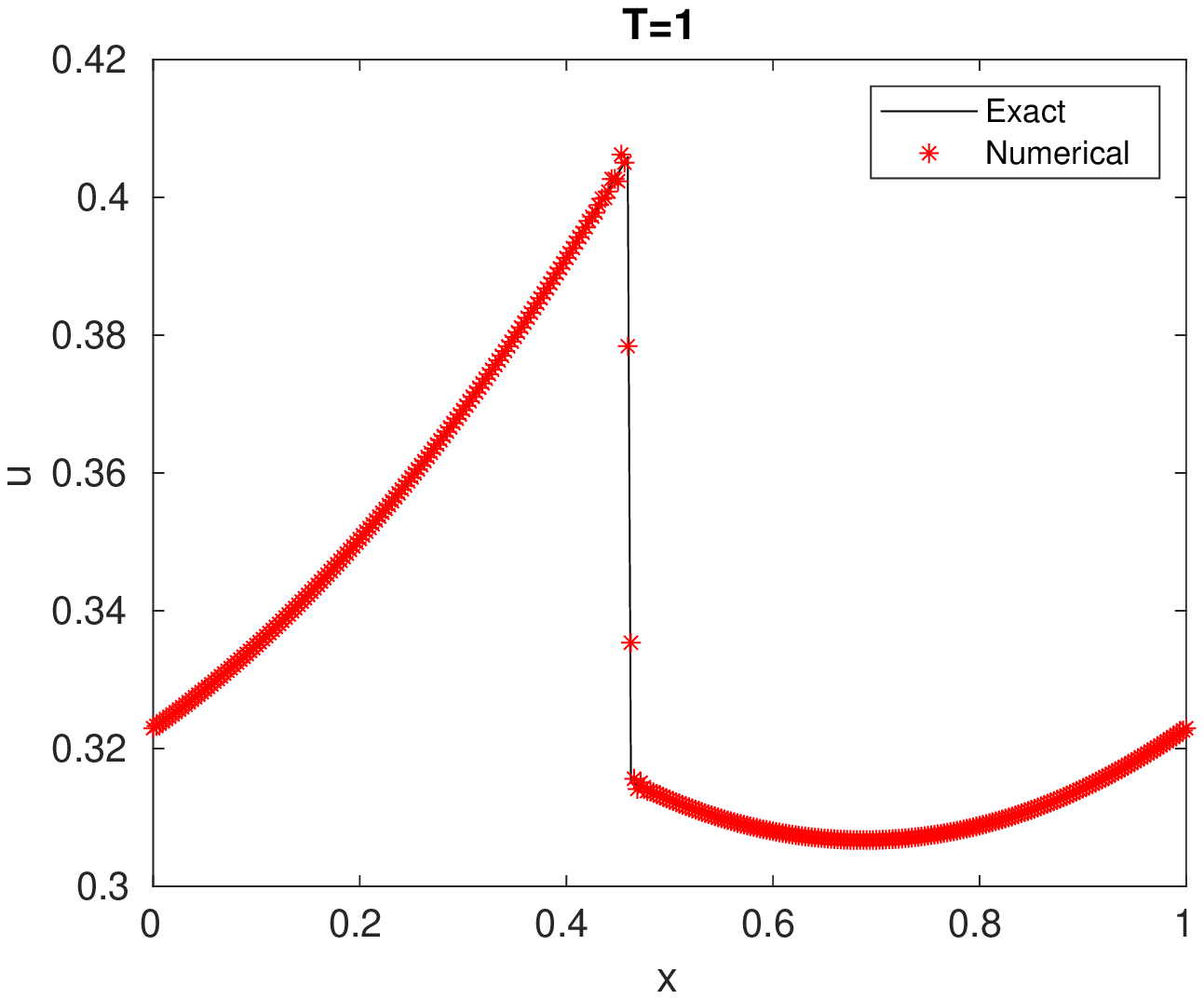}}
\quad
\mbox{
\includegraphics[width=0.45\textwidth]{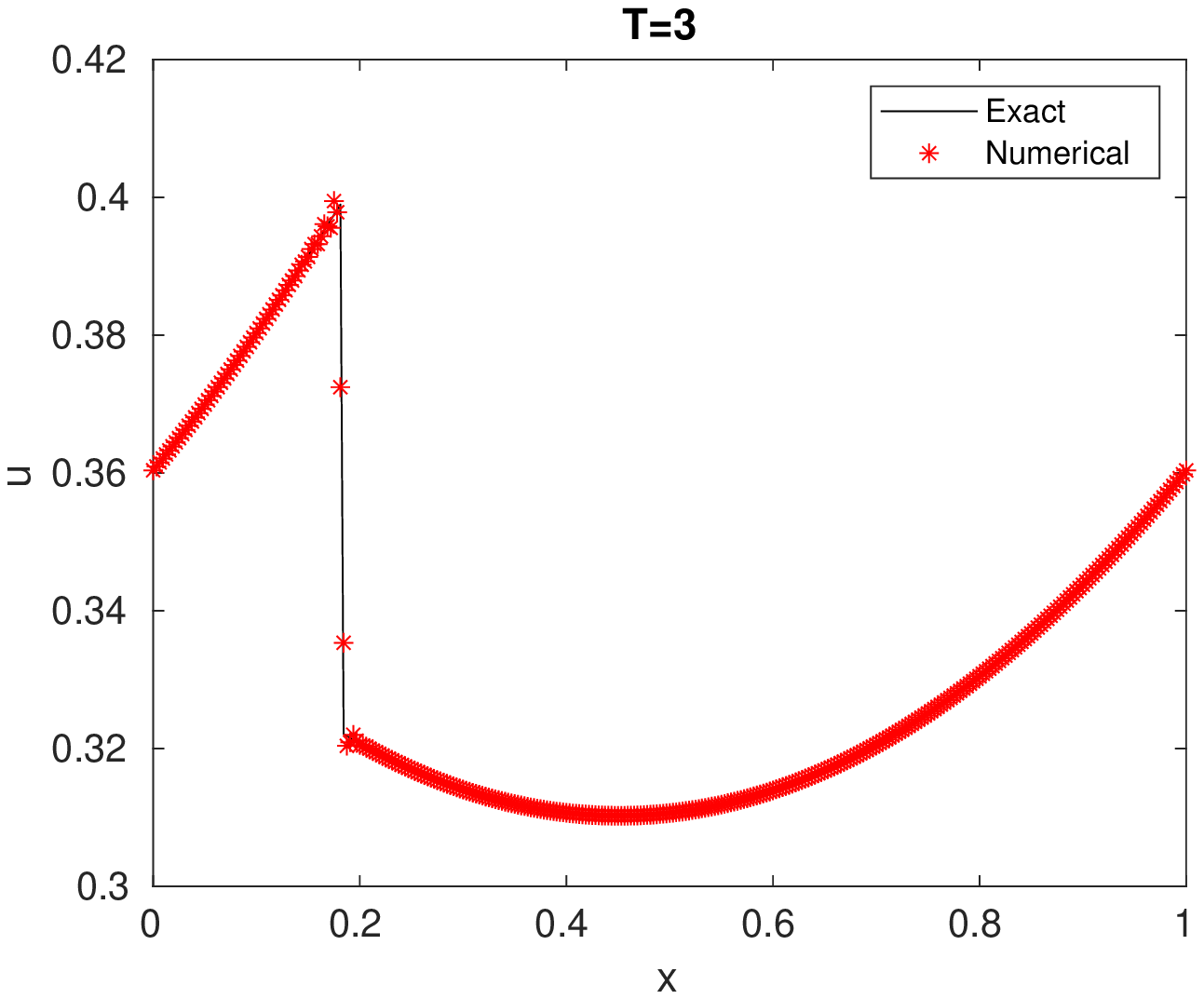} 
\quad 
\includegraphics[width=0.45\textwidth]{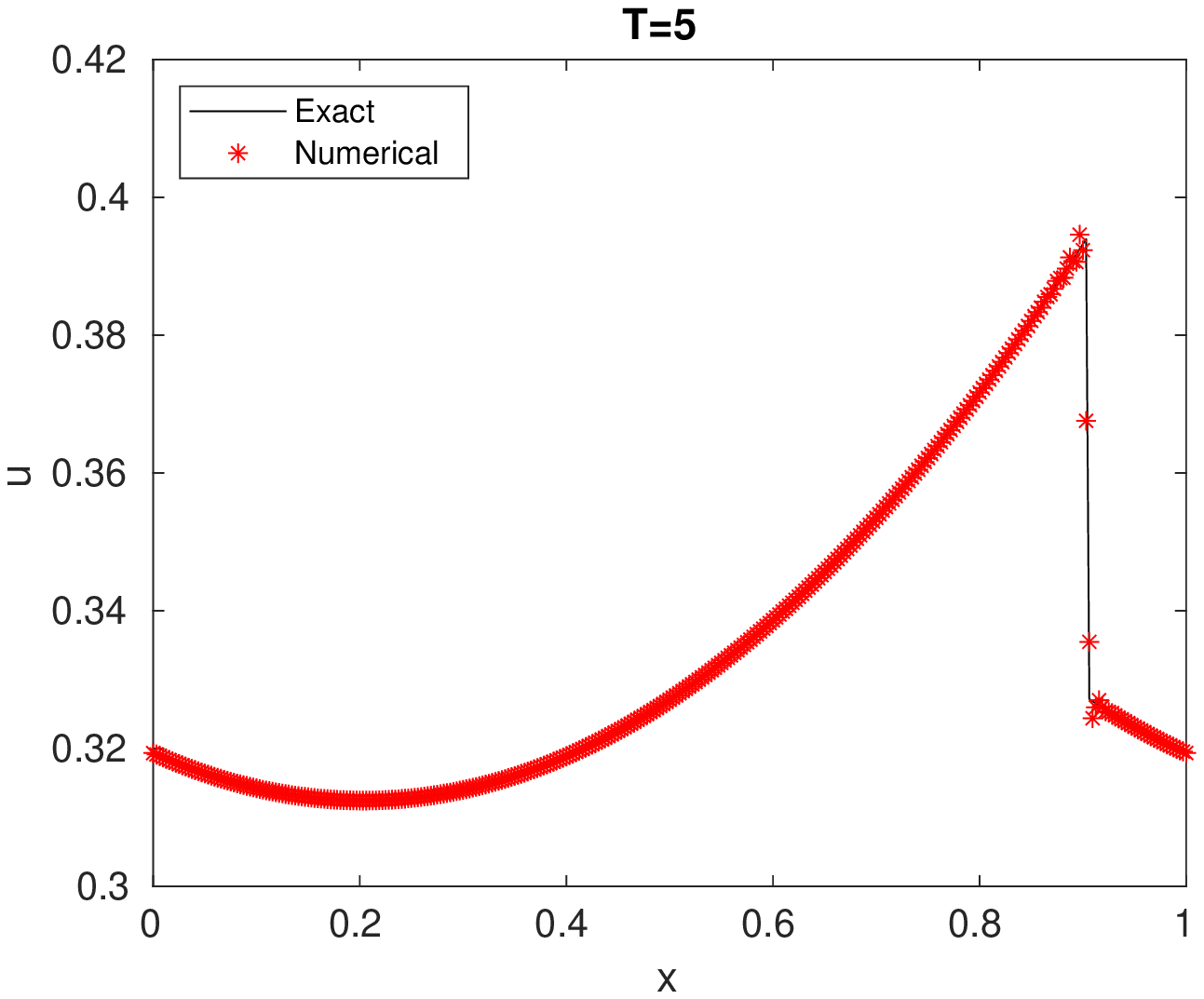}}

\caption{One-shock solution of the $\mu$DP equation in Example \ref{Ex:shock_mudp}. $N = 320$. MR-WENO7.}
\label{Fig:mu-DP_mrweno7_s1}
\end{figure}

\begin{figure}[!htbp]
\centering%
\mbox{
\includegraphics[width=0.45\textwidth]{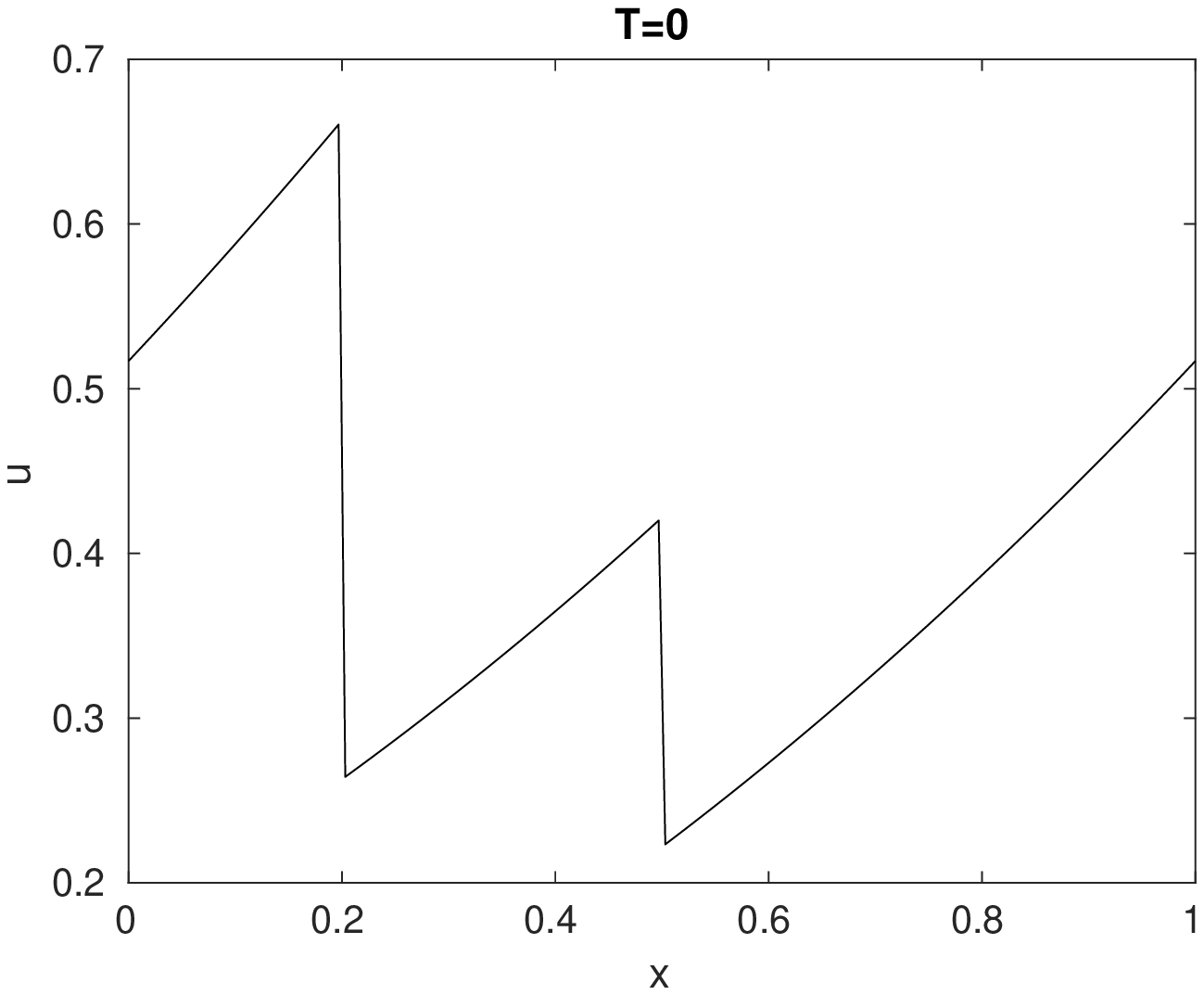} 
\quad 
\includegraphics[width=0.45\textwidth]{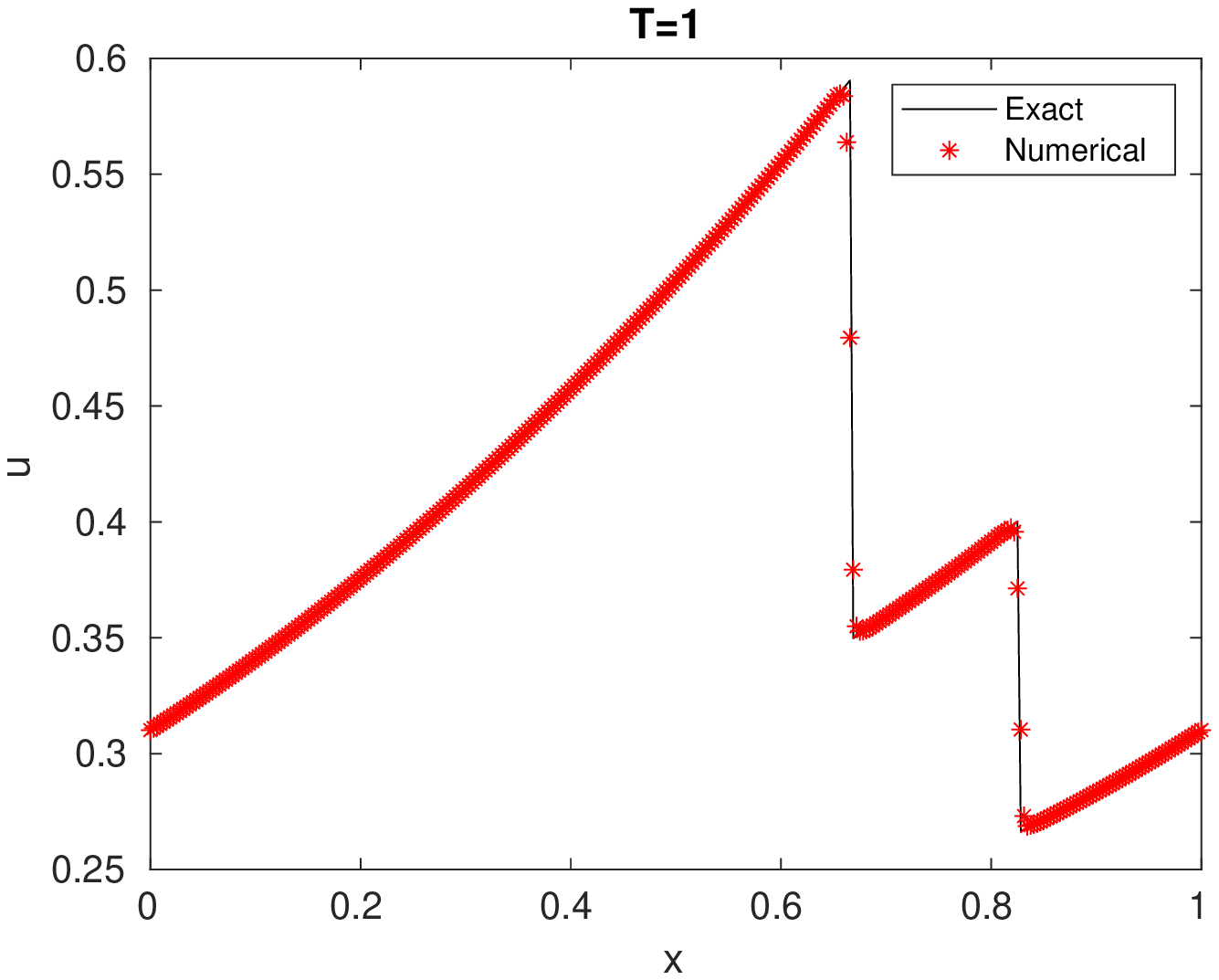}}
\quad
\mbox{
\includegraphics[width=0.45\textwidth]{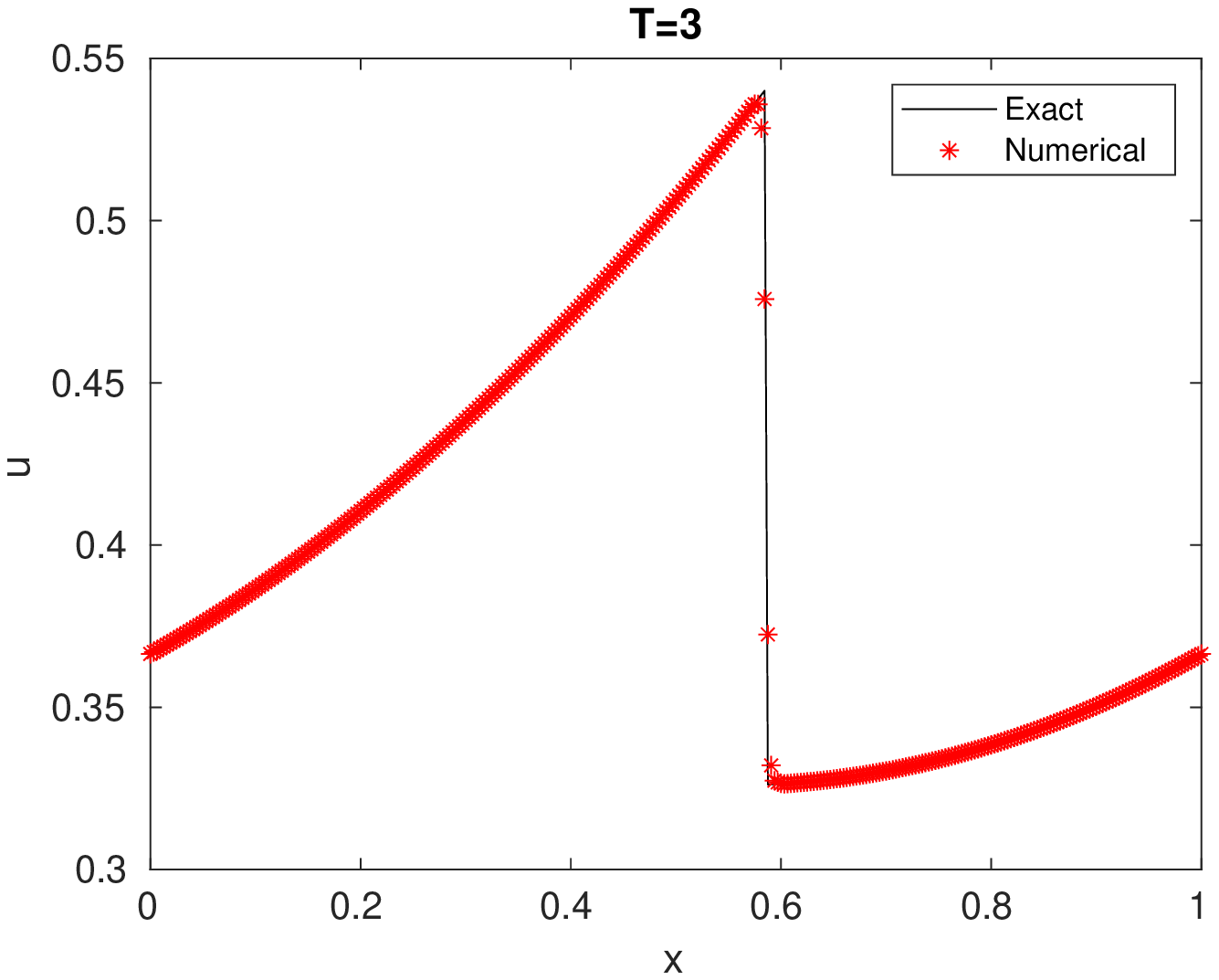} 
\quad 
\includegraphics[width=0.45\textwidth]{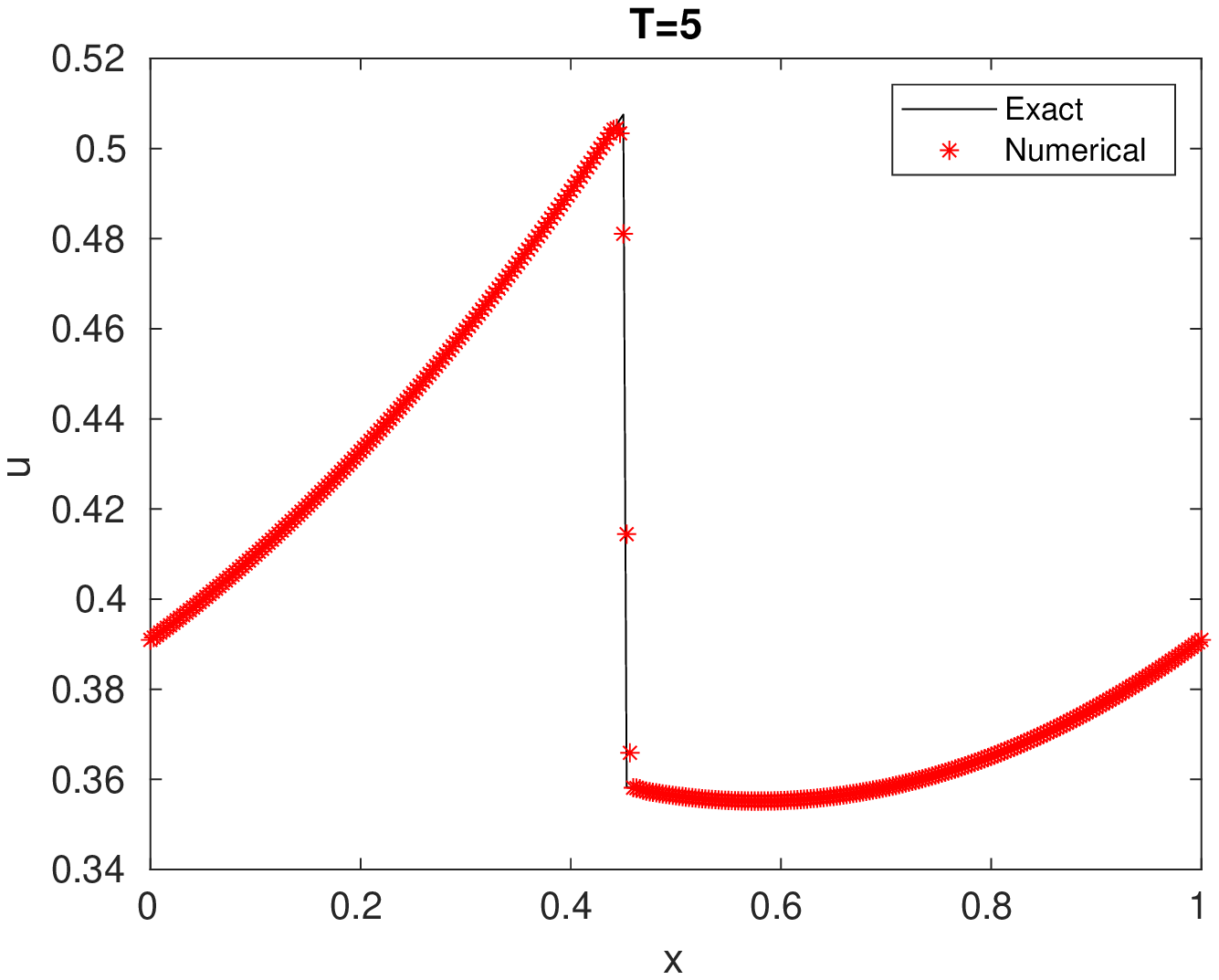}}

\caption{Two-shock solution of the $\mu$DP equation in Example \ref{Ex:shock_mudp}. $N = 320$. WENO5.}
\label{Fig:mu-DP_weno5-zq_s2}
\end{figure}

\begin{figure}[!htbp]
\centering
\mbox{
\includegraphics[width=0.45\textwidth]{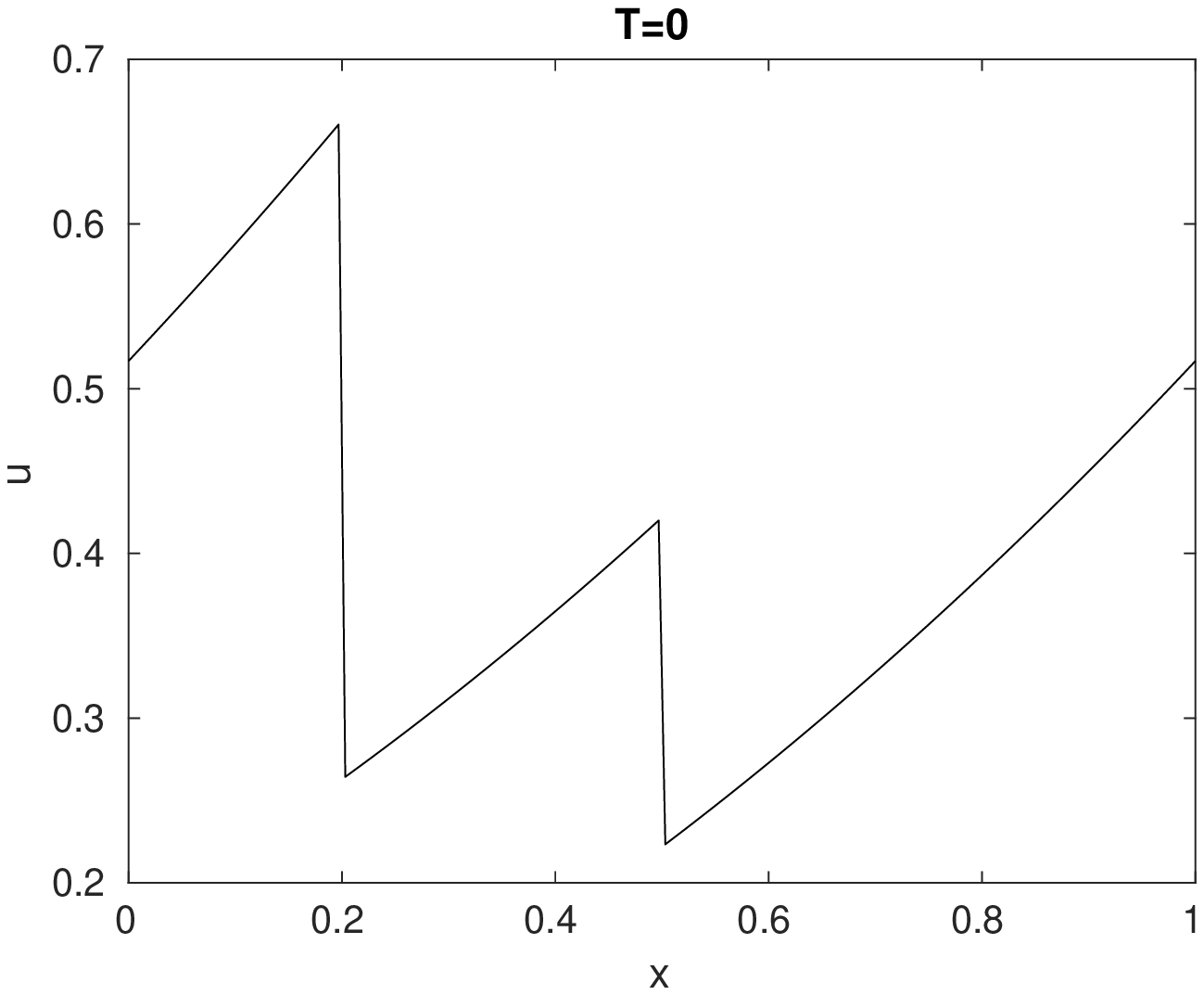} 
\quad 
\includegraphics[width=0.45\textwidth]{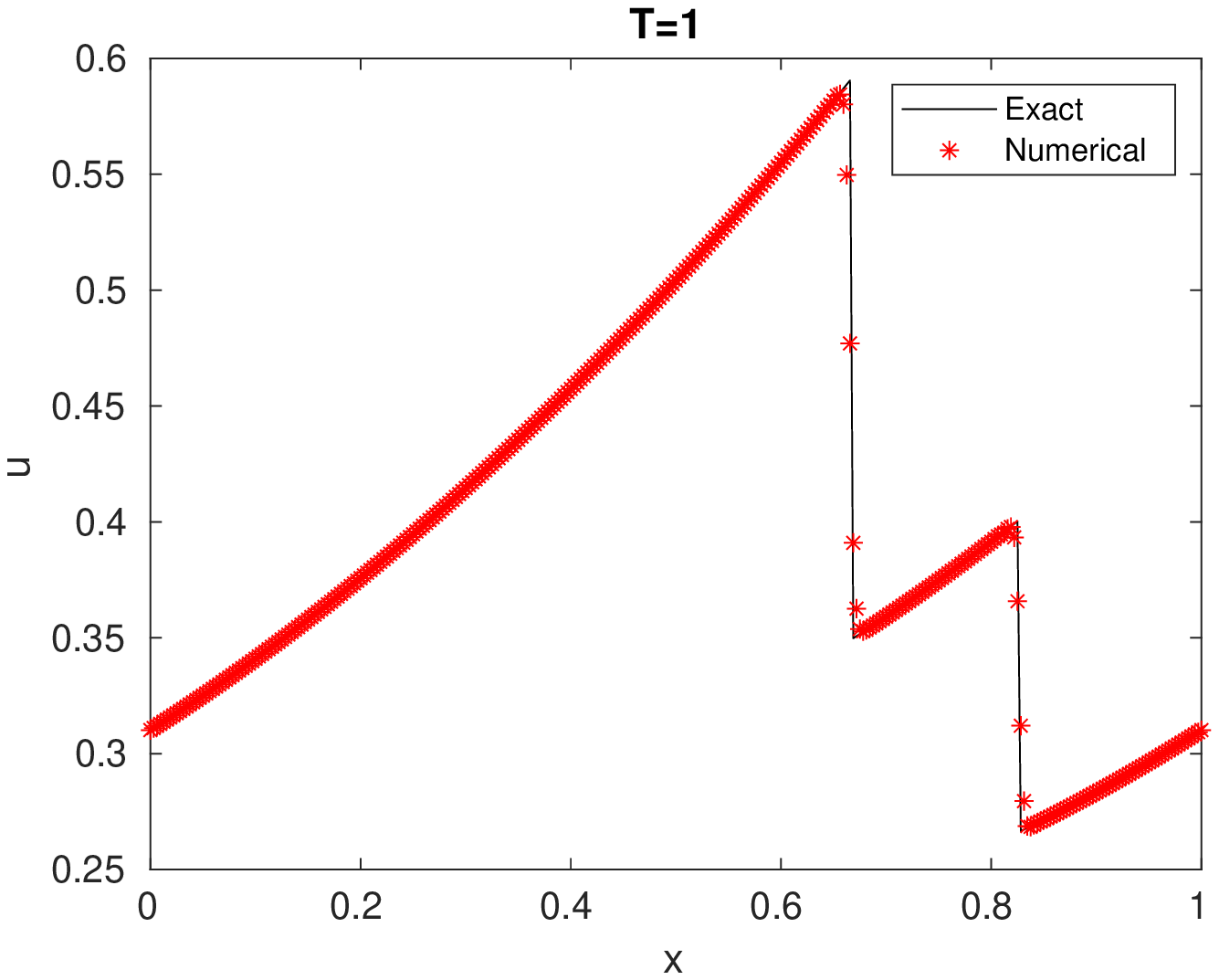}}
\quad
\mbox{
\includegraphics[width=0.45\textwidth]{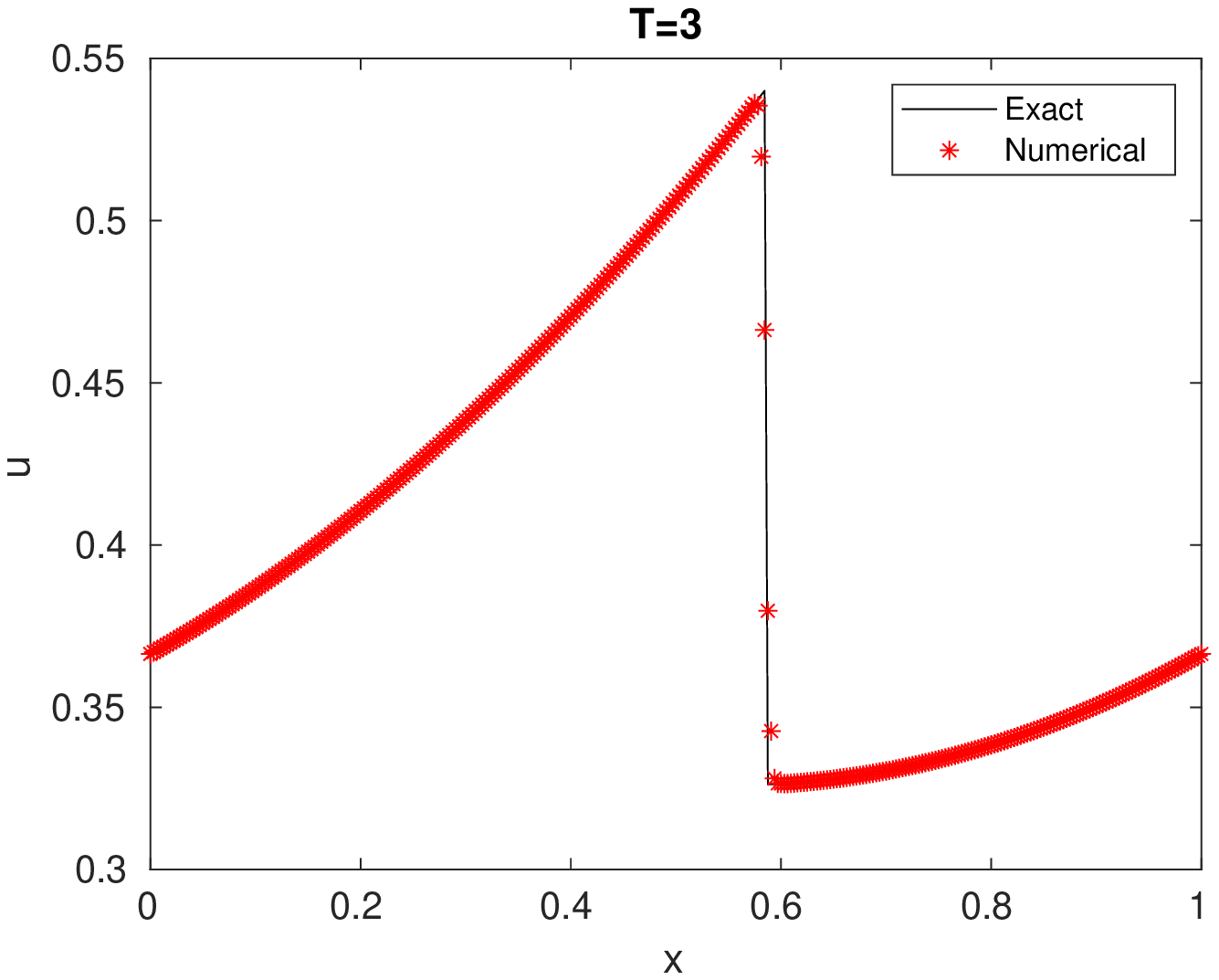} 
\quad 
\includegraphics[width=0.45\textwidth]{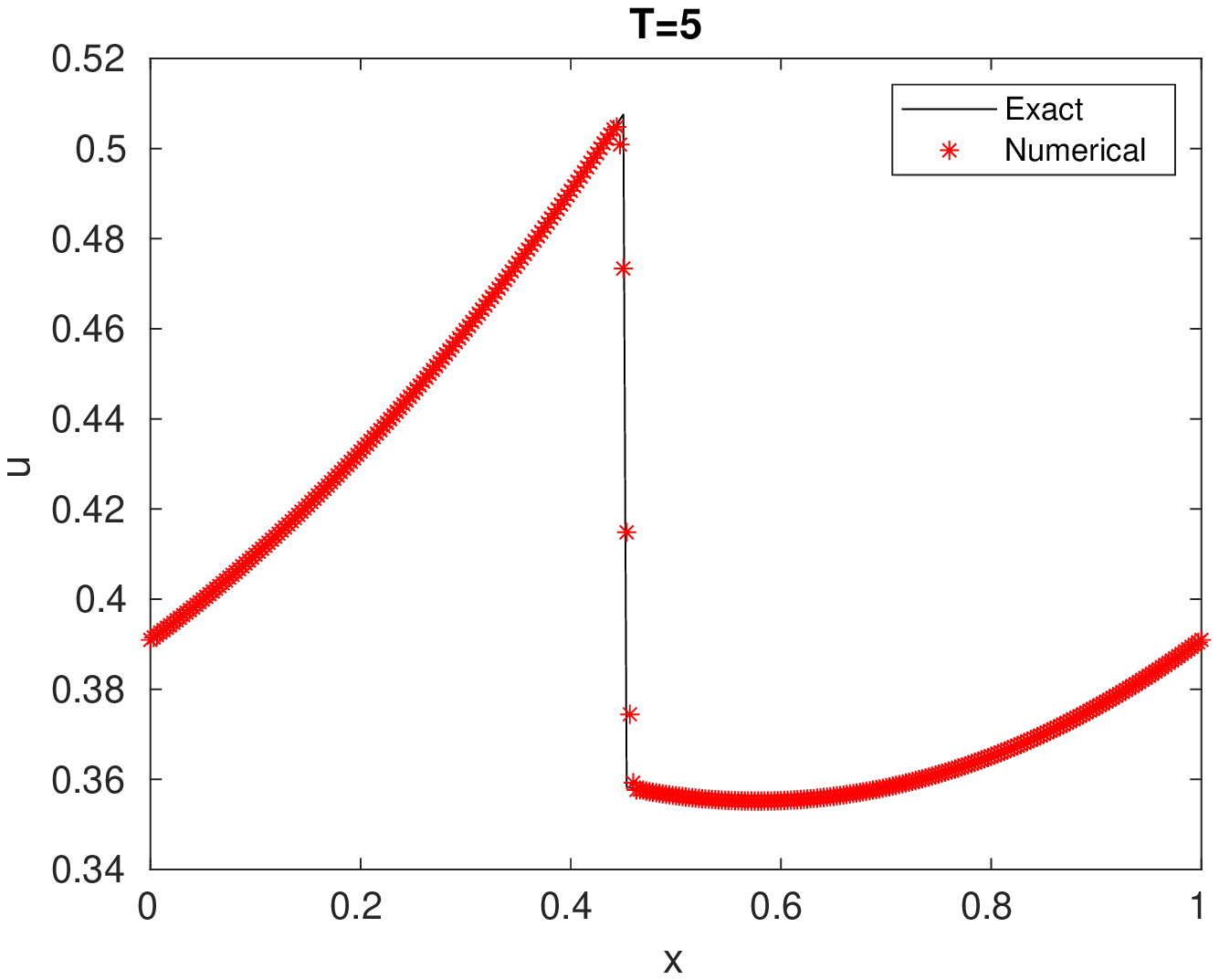}}

\caption{Two-shock solution of the $\mu$DP equation in Example \ref{Ex:shock_mudp}. $N = 320$. MR-WENO5.}
\label{Fig:mu-DP_mrweno5_s2}
\end{figure}

\begin{figure}[!htbp]
\centering%
\mbox{
\includegraphics[width=0.45\textwidth]{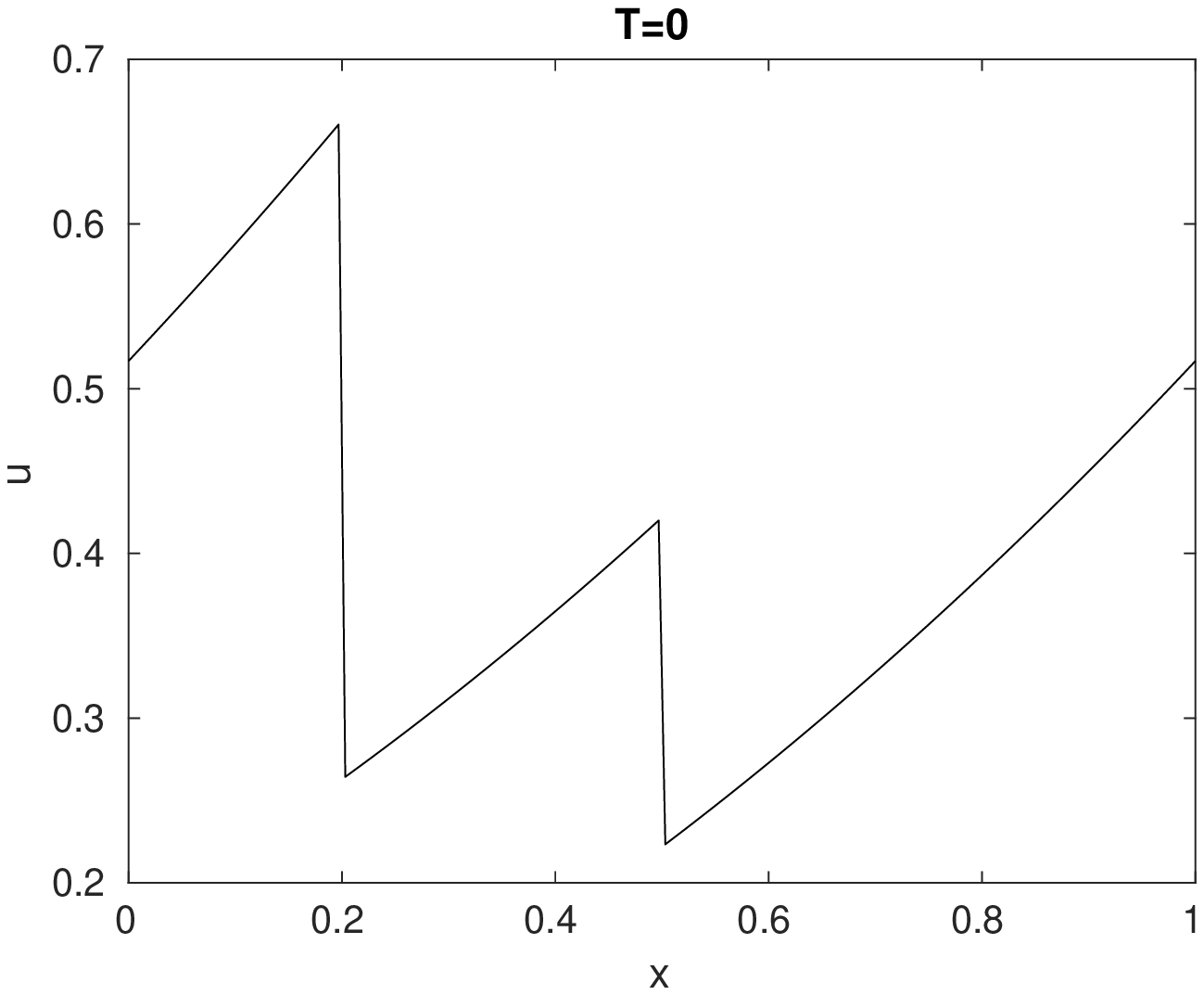} 
\quad 
\includegraphics[width=0.45\textwidth]{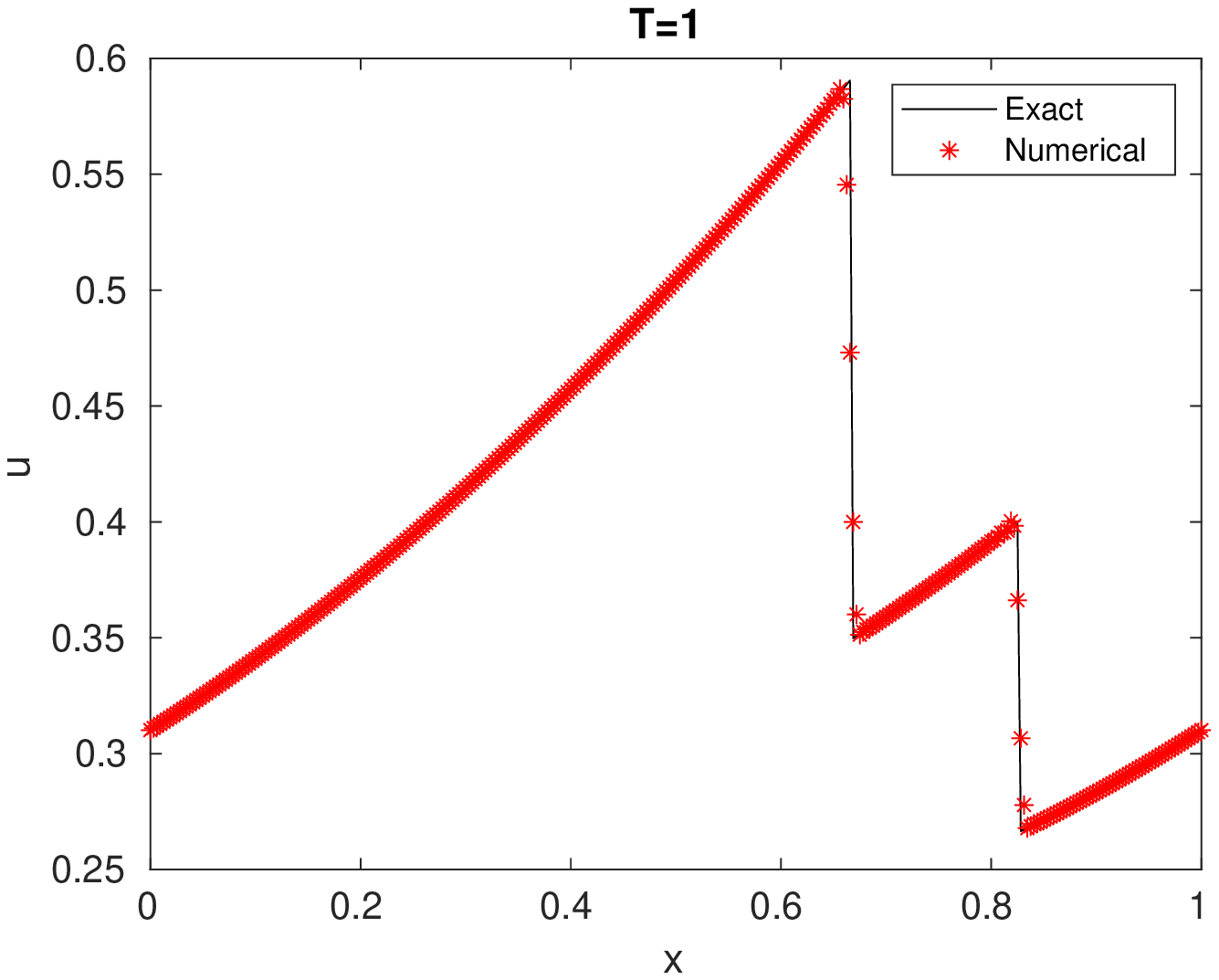}}
\quad
\mbox{
\includegraphics[width=0.45\textwidth]{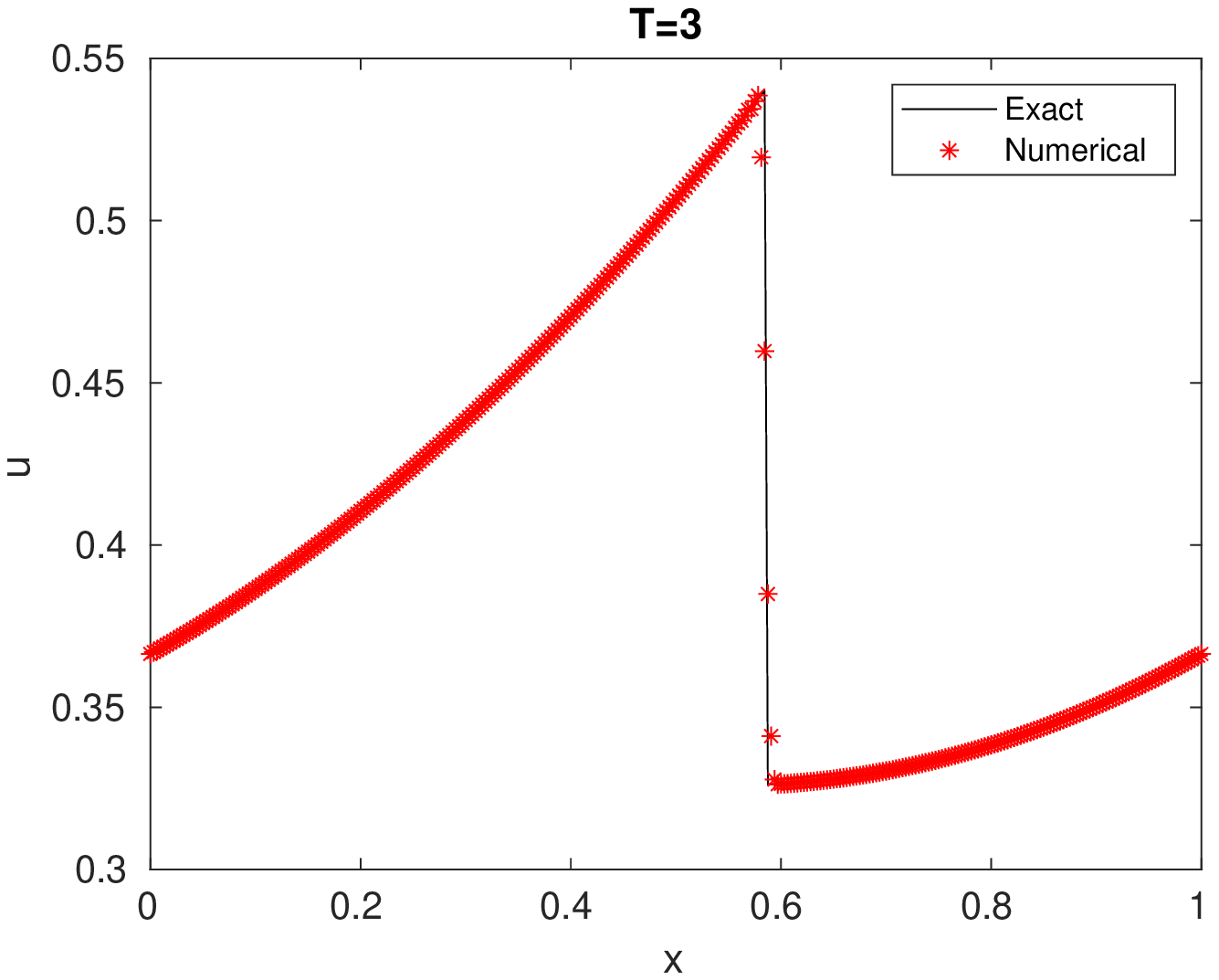} 
\quad 
\includegraphics[width=0.45\textwidth]{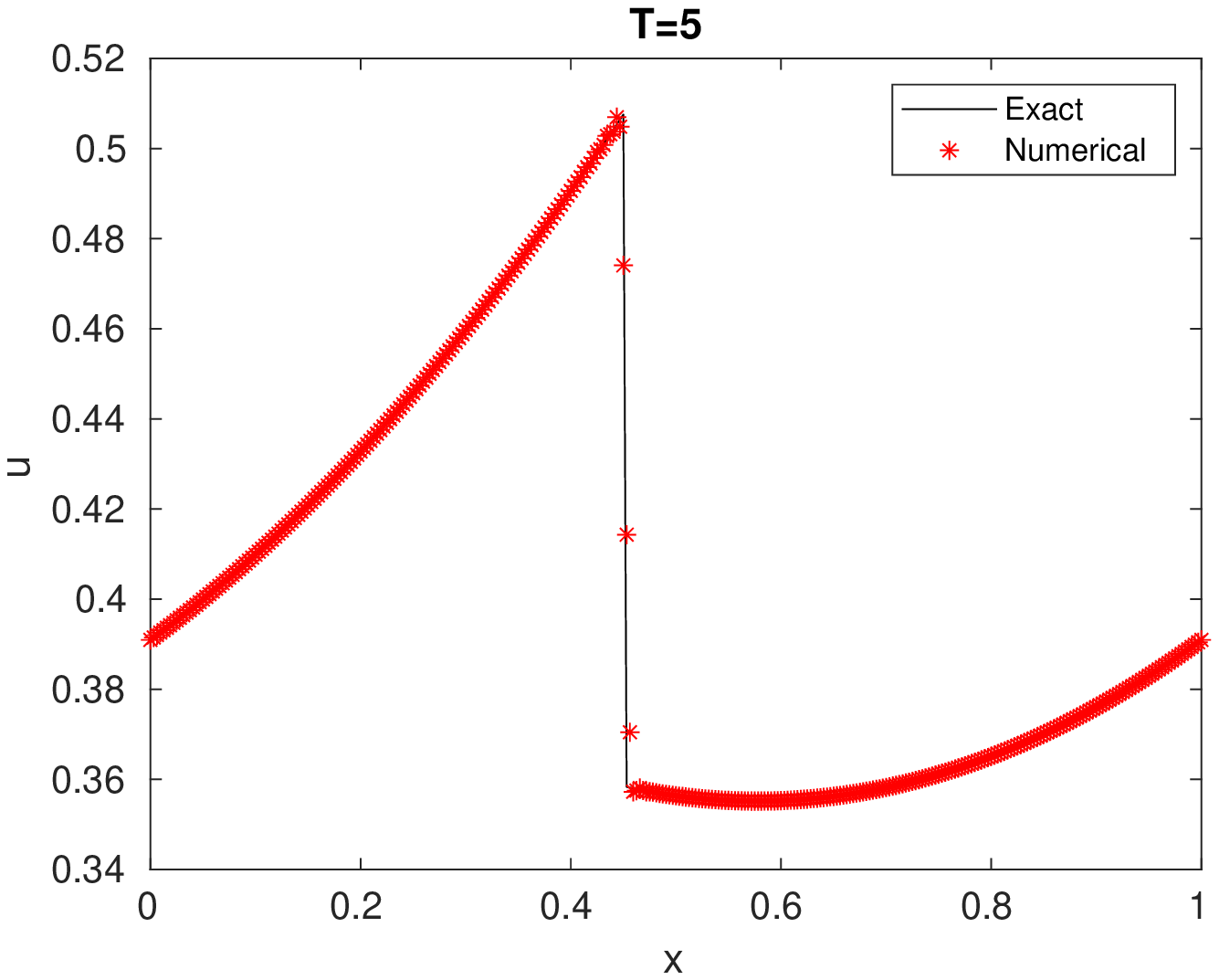}}

\caption{Two-shock solution of the $\mu$DP equation in Example \ref{Ex:shock_mudp}.   $N = 320$.  MR-WENO7.}
\label{Fig:mu-DP_mrweno7_s2}
\end{figure} 
\end{exa}

{\color{red}To conclude this section, we present the CPU time solving the DP and $\mu$-DP equations 
	in Table \ref{cpu}  with WENO5, MR-WENO5 and MR-WENO7 discussed  in this paper, as well as the finite difference WENO scheme discussed in \cite{xia_weighted_2017}, denoted as WENO5-JS,   based on the classical WENO scheme proposed by Jiang and Shu in \cite{jiang_efficient_1996}. The CPU time for each case is recorded  as the average of 5 runs on a ThinkCentre computer with an Intel core i7-6700H 3.40 GHz  and 16 GB RAM. 
	For simplicity, 
	Example \ref{Ex:peakon_one}-1 and \ref{Ex:peakon_one}-2 denote the single peakon and anti-peakon cases of Example \ref{Ex:peakon_one}, respectively. Example \ref{Ex:twopeakon}-1 and \ref{Ex:twopeakon}-2 denote the two-peakon interaction and two-anti-peakon interaction cases of Example \ref{Ex:twopeakon}, respectively. Example \ref{Ex: wavebreak}-1 and \ref{Ex: wavebreak}-2 denote the wave breaking cases of Example \ref{Ex: wavebreak} with initial conditions \eqref{init:wavebreak1} and \eqref{init:wavebreak2},  respectively. Example \ref{Ex:mudp_peakon}-1 and \ref{Ex:mudp_peakon}-2 denote one-peakon and two-peakon cases of  Example \ref{Ex:mudp_peakon}, respectively.  Example \ref{Ex:shock_mudp}-1 and \ref{Ex:shock_mudp}-2 denote one-shock and two-shock cases of Example \ref{Ex:shock_mudp}, respectively.}


\begin{table}[!htbp]
\centering
\caption{{\color{red}CPU time (seconds) of WENO5, MR-WENO5, MR-WENO7 and WENO5-JS on  a ThinkCentre computer with an Intel core i7-6700H 3.40 GHz  and 16 GB RAM. }}
\label{cpu}
\smallskip
\begin{tabular}{lrrcccc}
\toprule
\textbf{} &$N$&$T$&{WENO5}&{MR-WENO5}&{MR-WENO7}&{WENO5-JS}
\\ 
\midrule
{\em DP equation} &&&&&&   \\
\cmidrule{1-1}
Example \ref{Ex:peakon_one}-1 & 640 & 16 & 0.20 & 0.27 & 0.39 & 0.15\\
Example \ref{Ex:peakon_one}-2 & 640 & 16 & 0.20 & 0.26 & 0.39 & 0.17\\
Example \ref{Ex:twopeakon}-1 & 1280 & 12 & 0.57 & 0.81 & 1.16 & 0.46\\
Example \ref{Ex:twopeakon}-2 & 1280 & 12 & 0.57 & 0.81 & 1.14 & 0.47\\
Example \ref{Ex:shockpeakon} & 640 & 6 & 0.14 & 0.18 & 0.24 & 0.12\\
Example \ref{Ex:peakon+anti-peakon} & 640 &7& 0.19 & 0.24 & 0.35 & 0.15\\
Example \ref{Ex:triple} & 640 & 7 & 0.19 & 0.25 & 0.34 & 0.14\\
Example \ref{Ex: wavebreak}-1  & 640 & 1.1& 0.28 & 0.38 & 0.56 & 0.23\\
Example \ref{Ex: wavebreak}-2 & 2560 & 30 & 3.00 & 4.00 & 5.56 & 2.48\\
\midrule
{\em $\mu$DP equation} &&&&&& \\	
\cmidrule{1-1}
Example \ref{Ex:mudp_peakon}-1 & 160 & 1 & 1.75 & 1.78 & 1.83 & 1.73\\
Example \ref{Ex:mudp_peakon}-2 & 160 & 1 & 1.76 & 1.78 & 1.83 & 1.75\\
Example \ref{Ex:shock_mudp}-1  & 320 & 1 & 23.77 & 23.99 & 24.49 & 23.60\\
Example \ref{Ex:shock_mudp}-2  & 320 & 1 & 23.85 & 23.86 & 24.53 & 23.75\\
\bottomrule
\end{tabular}
\end{table}


\section{Conclusion}
\label{sec:conclusion}
In this paper, we investigate two finite difference WENO schemes with unequal-sized sub-stencils for solving the DP and \mdp equations. 
We first rewrite the DP equation as a hyperbolic-elliptic system and the \mdp equation as a first order system, by introducing auxiliary variable(s). 
Then suitable numerical fluxes are chosen to ensure stability and correct upwinding. 
For the numerical fluxes of the auxiliary variable(s), we choose a linear finite difference scheme to approximate them with suitable order of accuracy. For the numerical fluxes of the primal variable, we adopt two finite WENO procedures with unequal-sized sub-stencils for the reconstruction, i.e. the simple finite difference WENO procedure \cite{zhu_new_2016, zhu_new_2018} or the multi-resolution WENO procedure \cite{zhu_new_2018,zhu_new_2019,zhu_new_2020}.
The simple WENO procedure uses one large stencil and several smaller stencils, while the multi-resolution WENO procedure uses a hierarchy of nested central stencils. 
in which all stencils are central and if the large stencil has 7 cells, then the following smaller stencils have 5, 3 and 1 cell(s), respectively. 
Comparing with the classical WENO procedure, both WENO procedures with unequal-sized choose linear weights to be any positive number on the condition that their sum is one. They provide a simpler way for WENO reconstruction.
Numerical examples are provided to demonstrate that our proposed schemes
can achieve high order accuracy in smooth regions, and resolve shocks or peakons sharply and in an essentially non-oscillatory fashion. 

%
%

\section*{Acknowledgements}

Y. Xu was partially supported by National Natural Science Foundation of China (Grant No. 12071455).
X. Zhong was partially supported by National Natural Science Foundation of China (Grant No. 11871428).
The authors appreciate Dr. Yinhua Xia and Dr. Jun Zhu for many helpful discussions.


\bibliographystyle{abbrv}
\bibliography{DP,xinghui}
\end{document}